
\magnification=\magstep1
\vsize=22truecm
\input amstex
\documentstyle{amsppt}
\leftheadtext{J. Jer\'onimo-Castro, E. Makai, Jr.}
\rightheadtext{Ball characterizations II}
\topmatter
\title 
\centerline{Ball characterizations in planes}
\centerline{and spaces of constant curvature, II}
\vskip.1cm
\centerline{\rm{This pdf-file is not identical with the printed
paper.}}
\endtitle
\author J. Jer\'onimo-Castro*, E. Makai, Jr.**\endauthor

\thanks *Research (partially) supported by CONACYT, SNI 38848
\newline
**Research (partially) supported by Hungarian National
Foundation for 
Scientific Research, grant nos. T046846, T043520, K68398,
K81146, Research supported by ERC Advanced Grant
\newline
``GeoScape'', No. 882971.
\endthanks
\subjclass {\it Mathematics Subject Classification} 2020. 
52A55
\endsubjclass
\abstract
High proved the following theorem. 
If the intersections of any two congruent
copies of a plane convex body are centrally symmetric, then
this body is a
circle. 
In our paper we extend the theorem of High to the sphere and
the hyperbolic plane, and partly to 
spaces of constant curvature.
We also investigate the dual question about the
convex hull of the
unions, rather than the intersections.

Let us have in $H^2$ proper
closed convex subsets $K,L$ with interior points,
such that 
the numbers of the connected
components of the boundaries of $K$ and $L$ are finite.
We exactly describe all 
pairs of such subsets $K,L$, whose any
congruent copies have an intersection with axial symmetry;
there are nine cases. (The cases of $S^2$ and ${\Bbb{R}}^2$
were described in Part I, i.e., \cite{5}.)

Let us have in $S^d$, ${\Bbb{R}}^d$ or $H^d$ proper
closed convex $C^2_+$ subsets $K,L$ with interior points,
such that all sufficiently small intersections of their
congruent copies are symmetric w.r.t.\ a particular
hyperplane. Then the boundary components of both $K$ and $L$
are congruent, and each of them is
a sphere, a parasphere or a hypersphere.

Let us have a pair of
convex bodies in $S^d$, ${\Bbb{R}}^d$ or $H^d$,
which have at any boundary points supporting
spheres (for $S^d$ of
radius less than $\pi /2$). If
the convex hull of the union of any congruent copies of
these bodies is
centrally symmetric, then our bodies are congruent balls
(for $S^d$
of radius less than $\pi /2$). An analogous statement holds
for symmetry w.r.t.\ a particular hyperplane.
For $d=2$, suppose the
existence of the above supporting circles
(for $S^2$ of radius less than $\pi /2$),
and, for $S^2$, smoothness of $K$ and $L$.
If we suppose axial symmetry of all
the above convex hulls, then our bodies are (incongruent)
circles
(for $S^2$ of radii less than $\pi /2$).
\endabstract
\endtopmatter\document

$^*$ 
Facultad de Ingenier\'\i a, Universidad Aut\'onoma de
Quer\'etaro, Centro
Uni\-ver\-si\-ta\-rio, 
\newline Cerro de las Campanas s/n C.P. 76010,
Santiago de Quer\'etaro, Qro. M\'exico, ME\-XI\-CO
\newline
{\rm{ORCID ID: https://orcid.org/0000-0002-6601-0004}}
\newline
$^{**}$ Alfr\'ed R\'enyi Mathematical Institute,
Hungarian Research Network (HUN-REN),
\newline
H-1364 Budapest, Pf. 127, HUNGARY
\newline
{\rm{http://www.renyi.hu/\~{}makai}}
\newline
{\rm{ORCID ID: https://orcid.org/0000-0002-1423-8613}}

\vskip.1cm

{\it{E-mail address:}}
$^*$ jeronimo\@cimat.mx, jesusjero\@hotmail.com
\newline
$^{**}$ makai.endre\@renyi.hu

\head 5. New results, Theorems 5--8 \endhead

$S^d$, ${\Bbb{R}}^d$, $H^d$, with $d \ge 2$,
are the $d$-dimensional spherical,
Euclidean and hyperbolic spaces, resp.
{\it{Convexity of a set $K \subset H^d$}} is defined
as for $K \subset {\Bbb{R}}^d$.
{\it{Convexity of $K \subset S^d$, with}} 
int\,$K \ne \emptyset $, is meant as follows:
for
any two non-antipodal points of $K$ the shorter great
circle arc connecting them belongs to $K$. Then for
two antipodal points of $K$ some great circle arc
connecting them belongs to $K$.
A {\it{convex body $Y \subset X$}} is a compact convex set
with nonempty interior. By a {\it{ball, or sphere in}} $S^d$,
we mean one {\it{of radius at most}} $\pi / 2$ (thus a ball is
convex).


In the following theorem, the base line of a straight line
is meant to be itself. We recall (from Part I, i.e.,
\cite{5}) 
$$
\cases
X {\text{ is }} S^d, \,\,{\Bbb R}^d {\text{ or }} H^d,
{\text{ with }} d \ge 2, 
{\text{ and }} K,L \subsetneqq X {\text{ are closed
convex sets}}
\\
{\text{with interior points. Moreover, }}
\varphi , \psi : X \to X, {\text{ sometimes with indices,}}
\\
{\text{are orientation preserving congruences, with }}
{\text{int}}\,[ (\varphi K) \cap (\psi L) ]
\ne \emptyset .
\endcases
\tag *
$$


\proclaim{Theorem 5}
Assume \thetag{*} with $d = 2$ and let $X=H^2$. 
Then we have $(3) \Longrightarrow (2) \Longrightarrow (1)$.
If
all connected components of the boundaries of both of $K$
and $L$ are 
hypercycles or straight lines, let 
their total number be finite. Then
we have $(3) \Longleftrightarrow (2) \Longleftrightarrow
(1)$. Here:
\roster
\item
For each $\varphi , \psi $
we have that $(\varphi K) \cap (\psi L)$
admits some non-trivial congruence.
\item
For each $\varphi , \psi $
we have that $(\varphi K) \cap (\psi L)$
is axially symmetric.
\item 
We have either {\rm{(A)}}, or {\rm{(B)}}, or {\rm{(C)}}, or
{\rm{(D)}} or {\rm{(E)}}, where 
\newline
{\rm{(A)}}: Any of \,$K$ and $L$ is a circle, 
a paracircle, a convex domain bounded by
a hypercycle, 
or a half-plane.
However, if one of \,$K$ and
$L$ is a convex set bounded by
a hypercycle or is a half-plane, 
then the other one is either a circle, or a congruent copy of
the first one.
\newline 
{\rm{(B)}}: One of \,$K$ and $L$ is a circle, and the other one is
bounded either by two hypercycles, whose base lines coincide, 
or by a hypercycle, and its base line.
\newline
{\rm{(C)}}: One of \,$K$ and $L$ is a circle, of radius $r$, say, and the other one is
bounded by at least two
hypercycles or straight lines (with all base lines different), 
whose mutual distances are at least
$2r$.
\newline
{\rm{(D)}}: One of \,$K$ and $L$ is a paracircle, and the
other one is a parallel domain
of some straight line, for some distance $l>0$.
\newline
{\rm{(E)}}: $K$ and $L$ are congruent, and both are parallel
domains of some straight lines, for some distance $l > 0$.
\endroster
\endproclaim


{\bf{Remark 2.}} Observe that in Part I
(i.e., \cite{5}) 
Theorem 1 dealt with 
intersections of sufficiently small diameters and
Theorems 3, 4, and this paper Theorem 5 dealt with all
intersections. There is a question ``in between'' these: the
case of all

\newpage

\noindent
compact intersections (observe that
in Theorem 2 the
intersections are compact).
The
case of central symmetry of all compact
intersections, in \cite{4} 
Theorem 3,
in $S^d$, ${\Bbb{R}}^d$ and $H^d$,
under some regularity hypotheses (weaker than $C^2_+$), 
was clarified. There were
found six cases. For central symmetry of all
intersections \cite{4} 
Theorem 2 gave that
in $S^d$, ${\Bbb{R}}^d$ and $H^d$, under the
above mentioned regularity hypotheses, $K$ and $L$
were congruent balls, of radius at most $\pi /2$.
This was just one of the above six
cases. 

Also here, for some non-trivial congruence, or axial symmetry
the ``in-between'' question about all compact intersections,
in $S^2$, ${\Bbb{R}}^2$ and $H^2$,
can be posed, supposing the finiteness hypothesis in
Theorem 5. Observe that for $X = S^2$ this is settled by
Theorem 2. For $X = {\Bbb{R}}^2$ the answer is Theorem 1, (3).
Namely, Theorem 3, (2) settles the case of all intersections
(there are five cases),
so these are examples also for all compact intersections.
On the other hand, suppose that all compact intersections, with
nonempty interior, are axially symmetric, or admit some
non-trivial congruence. Then also all intersections, of
sufficiently small diameter, are axially symmetric, or admit
some
non-trivial congruence. These are described in Theorem 1, (3)
(there are six cases).
There
is just one case in Theorem 1, (3), not covered by Theorem 3,
(2), namely the case of a half-plane and a parallel strip.
However, these
do not have a compact intersection with nonempty interior.
Hence in this case the above properties ''admitting some
non-trivial congruence'' and ``axial symmetry''
hold vacuously. 

\vskip.1cm


{\bf{Problem 3.}} 
Can one decribe in $H^2$
the cases of admitting some non-trivial congruence, or
axial symmetry of
all compact intersections, at least when the numbers of
connected components of ${\text{bd}}\,K$ and of
${\text{bd}}\,L$ are finite and when $K$ and $L$
are $C^2_+$?
Here a larger number of cases
can be expected. Namely, already our Theorem 5, about all
intersections in $H^2$,
has five cases (and even
Theorem 5, (A) comprises seven subcases). Moreover, 
\cite{4}, 
Theorem 3 about central symmetry of
compact intersections in
$S^d$, ${\Bbb{R}}^d$ and $H^d$, when restricted to $H^2$,
has five cases -- although three of these are contained in the
cases in our Theorem 5. Clearly in all of the cases in both
of these
theorems, all compact intersections have some non-trivial
symmetry.

\vskip.1cm


{\bf{Problem 4.}}
Is the finiteness hypothesis in Theorem 5 necessary?

\vskip.1cm


{\bf{Remark 3.}} 
As will be seen from the proof of Theorem 5, namely in the
proof of Lemma 5.10, rather than the finiteness hypothesis
in Theorem 5, we may suppose only the following. 
\newline
(1)
One of $K$ and $L$, e.g.,
the set $K$ has a boundary component $K_1$, 
such that some non-trivial arcs of $S^1$, each
with one endpoint at the infinite endpoints of $K_1$, resp.,
and lying on the
interior side of $K_1$ w.r.t.\ $K$, contain no infinite point
of any other boundary component of $K$. Moreover,
the other set $L$
has two neighbourly boundary components $L_1$ and $L_2$, with
$L_2$ following $L_1$ in the positive sense, such
that
passing on ${\text{bd}}\,L$, meant in $B^2$ containing the
model circle, in the positive sense from $L_1$
to $L_2$, we pass on no other connected components of
${\text{bd}}\,L$, meant in $H^2$.
\newline
As follows from Theorem 5, namely in the proofs of Lemmas
5.7--5.10, rather than the respective
finiteness hypotheses in Theorem 4 (in Part I, i.e., \cite{4})
and Theorem

\newpage

\noindent
5,
we may suppose the following.
\newline
(2)
One of $K$ and $L$ has two
boundary components with a common infinite point. (For Theorem
4 we have to consider that if all boundary components both of
$K$ and $L$ are straight lines, then from Theorem 5, (3) only
the case when both $K$ and $L$ are half-planes is possible.
However their intersection is not always
centrally symmetric.)


For the following Theorem 6 we
will need the following weakening of the $C^2$ property.
$$
\cases
{\text{Let for each }} x \in {\text{bd}}\,K, {\text{ and each }} 
y \in {\text{bd}}\,L,{\text{ there exist an }} 
\varepsilon _1 (x) > 0, \\
{\text{and an }} \varepsilon _1 (y)>0, {\text{ such that }}
K {\text{ and }}L {\text{ contain balls of radius }} \varepsilon _1 (x) \\
 {\text{and }} \varepsilon _1 (y), {\text{ containing }}
x {\text{ and }} y {\text{ in their boundaries, resp.}} 
\endcases
\tag **
$$
Moreover, we will need the following property, which
together with \thetag{**}
is a weakening of the $C^2_+$ property.
$$
\cases
{\text{Let for each }} x \in {\text{bd}}\,K, {\text{ and
each }} y \in {\text{bd}}\,L, {\text{ there exist an }}
\varepsilon _2 (x) > 0 {\text{ and}}
\\
{\text{an }} \varepsilon _2 (y) > 0, {\text{ such that
the following holds. The set of points of }} K
{\text{ and}}  
\\
L, {\text{ lying at a distance at most }}
\varepsilon _2 (x) {\text{ and }} \varepsilon _2 (y)
{\text{ from }} x {\text{ and from }} y, {\text{ is}}
\\
{\text{contained in a ball }} B {\text{ (for }} X=S^d,
\,\,{\Bbb R}^d), {\text{ or in a convex set }} B
{\text{ bounded}}
\\
{\text{by a hypersphere (for }} X = H^d). {\text{ Moreover,
bd\,}}B {\text{ has sectional curvatures}}
\\
{\text{at least }}\varepsilon _2 (x) {\text{ and }}
\varepsilon _2 (y), {\text{ and bd\,}} B
{\text{ contains }} x {\text{ or }} y, {\text{ resp.}} 
\endcases
\tag ***
$$
Clearly \thetag{**} implies
smoothness and \thetag{***} implies strict convexity,
resp.
Observe that both in \thetag{**} and \thetag{***}
$\varepsilon _i( x ) > 0$ and 
$\varepsilon _i(y) >
0$ can be decreased, 
and then \thetag{**} and \thetag{***} remain valid.
(\thetag{**} and \thetag{***} are (A) and (B) in
\cite{4}.) 


An {\it{inball of a convex body $Y \subset X$}} is a ball
of maximal radius contained in $Y$.


\proclaim{Theorem 6}
Assume \thetag{*}.
Let us assume $C^2$ for $K$ and $L$ (actually $C^2$ can be weakened to
\thetag{**}). 
For $X={\Bbb R}^d$ assume additionally 
that one of $K$ and $L$ has an extreme point. For
$X=H^d$ assume $C^2_+$ for $K$ and $L$ (actually $C^2_+$ can be weakened
to \thetag{**} and \thetag{***}).
Then we have $(1) \Longleftrightarrow (2) \Longleftrightarrow
(3) \Longleftrightarrow (4)$. Here
\roster
\item
There exists some $\varepsilon = \varepsilon (K,L) > 0$, such
that
for each $\varphi , \psi $, for 
which
${\text{\rm{diam}}} \, [ (\varphi K) \cap (\psi L) ]
< \varepsilon $, 
we have that $(\varphi K) \cap (\psi L)$
is centrally symmetric.
\item
For each $\varphi , \psi $ and each 
$x \in {\text{\rm{bd}}}\,K$
and $y \in {\text{\rm{bd}}}\,L$,
there exists some $\varepsilon = 
\varepsilon (K,L,x,y) > 0$,
such that
the following holds.
Suppose that $\varphi , \psi $ and $x, y$ 
satisfy the following
hypotheses {\rm{(A)}}, {\rm{(B)} and {\rm{(C)}}}.
\newline
{\rm{(A)}} 
${\text{\rm{diam}}} \, [ (\varphi K) \cap (\psi L) ]
< \varepsilon $.
\newline
{\rm{(B)}} $(\varphi K) \cap (\psi L)$ has a
unique inball $B_0$.
\newline
{\rm{(C)}} $\big[{\text{\rm{bd}}}\,[(\varphi K) \cap (\psi L)]
\big]
\cap ({\text{\rm{bd}}}\,B_0) =
\{ \varphi x, \psi y \} $, where
$\varphi x \in {\text{\rm{int}}}\,(\psi L)$ and 
$\psi y \in {\text{\rm{int}}}\,(\varphi K)$ are antipodal
on ${\text{\rm{bd}}}\,B_0$.
\newline
Then we have that $(\varphi K) \cap (\psi L)$
is centrally symmetric.
\item
The same as {\rm{(2)}}, but
``centrally symmetric'' in the conclusion replaced by
``symmetric w.r.t.\ the orthogonal bisector hyperplane
of $[\varphi x, \psi y]$''.

\newpage

\item
The connected components of the boundaries of both $K$ and $L$
are congruent 
spheres (for $X=S^d$ of radius at most $\pi /2$),
or paraspheres, or congruent hyperspheres (for ${\Bbb{R}}^d$
and $H^d$
degeneration to hyperplanes being not admitted). 
For the case of congruent spheres or
paraspheres we have that either $K$ and $L$ are congruent balls 
(for $X=S^d$ of radius at most $\pi /2$), or they are paraballs.
\endroster
\endproclaim


In Theorem 6, $(1) \Leftrightarrow (4)$ is the statement of
\cite{4}, 
Theorem 1.


\vskip.1cm


For Theorems 7 and 8, about (closed) convex hulls of unions
(which is the question
``dual'' to that about the intersection), 
we use the hypothesis
$$
\cases
X {\text{ is }} S^d, \,\,{\Bbb R}^d {\text{ or }} H^d,
{\text{ with }} d \ge 2. {\text{ Further, }} K,L \subsetneqq X
{\text{\, are closed convex}}
\\
{\text{sets with interior points. For }} S^2, {\text{ both }}
K {\text{ and }} L {\text{ lie in some open}}
\\
{\text{hemispheres. Moreover, }} \varphi , \psi : X \to X
{\text{ are orientation preserving congruences.}}
\\
{\text{For }} X = S^d {\text{ the set }} (\varphi K) \cup
(\psi L) {\text{ lies in an open hemisphere }} S {\text{ of }}
S^d .
\endcases
\tag ****
$$


The convex hull of a set $Y \subset H^d$ 
is defined as for ${\Bbb{R}}^d$.
For $Y \subset S^d$, 
since we will use only sets
$Y\,\,( = (\varphi K) \cup (\psi L) ) $ lying in an open
half-$S^d$, namely $S$ (cf.\ \thetag{****}), the convex hull of
such a $Y$ is the minimal convex set in $S$, containing $Y$.
We may suppose that $S$ is the open southern
hemisphere, thus the collinear model exists for it.
Thus, for any $X$,
the image in the collinear model of the convex hull of such
a $Y \subset X$ is the convex hull of 
the image in the collinear model of $Y$.

\vskip.1cm

We say that 
{\it{a convex body 
$Y$ in $S^d$, ${\Bbb{R}}^d$ or $H^d$, 
has at its boundary point $y$ 
a supporting sphere if the following holds.
There exists a ball containing $Y$, for $S^d$ of radius less
than $ \pi /2$, 
such that $y$ belongs to the boundary of
this ball}}. This boundary is called the {\it{supporting
sphere}}.
Observe that this implies compactness and strict convexity
of $Y$, and for $X = S^d$
that $Y$ is contained in an open hemisphere of $S^d$.
Therefore, if $K,L$ satisfy \thetag{****}, then
{\rm{conv}}\,$[ (\varphi K) \cup (\psi L) ] $ is
compact, hence closed. Therefore we will not need to take the
closure of this set.
Also, if $Y$ is
a convex
body, satisfying this italicized hypothesis for all $y \in $
bd\,$Y$, then for $S^d$,
${\Bbb{R}}^d$ and $H^d$,
any existing sectional curvature of $Y$ is greater
than $0$, $0$ and $1$, resp.  For $H^d$ this is a
serious geometric
restriction.

Let $D \subset X$ be a closed segment (to be specified later,
in Theorems 7 and 8),
with $D \subset S$ for $S^d$, and let 
$\varepsilon > 0$, with $\varepsilon $ sufficiently
small for $S^d$, cf. \thetag{*****} just below. 
Then we write $D(\varepsilon )$
for the (closed) {\it{convex hull of the union of two closed
$(d - 1)$-balls, of radius $\varepsilon $, lying
in hyperplanes orthogonal to $D$ at the
endpoints of $D$, and of centres the two
endpoints of $D$}}.
$$
\cases
{\text{For }} S^d {\text{ we choose }} \varepsilon > 0
{\text{ so that these two closed}}
\\
(d - 1){\text{-balls lie
in }} S {\text{ (or, equivalently, }} D(\varepsilon )
\subset S).
\endcases
\tag *****
$$
Observe that $D(\varepsilon )$
is centrally
symmetric w.r.t.\ the midpoint
of $D$, and also is symmetric
w.r.t.\ the orthogonal bisector hyperplane of $D$. 

Suppose that the midpoint of the segment $D$ is the south
pole $o$. 
Then the ball $B(D)$,
having $D$ as a diametral chord, has as image
in the collinear model a ball

\newpage

\noindent
about the centre $0$
of the model.
Further, the two bases of $D(\varepsilon )$ have as images
in the collinear model two congruent
closed $(d - 1)$-balls. These lie
in opposite tangent hyperplanes of the image of
$B(D)$ in the collinear model. Moreover, they
have centres the images 
of the two endpoints of $D$ in the collinear model.
Moreover, the image of 
$D(\varepsilon )$ in the collinear model is a right cylinder
over a $(d - 1)$-ball, with axis of rotation
(an axis of symmetry for $d = 2$) spanned by the image of $D$.


\proclaim{Theorem 7}
Assume \thetag{****}.
Let 
both $K$ and $L$ have supporting spheres 
at any of their
boundary points, for $S^d$ of radii less than $\pi /2$.
Then $(1) \Longleftrightarrow (2) \Longleftrightarrow (3)
\Longleftrightarrow (4) \Longleftrightarrow (5)
\Longleftrightarrow (6) \Longleftrightarrow (7)$. Here:
\roster
\item
For each $\varphi , \psi $ we have that 
{\rm{conv}}\,$[ (\varphi K) \cup (\psi L)] $ is 
centrally symmetric.
\item
For each $\varphi , \psi $,
provided that {\rm{conv}}\,$[ (\varphi K) \cup (\psi L)] $
has a unique diametral
segment $D$, we have that
{\rm{conv}}\,$[ (\varphi K) \cup (\psi L)] $ is 
symmetric w.r.t.\ the orthogonal bisector hyperplane of $D$.
\item
Suppose that $\varphi , \psi $ and
$x \in {\text{\rm{bd}}}\,K$
and $y \in {\text{\rm{bd}}}\,L$ 
satisfy the following
hypotheses {\rm{(A)}} and {\rm{(B)}}.
\newline
{\rm{(A)}} For $S^d$
we have that 
{\rm{diam}}\,$\big[ ${\rm{conv}}\,$[ (\varphi K) \cup (\psi L)
] \big] $ 
is sufficiently close to $\pi $. 
For ${\Bbb {R}}^d$ and $H^d$ we have
that this diameter is sufficiently large.
\newline
{\rm{(B)}}
{\rm{conv}}\,$[ (\varphi K) \cup (\psi L)
] $ has $D := [ \varphi x, \psi y]$ as its
unique diametral
segment.
\newline
Then 
{\rm{conv}}\,$[ (\varphi K) \cup (\psi L) ] $ is centrally symmetric.
\item
The same as {\rm{(3)}}, but
``centrally symmetric'' changed to
``symmetric w.r.t.\ the orthogonal bisector hyperplane
of $D$''.
\item
Suppose that $\varphi , \psi $ and
$x \in {\text{\rm{bd}}}\,K$
and $y \in {\text{\rm{bd}}}\,L$ 
satisfy 
hypotheses {\rm{(A)}} and {\rm{(B)}} from {\rm{(3)}}.
Then there exists an
$\varepsilon (K,L,x,y, \varphi , \psi ) > 0$
(for $S^d$ an $\varepsilon (K,L,x,y, \varphi , \psi, S) > 0$,
satisfying \thetag{*****} in place of $\varepsilon $
there),
such that the following holds. For
$0 < \varepsilon < \varepsilon (K,L,x,y, \varphi , \psi )$
(for $S^d$ for $0 < \varepsilon <
\varepsilon (K,L,x,y, \varphi , \psi, S)$)
we have the following. For $D := [\varphi x, \psi y]$ and for
$D(\varepsilon )$ defined before this theorem,
$\big[ ${\rm{conv}}\,$[ (\varphi K) \cup (\psi L) ] \big] \cap
D(\varepsilon )$ is centrally symmetric.
\item
The same as {\rm{(5)}}, but ``centrally symmetric'' changed to
``symmetric w.r.t.\ the orthogonal bisector hyperplane
of $D$''.
\item
$K$ and $L$ are congruent balls (for the case of $S^d$ of radius less than
$ \pi /2$).
\endroster
\endproclaim


The equivalence of (1) and (7) in 
Theorem 7
slightly improves Theorem 4 of \cite{4}. 
There we
needed the
additional hypothesis of smoothness of $K$ and $L$ for the
same conclusion (but could add hypothesis (A) from Theorem 7,
(3), to (1)).


\proclaim{Theorem 8}
Assume \thetag{****} with $d = 2$.
Also assume that both $K$ and $L$ 
have supporting circles
at any of their
boundary points,
for $S^2$ of radii less than $\pi /2$.
For $S^2$ also assume that $K$ and $L$ are $C^1$.

Then $(1) \Longleftrightarrow (2) \Longleftrightarrow (3)
\Longleftrightarrow (4) \Longleftrightarrow (5)
\Longleftrightarrow (6) \Longleftrightarrow (7)$. Here:
\roster
\item
For each $\varphi , \psi $ we have that
{\rm{conv}}\,$[ (\varphi K) \cup (\psi L) ] $ is axially symmetric.
\item
The same as {\rm{(1)}}, but ``is axially symmetric'' 
changed to ``admits some non-trivial congruence''.
\item
Suppose that $\varphi , \psi $
and
$x \in {\text{\rm{bd}}}\,K$
and $y \in {\text{\rm{bd}}}\,L$ 
satisfy hypotheses {\rm{(A)}}
and {\rm{(B)}}
from {\rm{(3)}} of Theorem {\rm{7}}. Then
{\rm{conv}}\,$[ (\varphi K) \cup (\psi L) ] $ is axially symmetric.

\newpage

\item
The same as {\rm{(3)}}, but ``is axially symmetric'' 
changed to ``admits some non-trivial congruence''.
\item
Suppose that $\varphi , \psi$
and
$x \in {\text{\rm{bd}}}\,K$
and $y \in {\text{\rm{bd}}}\,L$ 
satisfy hypotheses {\rm{(A)}}
and {\rm{(B)}} from {\rm{(3)}} of Theorem
{\rm{7}}.
Then there exists an $\varepsilon (K,L,x,y, \varphi , \psi)
> 0$
(for $S^2$ an $\varepsilon (K,L,x,y, \varphi , \psi, S) > 0$,
satisfying
\thetag{*****} in place of $\varepsilon $ there),
such that the following holds.
For $0 < \varepsilon < \varepsilon (K,L,x,y, \varphi , \psi )$
(for $S^2$ for $0 < \varepsilon <
\varepsilon (K,L,x,y, \varphi , \psi, S)$)
we have the following. For $D := [\varphi x, \psi y]$ and for
$D(\varepsilon )$ defined before Theorem {\rm{7}},
$\big[ ${\rm{conv}}\,$[ (\varphi K) \cup (\psi L) ] \big]
\cap D(\varepsilon )$ is axially symmetric.
\item
The same as {\rm{(5)}}, but ``is axially symmetric'' changed
to ``admits some non-trivial congruence''.
\item
$K$ and $L$ are {\rm{(in general incongruent)}} circles 
(for the case of $S^2$ of radii less than $ \pi /2$).
\endroster
\endproclaim


{\bf{Remark 4.}}
Probably smoothness for $S^2$ is not necessary for Theorem 8,
$(1) \Longleftrightarrow (7)$. 
However, without the supporting circle hypothesis,
Theorem 8, $(1) \Longleftrightarrow (7)$
does not remain true. Namely, besides the pairs $K,L$
in Theorem 8, (7),
also for the pairs $K,L$ 
in Part I (i.e., \cite{5}), 
Theorem 3, (2) (for ${\Bbb{R}}^2$), and in 
Theorem 5, (3),
(A), (B), (D), (E) (for $H^2$), the set
{cl\,conv}\,$[ (\varphi K) \cup (\psi L) ] $
admits some non-trivial congruences (actually, except for two
parallel strips in ${\Bbb{R}}^2$, axial symmetries).
Still there are some evident examples, with closed convex
hulls rather than convex hulls, satisfying \thetag{****},
except for $S^2$ its last statement:
for $S^2$ one of $K$ and $L$ is a halfsphere,
and for ${\Bbb{R}}^2$ one of $K$ and $L$ is either a parallel
strip, or a halfplane.

\vskip.1cm


{\bf{Problem 5.}} Are these the only pairs $K,L$ in $S^2$,
${\Bbb{R}}^2$ and $H^2$, with this property?

\vskip.1cm


Observe that (5) and (6) of Theorem 7 are 
``localizations'' of (3) and (4) of Theorem 7, resp.,
as well as
(5) and (6) of Theorem 8 are ``localizations'' of (3) and (4)
of Theorem 8, resp.
Thus these are analogues
of the ``localizations'' in
Part I (i.e., \cite{5}), 
Theorem 2, where (2), (6) and (7)
are ``localizations'' of (1), (4) and (5), resp.

\vskip.1cm


In the proofs of our Theorems we will use some ideas of
\cite{3}. 


\head 6. Preliminaries\endhead

For 1. Introduction, and 3. Preliminaries we refer to Part I
(i.e., \cite{5}). 
Most of our notations are standard. However,
we have to repete some notations from them.

$X$ is $S^d$, ${\Bbb R}^d$ or $H^d$, with $d \ge 2$. For
$x,y \in X$ we write $d(x,y)$ for their distance.
In $S^d$, when saying {\it{ball}}, or {\it{sphere}}, 
we always mean one with radius at most $\pi /2$.
A {\it{paraball}} is the
closed convex set in $H^d$, bounded by a parasphere. For
$d = 2$ we say {\it{paracircle}}.
For $d = 2$, and
$x_1,x_2$ on the boundary of a closed
convex set $K \subset X$ with interior points,
``close'' to each other,
we write
${\widehat{x_1x_2}}$ for the (shorter, or unique)
arc of ${\text{bd}}\,K$, which
set $K$ will be clear from the context.
We will use the collinear and the conformal models of $H^d$,
in the interior of the unit ball $B^d \subset {\Bbb{R}}^d$.   
We will speak about collinear and conformal
models of $S^d$ in ${\Bbb R}^d$. We mean by this the ones obtained by central
projection (from the centre), 
or by stereographic projection (from the north pole), to the
tangent hyperplane of $S^d$,
at the south pole, in ${\Bbb R}^{d+1}$. These exist of
course only on the open southern half-sphere, or on

\newpage

\noindent
$S^d$ minus the north pole,
resp. Their images are ${\Bbb R}^d$. We  call the
{\it{centre of the model}} the south pole of $S^d$.
The collinear and conformal models of ${\Bbb R}^d$ are meant as
itself, with {\it{centre}} the origin.
Sometimes we will consider the (collinear or conformal) model
circle of $H^2$ as the unit circle of the complex plane
${\Bbb{C}}$. Thus
we will speak about its points $1$, $i$, etc.
{\it{Smooth}} will mean differentiable, which for closed
convex sets with interior points is equivalent to $C^1$.
For a topological space $Y$ we say that some property of a
point
$y \in Y$ holds {\it{generically}}, if it holds outside a
nowhere dense closed subset.

A {\it{convex surface in}} $X$ is the boundary of a proper
closed convex subset of $X$
with interior points.
We will call a convex surface in ${\Bbb{R}}^d$
at some of its points {\it{twice
differentiable}} if the following holds.
It is the graph, in a suitable
rectangular coordinate
system, of a function having a Taylor series expansion of
second degree at this point,
with an error term $o( \| \cdot \| ^2 \| )$. 
By \cite{6}, 
pp. 31-32 (in both editions), convex surfaces in
${\Bbb{R}}^d$ are
almost everywhere twice differentiable.
This extends to $S^d$ and $H^d$ by using their collinear
models.

Part I (i.e., \cite{5}) 
contains Abstract,
1. Introduction, 2. New results: Theorems
1--4, 3. Preliminaries and 4. Proofs of Theorems 1--4.
Part II contains Abstract (about Part II), 5. New results:
Theorems 5--8, 6. Preliminaries and
7. Proofs of Theorems 5--8. The References in Part II contain
only items referred to in Part II; more detailed references, in
particular about hyperbolic geometry, cf.\ in Part I.


\head 7. Proofs of Theorems 5--8 \endhead


Before passing to the proof of Theorem 5, we introduce some
terminology. 

Suppose that $Y$ is a connected real analytic manifold, and
$f,g:Y \to
{\Bbb R}$ are analytic functions. Then either $f$ and $g$ coincide, or else
they cannot coincide on any subset of $Y$, whose closure
contains some
nonempty open subset of $Y$. This is the 
{\it{principle of analytic continuation}}.
Otherwise said, in the 
second case, generically, for $y \in Y$ (i.e., except on
a nowhere dense closed
subset of $Y$), we have $f(y) \ne g(y)$.

Recall that a finite union of nowhere dense closed subsets is itself nowhere
dense and closed. In {\bf{2}} of the proof of Lemma 5.4, in
{\bf{3}}--{\bf{5}} of the proof of Lemma 5.7, and in
{\bf{3}}--{\bf{6}} of Lemma 5.8
we will have the following
situation. On a connected real analytic manifold $Y$
(in fact, on
$H^2$, or on the orientation preserving congruences of $H^2$)
there are finitely many pairs of
analytic functions, $f_1,g_1; \ldots ;f_n,g_n:Y \to {\Bbb R}$,
say,
where $f_i$ and $g_i$ are different for all $1 \le i \le n$.
Then generically, for $y \in Y$, 
we have that $f_i(y) \ne g_i(y)$ for all $1 \le i \le n$.

\vskip.1cm

Before the proof of Theorem 5
we show a formula of trigonometrical type in $H^2$. It is in
a sense an analogue of the law of cosines for an angle of a triangle
in $H^2$. Namely,
the law of cosines allows us, for two circles, or radii $r, R$, and distance
of centres $c$, to determine the half central angle of the arc of the circle
of radius $r$, lying in the circle of radius $R$. We will need an analogous
formula, for a circle of radius $r$, and a hypercycle, with distance $l$ from
its base line, for the half central angle of the arc of the circle
of radius $r$, lying in one of the sets bounded by the
hypercycle. Here the
distance $c$ of the centre of the circle and the base line of the 

\newpage

\noindent
hypercycle is
given. We consider $c$ and $l$ as signed distances.
We admit degeneration to a straight line, i.e.,
{\it{we admit}} $l=0$. This formula is
surely known, but we could not locate a proof for it.
Therefore we give its simple proof.


\proclaim{Lemma 5.1}
Let $H \subset H^2$ be a hypercycle, with base line $B$, and
with
signed distance $l$ from $B$. Let $H^- \subset H^2$ be the set
of points at a signed distance at most $l$ from $B$.
Let $C \subset H^2$ be a circular line,
with centre $o$, radius $r$, and with $o$ having a signed
distance $c$ from $B$. Then $C \cap H \ne \emptyset $ 
if and only if $|c - l| \le r$. If this inequality is
satisfied, then for the half central angle $\omega $ of the
arc $C \cap H^-$ of $C$ we have
$$
{\text{\rm{sinh}}}\,l = {\text{\rm{cosh}}}\,r \cdot
{\text{\rm{sinh}}}\,c 
- {\text{\rm{sinh}}}\,r \cdot {\text{\rm{cosh}}}\,c \cdot
\cos \omega \,.
\tag 5.1.1
$$
\endproclaim


\demo{Proof}
We use the conformal model.
This shows that $C$ and $H$ have either two common points,
or they are tangent to each other, or they are disjoint.
(Their images are a circle, and a circular arc or a
segment which
cuts the model into two connected parts.) We suppose that
$B$ is a horizontal line containing the centre of the model,
and the signed distance from $B$
is positive in the open upper half-circle of the model.

Clearly we have $C \cap H \ne \emptyset $ if and only if 
$|c - l| \le r$. {\it{From now on we suppose this inequality.}}

Let $x$ be
one of the points of $C \cap H$, and let 
$y$ and $z$ be the orthogonal projections of $o$ and $x$
to $B$. Thus $d(x,z) = l$.
We let $d := d(x,y)$.
So we have to determine the angle $\omega = \angle xoy$
(for $o \in B$ the angle
$\omega $ is defined by the evident limit procedure).

By the law of cosines for sides
we have 
${\text{cosh}}\,d= {\text{cosh}}\,r \cdot {\text{cosh}}\,c
- {\text{sinh}}\,r
\cdot {\text{sinh}}\,c \cdot \cos \omega $.
Now we calculate the angle $\alpha := \angle oyx $ (for $o
\in B$ 
defined as a limit).
{\it{Preliminarily let us suppose $l \ne 0$}}, that
implies $d \ne 0$. Then
by the law of sines we have
$
\sin ^2 \alpha =\sin ^2 \omega \cdot {\text{sinh}}^2 r /
{\text{sinh}}^2 d \,.
$
Last, from the right triangle $yxz$ we have 
$
{\text{sinh}}^2 d(x,z) = \sin ^2 ( \pi /2 - \alpha ) \cdot
{\text{sinh}}^2 d\,.
$
So, fixing $r,c$, and supposing $\cos \omega $ as given, we determine, by
substitutions, successively, first cosh\,$d$, then $\sin ^2 \alpha $, then
${\text{sinh}}^2 d(x,z)$. This last expression should equal ${\text{sinh}}^2 l$.
Solving this last equation for $\cos \omega $ (which is a quadratic equation), 
we obtain, by rearranging,  
$$
\cases
\pm {\text{sinh}}\, l = {\text{cosh}}\,r \cdot {\text{sinh}}\,c 
- {\text{sinh}}\,r \cdot {\text{cosh}}\,c \cdot \cos \omega \\
={\text{cosh}}\,r \cdot {\text{cosh}}\,c \cdot
({\text{tanh}}\,c - {\text{tanh}}\,r
\cdot \cos \omega )
\,.
\endcases
\tag 5.1.2
$$
{\it{The same equation holds also for $l = 0$.}}
Namely then we have 
${\text{tanh}}\,c = {\text{tanh}}\,r \cdot \cos \omega $,
which is a trigonometric formula for a right triangle in
$H^2$.

We will show that \thetag{5.1.2} is satisfied with the plus sign,
as stated in \thetag{5.1.1} in the lemma. 
Recall that $|c - l| \le r$, or, equivalently, $C \cap H \ne
\emptyset $, was assumed.

First suppose $c \ge 0$. We distinguish two cases.
\newline
(1) We have $0 \le r \le c$.
\newline
(2)  We have $r > c \ge 0$.

In case (1) the entire circle $C$ lies (not strictly) above
$B$.
By $C \cap H \ne \emptyset $ then $H$ has some point at a
non-negative signed distance from $B$. Thus $l \ge 0$. On the
other hand, the expression in the middle of \thetag{5.1.2} 
lies in $[{\text{sinh}}\,(c-r), {\text{sinh}}\,(c+r)] \subset
[0, \infty )$.

\newpage

\noindent
Hence both sides of \thetag{5.1.1} are
non-negative, which means that in \thetag{5.1.2} the $\pm $ sign
is a $+$ sign.

In case (2) $C \cap B$ consists of two points, and
the line $B$ cuts $C$ to two non-degenerate arcs, both with
these two points as endpoints.
One of these arcs lies above $B$, the other one below $B$.
For these two points we have $l = 0$, and from \thetag{5.1.2}, the
angle $\omega _0$ corresponding to these points satisfies
${\text{tanh}}\,c = {\text{tanh}}\,r \cdot \cos \omega _0$.
Moreover, the points of $C$ corresponding to an angle $\omega
\in [0, \omega _0)$, or $\omega \in (\omega _0, \pi / 2]$,
lie strictly below or above $B$, resp.
Then $l < 0$, or $l > 0$,
resp. On the other hand, in the last expression in
\thetag{5.1.2}, the factor ${\text{tanh}}\,c - {\text{tanh}}\,r
\cdot \cos \omega $ vanishes for $\omega = \omega _0$, and
for $\omega
\in [0, \omega _0)$ or $\omega \in (\omega _0, \pi / 2]$ it
is negative or positive, resp. Hence once more, the
left and right hand sides of \thetag{5.1.1} have the same signs, 
which means that in \thetag{5.1.2} the $\pm $ sign
is a $+$ sign.

Last we extend the validity of \thetag{5.1.1} 
to $c < 0$. Let us apply \thetag{5.1.1} to $-c,-l, \pi - \omega
$ rather than $c,l, \omega $. Then the validity of
\thetag{5.1.1} for these
values implies its validity for $c,l, \omega $. 
$\blacksquare $
\enddemo


Later, in the proof of Theorem 5, {\it{we will consider the
case when $l \ge 0$;
then, of course, $c$ will vary in ${\Bbb R}$}}.

\vskip.1cm


{\bf{Remark 5.}}
There is a well-known formalism between spherical and
hyperbolical trigonometric relations. Namely,
the hyperbolical ones
can be obtained from the spherical ones by
``choosing the radius of the ball as $i$''. Applying this in
the converse direction,
\thetag{5.1.1} goes over to the following. Distance lines on
$S^2$ with distance $l$ are circles of radius $\pi /2 - l$,
and similarly distance $c$ from a line turns to distance
$\pi /2 - c$ from the centre of this circle, while $\omega $
turns to $\pi - \omega $. Thus we obtain
$\cos l = \cos r \cdot \cos c + \sin r \cdot \sin c \cdot \cos
\omega $, i.e., the cosine law for sides.

\vskip.1cm


In the proof of Theorem 5 we will allow that hypercycles
degenerate to straight lines, i.e., {\it{hypercycle will
mean either a proper hypercycle, or a straight line; i.e., the
curvature is allowed to be $0$ as well}}. The
base line of a straight line is meant as itself.

\vskip.1cm


\demo{Proof of Theorem {\rm{5}}}
{\bf{1.}}
The implication $(2) \Longrightarrow (1)$ is trivial.


\proclaim{Lemma 5.2}
Assume \thetag{*} with $d = 2$ and $X = H^2$.
Then the
implication $(3) \Longrightarrow (2)$ of Theorem {\rm{5}} holds.
\endproclaim


\demo{Proof}
{\bf{1.}}
We begin with case (A).
A circle is axially symmetric w.r.t.\ any straight line
passing through its centre. 
A paracircle is axially symmetric w.r.t.\ any straight
line passing through its centre (its point at infinity). 
A half-plane, or a convex domain bounded by
a hypercycle is axially symmetric w.r.t.\ 
any straight line, that intersects its boundary, or its
base line
orthogonally, resp. These imply, that if any of $\varphi K$ 
and $\psi L$ is either a circle or a
paracircle, then their intersection is axially symmetric
w.r.t.\ 
(any) straight line joining their centres. Suppose
that one of $\varphi K$ and $\psi L$ 
is a circle, and the other one is
a convex set bounded by a hypercycle, or is a half-plane.
Then the straight
line passing
through the centre of the circle, and orthogonal to the base line of the
hypercycle, or to the boundary of the half-plane, 
is an axis of symmetry of the intersection. 

Last, let 
$\varphi K$ and $\psi L$ be congruent convex sets, each bounded by one
hypercycle,

\newpage

\noindent
or let them be two half-planes. 
Consider the base lines of these hypercycles, or the boundaries of 
these half-planes, resp. There
are four cases. These base lines
\newline 
(a) may coincide; or
\newline 
(b) may transversally intersect; or
\newline 
(c) may be parallel but distinct; or 
\newline 
(d) may be ultraparallel.
\newline In case (a) any straight line orthogonal to the
common base line is an axis of symmetry. In case (b),
$\varphi K \ne \psi L$, and
${\text{bd}}\,(\varphi K)$ and ${\text{bd}}\,(\psi L)$
intersect
transversally at some point $p$ (for this use the conformal model). 
Then $(\varphi K) \cap (\psi L)$ has an inner
angle at $p$, of measure less than $ \pi $, and the bisector
of this angle is an axis of symmetry of 
$(\varphi K) \cap (\psi L)$. In case (c), if one of $\varphi K$ and $\psi L$
contains the other, the intersection is evidently axially
symmetric. Otherwise, the symmetry axis of the base lines is an axis of
symmetry of the intersection.
In case (d), we consider the pair of points on the base lines, realizing
the distance of these lines. The straight line connecting these points
is orthogonal to both
lines, and is an axis of symmetry of $(\varphi K) \cap (\psi L)$.

{\bf{2.}}
We continue with case (B). Suppose that $K$
is a circle, and $L$ is bounded
by two hypercycles, whose base lines coincide
(one of them possibly degenerating to a straight line).
Then the straight line passing through
the centre of $\varphi K$, and orthogonal to the above base line, is an axis of
symmetry of $(\varphi K) \cap (\psi L)$.

{\bf{3.}}
We continue with case (C). Suppose that $K$ is a circle of radius $r$, and
the boundary hypercycle or straight line 
components of $L$
have pairwise distances at least $2r$. Then 
int\,$(\varphi K)$ can intersect at most one boundary component of $\psi L$. 

If int\,$(\varphi K)$ does not intersect any boundary component of $\psi L$
(and, by hypothesis, int\,$[(\varphi K) \cap (\psi L)] \ne \emptyset $),
then $(\varphi K) \cap (\psi L) = \varphi K$ is a circle.
Hence it is axially
symmetric.

Suppose that int\,$(\varphi K)$ intersects exactly one boundary component $\psi L_1$ 
of $\psi L$.
Then $(\varphi K) \cap (\psi L)$ is the same as the intersection of 
$\varphi K$ and of 
the closed convex set, bounded by $\psi L_1$,
and containing $\psi
L$. This has an axis of symmetry, cf.\ case (A).

{\bf{4.}}
We continue with case (D). Let $K$ be a paracircle,
and $L$ a parallel
domain of some straight line, with some distance $l>0$. 
Consider the common base line of the two
hypercycles, bounding $\psi L$. If the infinite point of the paracircle
$\varphi K$ lies on this common base 
line, then this straight
line is an axis of symmetry of $(\varphi K) \cap (\psi L)$. 
Suppose that the infinite point of $\varphi K$ 
does not lie on this common base line. Then there is a unique
straight line that passes through the infinite point of $\varphi K$, 
and is orthogonal
to the common base line. Then this unique straight line is an axis of
symmetry.

{\bf{5.}}
Last we turn to case (E). Consider the common base lines of the two
hypercycles bounding $\varphi K$, and of the two
hypercycles bounding $\psi L$. These two straight 
lines can coincide, or can transversally intersect, or can be
parallel but distinct, 
or can be ultraparallel. In any case there is
an axial symmetry interchanging these two straight lines. This axial symmetry
interchanges the parallel domains of these straight lines, with distance $l$,
as well.
Hence it is an axial symmetry of the intersection of these parallel domains, i.e., of 
$(\varphi K) \cap (\psi L)$.
$\blacksquare $
\enddemo


\newpage

{\it{Proof of Theorem}} 5, {\bf{continuation.}} {\bf{2.}}
Last we turn to the proof of $(1) \Longrightarrow (3)$.
By Part I (i.e., \cite{5}), 
Theorem 1, $(2) \Longrightarrow (3)$ 
we know that each boundary component of 
both $K$ and $L$ is 
either a cycle, or a straight line. Thus,
for both $K$ and $L$, we have the following
possibilities: it is a
circle, or a paracircle, or its boundary components are 
hypercycles and straight lines. 

We make a case distinction. Either both bd\,$K$ and bd\,$L$
are connected, or
one of them has several connected components.


\proclaim{Lemma 5.3}
Assume \thetag{*} with $d = 2$ and $X = H^2$.
Suppose {\rm{(1)}} of Theorem {\rm{5}}. If
both {\rm{bd}}\,$K$ and {\rm{bd}}\,$L$ are connected, 
then we have {\rm{(A)}} of {\rm{(3)}} of Theorem {\rm{5}}.
\endproclaim


\demo{Proof}
We have to investigate the cases when 
\newline 
(a) $\varphi K$ and $\psi L$ are 
one paracircle and one convex
set bounded by a hypercycle or a straight line, or 
\newline 
(b) $\varphi K$ and $\psi L$ are two incongruent convex sets, both
bounded by a hypercycle or a straight line,
\newline
and in both cases we have to find a contradiction.

Now it will be convenient to use the conformal model.
In case (a), let the centre of the
paracircle be one endpoint of the base line of the hypercycle,
or one endpoint
of the straight line. Then
${\text{bd}}\,(\varphi K)$ and ${\text{bd}}\,(\psi L)$
transversally intersect, at a single point, and this
point $p$ is the only non-smooth point of $(\varphi K) \cap 
(\psi L)$.
In case (b), let the base lines of the two
hypercycles, or the base line of the hypercycle and the
straight line transversally
intersect, resp.\ (two straight lines cannot occur).
Then, also 
${\text{bd}}\,(\varphi K)$ and ${\text{bd}}\,(\psi L)$
transversally intersect, at a single point, and this
point $p$ is the only non-smooth point of $(\varphi K) \cap 
(\psi L)$. 

Both in case (a) and (b),
any non-trivial congruence admitted
by $(\varphi K) \cap (\psi L)$
would have $p$ as a fixed point. Moreover, the pair of the
half-tangents at $p$
would be preserved by this congruence. Thus, it would be an axial
symmetry, 
w.r.t.\ the 
angle bisector of the inner angle of $(\varphi K) \cap (\psi L)$ at $p$. Thus,
this axial symmetry should
interchange the portions of the boundaries of $\varphi K$ and $\psi L$, 
bounding $(\varphi K) \cap (\psi L)$. However, these portions of boundaries
have different
curvatures, which is a contradiction. 
$\blacksquare $
\enddemo


{\it{Proof of Theorem}} 5, {\bf{continuation.}} {\bf{3.}}
By Lemma 5.3 {\it{there remained the case when one of $K$
and $L$ has at least two boundary components.}}
Observe that this rules out the cases when
$K,L$ are two circles, or two paracircles, or one circle and one 
paracircle. There remain the cases when one of $K$ and $L$ is
bounded by several but finitely many
hypercycles and straight lines, and the other one either is
a circle, of some radius $r$, or is
a paracircle, or also is bounded
by finitely many (possibly one)
hypercycles and straight lines. We will investigate these
three cases separately, in Lemma 5.4, Lemma 5.5 and Lemmas
5.7--5.10, resp.

Suppose that one of $K$ and $L$ is bounded
by hypercycles and straight lines. Then the boundary
components $\varphi K_i$ of
$\varphi K$, or the boundary components $\psi L_i$ of
$\psi L$ have a natural cyclic
order, in the positive sense, on bd\,$(\varphi K)$, or
bd\,$(\psi L)$,
resp. We associate to $\varphi K$, or to $\psi L$ a graph.
Its
vertices
are the infinite points of the $(\varphi K_i)$'s, or
$(\psi L_i)$'s, resp., and between two such
points there is an edge, if they are the two infinite points of some $\varphi K_i$, or
$\psi L_i$, resp. We say that this edge is
$\varphi K_i$,

\newpage

\noindent
or $\psi L_i$,
resp. This graph either
is a union of vertex-disjoint paths, or is
a cycle. Here we admit a cycle of length $2$, shortly a
$2$-cycle.
Then the graph consists
of two vertices, and two edges between these two vertices,
which are
two $(\varphi K_i)$'s ($(\psi L_i)$'s) with both infinite
points in common. 

If we have two edges in these graphs with a common vertex, and they are e.g.,
$\varphi K_1$ and $\varphi K_2$, then by this notation we will mean 
that $\varphi K_2$ follows $\varphi K_1$ on bd\,$(\varphi K)$
in the positive sense.
If $\varphi K_1$ and $\varphi K_2$ 
form a $2$-cycle, then the notation is fixed some way.
(A similar convention holds for $\psi L$).

$$
\cases
{\text{Later in the proof of Theorem 5 we will allow that
hypercycles}}
\\
{\text{degenerate to straight lines, i.e., {\it{hypercycle
will mean either a}}}}
\\
{\text{{\it{proper hypercycle, or a straight line; i.e., the
curvature is allowed}}}}
\\
{\text{{\it{to be $0$ as well}}. The base line of a straight
line is meant as itself.}}
\endcases
\tag 5.1
$$


\proclaim{Lemma 5.4}
Assume \thetag{*} with $d = 2$ and $X = H^2$.
Suppose {\rm{(1)}} of Theorem {\rm{5}}.
Let $K$ be a circle of radius $r$ and centre $o$, and
let the connected boundary components of $L$ be
at least two
hypercycles or straight lines. Then {\rm{(B)}} or
{\rm{(C)}} of {\rm{(3)}} of Theorem {\rm{5}} holds.
\endproclaim


\demo{Proof}
{\bf{1.}}
We have to show that either $L$ is bounded by
two hypercycles with common base line, or
the at least two boundary components of $L$
have pairwise distances at least $2r$.
Let us
suppose the contrary. That is, we have both that
dist\,$(L_1,L_2) < 2r$, for some different boundary
components $L_1$ and $L_2$ of $L$, and
that $L$ is not bounded by
two hypercycles with a common base line.
By dist\,$(\psi L_1,\psi L_2) < 2r$,
we have, 
for some choice of $\varphi $, for $\varphi _0$, say, that 
int\,$(\varphi _0 K)$ intersects both $\psi L_1$ and
$\psi L_2$. Then
the same statement holds for all sufficiently small
perturbations $\varphi $ of $\varphi _0$,
hence also of the {\it{centre
$\varphi _0 o$ of}} $\varphi _0 K$. 

$$
\cases
{\text{We are going to show that for an arbitrarily small,
generic perturbation }} 
\\
\varphi o {\text{ of }} \varphi _0 o {\text{ we have that }}
(\varphi K) \cap (\psi L) {\text{ admits no non-trivial
congruence.}}
\endcases
\tag 5.4.1
$$
Observe that \thetag{5.4.1} contradicts the hypotheses of
this lemma, thus our above indirect assumption will lead to a
contradiction.

{\bf{2.}}
We have that int\,$(\varphi _0 K)$ intersects $\psi L_i$
for $i \in I(0,{\text{i}})$ (where $|I(0,{\text{i}})| \ge 2$
by $\{ 1,2 \} \subset I(0,{\text{i}})$), and $\varphi _0
K$ touches $\psi L_i$ from outside (w.r.t.\ $\psi L$) 
for $i \in I(0,{\text{t}})$, and
$\varphi _0 K$ is disjoint to $\psi L_i$ for 
$i \in I(0,{\text{d}})$.
For a sufficiently small perturbation $\varphi $ of
$\varphi _0$ we have for the analogously defined
$I({\text{i}})$, $I({\text{t}})$
and $I({\text{d}})$ that
$\{ 1,2 \} \subset I(0,{\text{i}}) \subset
I({\text{i}})$ and $I(0,{\text{d}}) \subset
I({\text{d}})$,
hence also
$$
I({\text{i}}) \cup I({\text{t}})
\subset I(0,{\text{i}}) \cup I(0,{\text{t}}) .
\tag 5.4.2
$$
(However, generically $I({\text{t}}) = \emptyset $.)
Moreover, each of $I(0,{\text{i}})$, $I(0,{\text{t}})$,
$I({\text{i}})$
and $I({\text{t}})$ is finite.
(Observe that any compact set in $H^2$, thus, e.g., a closed
circle concentric with $\varphi _0
K$, of radius just a bit larger than $r$,
intersects only finitely many
$(\psi L_i)$'s.)

\newpage

Let $\varphi $ be a small generic perturbation of
$\varphi _0$, thus $\varphi o$ be a small generic
perturbation of $\varphi _0 o$. Then
$(\varphi K) \cap (\psi L)$ is 
a convex body. It is
bounded, alternately, by finitely many, but
at least two non-trivial arcs of
bd\,$(\varphi
K)$, and by same many non-trivial arcs of some $(\psi L_i)$'s,
for different
$(\psi L_i)$'s. (The case when some $\varphi K$ touches
$\psi L_i$
at some point, from inside
(w.r.t.\ $\psi L$), gives that $\psi L_i$ does not cut off
anything from $\varphi K$, hence
can be left out
of consideration. The case when some $\varphi K$ touches
$\psi L_i$
at some point, from outside (w.r.t.\ $\psi L$), gives that
$|(\varphi K) \cap (\psi L)| = 1$, contradicting  
$\left( {\text{int}}\, (\varphi K) \right) \cap (\psi L_1)
\ne \emptyset $.
Thus,
by the conformal model, we may
suppose that the neighbouring boundary arcs of
$(\varphi K) \cap (\psi L)$ intersect transversally.) 
The curvatures of these arcs are greater than $1$, or
smaller than
$1$, resp., so each congruence admitted by $( \varphi K) \cap (
\psi L)$
preserves both above types of arcs, separately.
Thus the union of those boundary arcs of $(\varphi K) \cap (
\psi L)$, which are arcs of ${\text{bd}}\,(\varphi K)$, i.e.,
the set $[ {\text{bd}}\,(\varphi K) ] 
\cap (\psi L)$ is preserved
by each congruence admitted by
$( \varphi K) \cap ( \psi L)$. Therefore
also conv\,$\left[ [ {\text{bd}}\, (\varphi K) ] 
\cap (\psi L) \right] $ is preserved
by each congruence admitted by
$( \varphi K) \cap ( \psi L)$. Then
{\it{{\rm{conv}}\,$\left[ [ {\text{\rm{bd}}}\, (\varphi K) ] 
\cap (\psi L) \right] $ is 
obtained from the circle $\varphi K$, by cutting off disjoint
circular segments, relatively open in $\varphi K$,
by finitely many, but at least two disjoint
non-trivial chords. These 
have endpoints the points of intersection of 
the single $(\psi L_i)$'s with {\rm{bd}}\,$(\varphi K)$, 
for all $(\psi L_i)$'s with $i \in I({\text{\rm{i}}})$.}} 
$$
\cases
{\text{We will attain that all these chords have
different lengths --}}
\\
{\text{moreover, this can be attained by a
small, generic motion}}
\\
{\text{of the centre }} \varphi o {\text{ of }} \varphi K ,
{\text{ from its original position }} \varphi _0 o.
\endcases
\tag 5.4.3
$$ 

By the law of sines,
two such chords have equal length if and only if the
corresponding half central angles, belonging to $(0, \pi )
\subset [0, \pi ]$ have equal sines. 
That is, the absolute values of the cosines of these angles
coincide.  
We use formula \thetag{5.1.1}
from Lemma 5.1. However, in Lemma 5.1 we dealt with the half
central angle $\omega $ of the arc $C \cap H^-$ of $C$, while  
here we are interested in the central angle $\chi $ of the
arc of ${\text{bd}}\,(\varphi K)$, ``cut off'' by $\psi L_i$.
Then $\chi = \pi - \omega $.
We have several such equations like \thetag{5.1.1},
corresponding to the
$(\psi L_i)$'s, for $i \in I({\text{i}})$,
with respective fixed values $r$ and
$l_i$, and variable values
$c_i$. We choose
$l_i \ge 0$ for each $i$ (this is the choice of the sign
of the signed
distance from the base line of $\psi L_i$),
and then $c_i$ will vary in ${\Bbb R}$. For $l_i = 0$ we
should have
positive signed
distance from $\psi L_i$ on the exterior side of the
straight line
$\psi L_i$, w.r.t.\ $\psi L$.
We will use arbitrarily small generic
perturbations $\varphi o$ of the centre $\varphi _0 o$ of our
circle $\varphi _0 K$.

The half central angles of the arcs of
${\text{bd}}\,(\varphi K)$, ``cut off'' by
$\psi L_i$, for $i \in I({\text{i}})$,
are denoted by
by $\chi _i$.
Then $\sum _{i \in I({\text{i}})} \chi _i < \pi $. Then
for each
distinct $i,j \in I({\text{i}})$ (recall from the paragraph
following \thetag{5.4.1} that $| I({\text{i}}) | \ge
| I({\text{0,i}}) | \ge 2$)
we have
$\chi _i +  
\chi _j \le \sum _{i \in I({\text{i}})} \chi _i < \pi $.
Hence
for $| \cos \chi _i | = | \cos \chi _j |$ we
have $\chi _i = \chi _j < \pi /2$. (Observe that
$|I(0,{\text{i}})| \ge
2$ was essential in this proof. For $I({\text{i}}) =
\emptyset $ and
$I({\text{t}}) = \{ 1 \} $ we have the example when
$\varphi K$ touches
$\psi L_1$ from outside w.r.t.\, $\psi L$,
hence it does not intersect any other
$\psi L_i$, and where the analogously defined $\chi _1$
equals $\pi $, hence $\sum _{i \in I({\text{i}}) \cup
I({\text{t}})} \chi _i = \pi $.)

Hence, it suffices to exclude all pairwise equalities
of finitely many expressions 

\newpage

\noindent
for $\cos \chi = - \cos \omega $
-- obtained from solving equations \thetag{5.1.1}, 
for all $i \in I({\text{0,i}}) \cup I({\text{0,t}})$ (which
set contains $I({\text{i}})$ by \thetag{5.4.2}) --
namely those of the form on the right hand side of the
following formula:
$$
\cos \chi _i =
({\text{{sinh}}\,l_i -
{\text{cosh}}\,r \cdot {\text{sinh}}\,c_i) /
({\text{sinh}}\,r \cdot \text{cosh}}\,c_i) \,.
\tag 5.4.4
$$
Observe that all these expressions 
are analytic in
$\varphi o$, since the $c_i$'s are
analytic in $\varphi o$, and $r$ and the $l_i$'s are fixed.
(In fact, in the collinear model the signed distance of a
point
$(\xi , \eta )$ from a horizontal line of equation
$\eta = \eta _0 \in
(-1, 1)$ is the signed
hyperbolic
length of $[(\xi , \eta _0), (\xi ,\eta )]$,
which can be expressed by the
cross-ratio formula, cf.\ \cite{7}.) 

Moreover, none of these equalities, for distinct $i,j$ in 
$I(0,{\text{i}}) \cup I(0,{\text{t}})$, is an identity.
Namely, let $\varphi ' K$ touch a 
boundary component $\psi L_i$ of $\psi L$
from outside w.r.t.\ $\psi L$.
Then $(\varphi ' K) \cap (\psi L)$ consists of one point. So, in
particular, $\varphi K$ {\it{does not intersect any other
boundary
component $\psi L_j$ of $\psi L$}}.
Then in \thetag{5.4.4} $\cos \chi _i = -1$, however for $j \in
[I({\text{0,i}}) \cup I({\text{0,t}})]
\setminus \{ i \} $ the value for $\cos \chi _j$ from the
right hand side of
\thetag{5.4.4} is not in $[-1,1]$. 
(Else
\thetag{5.1.1} and the second paragraph of the proof
of Lemma 5.1
would give $l_j \in [c_j - r, c_j + r]$, so 
$|c_j - l_j| \le r$, hence
$[ {\text{bd}}\,(\varphi K) ]
\cap (\psi L_j) \ne \emptyset $, which is a contradiction.)
Hence
the $i$'th expression 
in the right hand side of
\thetag{5.4.4}, and any other $j$'th expression in
the right hand side of
\thetag{5.4.4}, for
$i,j \in I({\text{0,i}}) \cup I({\text{0,t}})$, 
are not identical. 

Therefore all our finitely many analytic equations,
that the right hand sides of \thetag{5.4.4} are equal for
distinct $i,j \in 
I({\text{0,i}}) \cup I({\text{0,t}})$, are not identities.
Hence each of them holds only for $\varphi o$ belonging to a 
nowhere dense closed subset. Therefore, supposing that
$\varphi o$ belongs to some small closed
neighbourhood $A$ of $\varphi _0 o$,
except for $\varphi o$ belonging to
a nowhere dense closed subset $B$ of $A$, 
none of these equations hold.
(Hence, by \thetag{5.4.2}, also
for any sufficiently small generic perturbation
$\varphi $ of $\varphi _0$, and for distinct $i,j \in 
I({\text{i}})$, none of these equations holds.) 
That is, we have
proved what was claimed in \thetag{5.4.3}.

{\bf{3.}}
From now on we will suppose \thetag{5.4.3}.
There are 
two possibilities. Either
\newline
(1) we can have at least three such chords, as in
\thetag{5.4.3},
or
\newline
(2) we always have exactly two such chords (then these
correspond to the above $L_1, L_2$).

Any congruence admitted by $(\varphi K) \cap (\psi L)$ 
preserves $\varphi K$, 
and conv\,$\big[ [  {\text{bd}}\, (\varphi K) ] \cap $
\newline
$ (\psi L) \big] $, and
also the above mentioned at least three, or exactly two 
disjoint chords, since their lengths are different. However, a
non-trivial
congruence admitted by $(\varphi K) \cap (\psi L)$,
preserving a single such
chord (and that side of the chord, on which
conv\,$\big[ [  {\text{bd}}\, (\varphi K) ] \cap (\psi L) \big]
$ lies), 
is an axial symmetry w.r.t.\ the
orthogonal bisector straight
line of this chord, which line
contains $\varphi o$. Since this is the unique non-trivial
congruence admitted by $(\varphi K) \cap (\psi L)$, any other
one of these chords gives the same axial symmetry, hence also
the same axis of symmetry of $(\varphi K) \cap (\psi L)$.
However, there are no three disjoint 
circular segments, relatively open in $\varphi K$,
cut off by disjoint
chords orthogonal to a single
straight line, their common orthogonal bisector.
Thus in case (1) we have a
contradiction.

There remains case (2) when we can have only exactly two such
chords.
These

\newpage

\noindent
must 
correspond to the
above considered hypercycles $\psi L_1$ and $\psi L_2$. 

Then $\varphi K$ does not even touch 
any other $\psi L_i$ from outside w.r.t.\ $\psi L$,
since then by a small motion of $\psi L$
we could attain that 
int\,$(\varphi K)$ intersects at least three $(\psi L_i)$'s,
which
case was settled just above.

The above reasoning gives that, in this case, the orthogonal
bisecting straight
lines
of the two chords coincide, furthermore, contain $\varphi o$.
However, the orthogonal bisecting straight
lines of these two chords are orthogonal also to
$\psi L_1$ and $\psi L_2$, and hence also to the
base lines of $\psi L_1$ and $\psi L_2$. We have that these base lines are different,
since their coincidence was excluded in the first paragraph
of {\bf{1}}
of the proof of this lemma. 
Then they
do not intersect, and thus are either parallel and distinct,
or ultraparallel.

If they are parallel but distinct, then they
admit no common orthogonal straight line, so we have a
contradiction. 

If they are ultraparallel,
then they have exactly one common orthogonal straight
line. This straight line depends on $\psi L$ only -- let it be
$l(\psi L)$. From above $\varphi o \in l(\psi L)$. Then
in case (2) for
$\varphi o \in A \setminus [ B \cup \left( A \cap
l(\psi L) \right) ]$ we have a contradiction.

Clearly $B \cup \left( A \cap
l(\psi L) \right) $ is a nowhere dense closed subset of $A$. 
Hence \thetag{5.4.1} holds for an arbitrarily small generic
perturbation $\varphi o$ of $\varphi _0 o$, contrary to the
hypotheses of this lemma. 
This shows that our indirect hypothesis in the first
paragraph of {\bf{1}} led to a contradiction,
also in case (2).

Since both cases (1) and (2) led to a contradiction, therefore
our indirect hypothesis was false.
In other words, the statement of this lemma holds.
$\blacksquare $
\enddemo


In the following lemmas, we will consider the (collinear or
conformal)
model circle of
$H^2$ as embedded in the complex plane ${\Bbb{C}}$. 


\proclaim{Lemma 5.5}
Assume \thetag{*} with $d = 2$ and $X = H^2$.
Suppose {\rm{(1)}} of Theorem {\rm{5}}.
Let $K$ be a paracircle, and let the connected boundary
components of $L$ be
at least two hypercycles or straight
lines. Then {\rm{(D)}} of {\rm{(3)}} of Theorem {\rm{5}} holds.
\endproclaim


\demo{Proof}
We consider the graph of $\psi L$. It has at least two edges,
by the hypothesis made after the proof of Lemma 5.3. Let 
$\psi L_1$ be one of its edges. Let one endpoint of $\psi L_1$,
say, $1$,
coincide with the centre of $\varphi K$. Then  
$(\varphi K) \cap (\psi L)$ has the unique infinite point $1$.
Then each congruence admitted by $(\varphi K) \cap (\psi L)$
preserves this point $1$, hence
a congruence is either an axial symmetry w.r.t.\ an axis
containing $1$, or is a rotation about the infinite point $1$.
(Observe that by Part I (i.e., \cite{5}), 
\thetag{2} 
it cannot be a nontrivial
translation or a glide-reflection which is not a reflection,
along a straight
line having $1$ as one of its infinite points.)
However, the second case is impossible by Part I (i.e.,
\cite{5}), 
\thetag{3}. 

In the first case let us consider the set 
$S = (\varphi K) \cap (\psi L_1) \subset {\text{bd}} \,[
(\varphi K) \cap
(\psi L) ] $, which is a closed semiinfinite arc of a
hypercycle, with
infinite point at $1$, and
finite endpoint on ${\text{bd}}\,(\varphi K)$, and else
lying in ${\text{int}}\,(\varphi K)$.
Then the
image $S'$ of $S$ by our axial symmetry
is a hypercycle arc with the same above listed properties. 
Since
$S' \subset
{\text{bd}} \,[ (\varphi K) \cap (\psi L) ] \subset
[{\text{bd}}\,(\varphi K) ]
\cup [ {\text{bd}}\,(\psi L) ] $,
therefore
$S'$, minus its finite endpoint, lies in ${\text{bd}}\,(\psi L)$.
Then taking closure also $S' \subset
{\text{bd}}\,(\psi L)$. Then $S' \subset \psi L_2$ for some
boundary component $\psi L_2$ of $\psi L$, congruent to 
$\psi L_1$.
Let $\psi L_2$ follow $\psi L_1$ at their common infinite
point $1$ in the positive sense. If also
their other infinite

\newpage

\noindent
points coincide, then $\psi L$ is a
parallel domain of their coincident base lines, and then
we have the statement of the lemma.
$$
\cases
{\text{Therefore, for contradiction, let us suppose that the
other endpoints of }} \psi L_1 {\text{ and}}
\\
\psi L_2
{\text{ are different, and we are going to show that this
leads to a contradiction.}}
\endcases
\tag 5.5.1
$$

Since any three
distinct points of the boundary circle of the model can be
taken over by (the extension of) an orientation preserving
congruence
to any other three distinct points of the boundary
circle of the model, of the same  orientation, we may
suppose that the other endpoints of 
$\psi L_1$ and $\psi L_2$ are, in the closure in ${\Bbb{C}}$
of the model circle, very close to $-1$. Then all
other boundary components of $\psi L$ are, in the
closure in ${\Bbb{C}}$ of the
(conformal or collinear) model circle, close to $-1$. 

We consider the conformal model. We fix the position of
$\psi L_1$ and $\psi L_2$, but 
we choose another $\varphi K$ so that 
it touches the boundary of the model circle
very close (in ${\Bbb{C}}$) to $1$, at a point
in the open upper
half-plane, and its interior
intersects both $\psi L_1$ and $\psi L_2$, and
its image in the model is a circle of very small radius (in
${\Bbb{C}}$).
Then $\varphi K$ 
is far (in ${\Bbb{C}}$)
from all other boundary components of $\psi L$. Hence,
$(\varphi K) \cap (\psi L)$ is 
a convex arc-quadrangle, bounded by an arc of $\varphi K$, an arc
of $\psi L_2$, another arc of $\varphi K$, 
and an arc of $\psi L_1$, in this cyclic order, in the positive
sense. The endpoints of these arcs
(the vertices of our arc-quadrangle)
are the only non-smooth
points of 
this arc-quadrangle, 
hence they are preserved by
any congruence admitted by this arc-quadrangle.
This arc-quadrangle has two vertices on $\psi L_1$, and
two vertices on $\psi L_2$,
hence its vertices are the
vertices of a strictly convex quadrangle (in the orientation
inherited from ${\text{bd}}\,[(\varphi K) \cap (\psi L)]$).
This arc-quadrangle is also compact, so by
Part I (i.e., \cite{5}), 
\thetag{2} 
and \thetag{3}, 
it can admit only congruences which are rotations,
or axial symmetries (according to whether they are
orientation preserving, or orientation reversing).

Since a paracycle and a hypercycle (straight line) have
different curvatures,
neither of the 
diagonals can be an axis of symmetry, and there is no
congruence which is a
combinatorial $4$-fold rotation. Therefore, our arc-quadrangle
can admit only a
rotation which is a combinatorial
central symmetry, and
by Part I (i.e., \cite{5}), 
\thetag{4}, 
this is also a geometrical central symmetry,
with a unique centre.
So, a non-trivial
congruence admitted by 
$(\varphi K) \cap (\psi L)$ is either a central symmetry, or
an axial symmetry
w.r.t.\ the common orthogonal bisector straight lines 
of the opposite arc-sides. 

We begin with the case of central symmetry. Then the
opposite arc-sides of 
$(\varphi K) \cap (\psi L)$ on the paracycle ${\text{bd}}\,
(\varphi K)$
are centrally symmetric images of
each other. Then the paracycles containing these arc-sides
have different infinite points (by the collinear model),
hence are different. However, both are equal to
${\text{bd}}\,(\varphi K)$, a contradiction.

We continue with the case of axial symmetry 
w.r.t.\ the common orthogonal bisector straight lines 
of the opposite arc-sides, lying on
bd\,$(\varphi K)$. However, a common orthogonal bisector straight line
to the two opposite (thus disjoint)
arc-sides of $(\varphi K) \cap (\psi L)$, lying on
bd\,$(\varphi K)$, cannot
exist. Namely, such orthogonal bisectors are different and
parallel, so they cannot coincide.

We continue with the case of axial symmetry 
w.r.t.\ the common orthogonal

\newpage

\noindent
bisector straight lines 
of the opposite arc-sides, lying on
$\psi L_1$ and $\psi L_2$. However, a
common orthogonal straight line
to $L_1$ and $L_2$ is a common orthogonal straight line
to their base
lines as well. These base lines are parallel, and by
\thetag{5.5.1} are different, hence admit no common
orthogonal straight line, which is a contradiction. This shows
that our indirect hypothesis \thetag{5.5.1} was false.
This ends the proof.
$\blacksquare $
\enddemo


{\it{Proof of Theorem}} 5, {\bf{continuation.}} {\bf{4.}}
Last, let both $K$ and $L$ be bounded
by finitely many
hypercycles and straight lines. 
Then, by the hypothesis made in the beginning of {\bf{3}} of
the proof of Theorem 5,
either $K$, or $L$ has at least two boundary components. 

We will show that the graphs of $\varphi K$ and $\psi L$ 
must have only a few edges, and we
will clarify the structure of these graphs, till we will obtain that we must
have case (E) of (3) of our theorem.

We will make the following case distinction, for the graphs
of $\varphi K$ and $\psi L$. 
$$
\cases
{\text{(1) Both graphs contain a pair of edges with at least
one common end-point.}}
\\
{\text{(2) One graph contains a pair of edges with at least
one common endpoint,}}
\\
{\text{but the other graph has only vertex-disjoint edges,
and the number of}}
\\
{\text{these edges is at least two.}}
\\
{\text{(3) One graph contains a pair of edges with at least
one common endpoint,}}
\\
{\text{and the other graph has a single edge.}}
\\
{\text{(4) Both graphs have only vertex-disjoint edges, and
at least one of the}}
\\
{\text{graphs contains at least two edges.}}
\endcases
\tag 5.2
$$
These cases are exhaustive, and mutually exclusive. They will
be treated successively
in the following Lemmas 5.7--5.10.
Still we
note that the finiteness hypothesis of Theorem 5 will not be
needed in Lemmas 5.7--5.9, therefore it will not be
included in their hypotheses.

{\bf{5.}}
In the proof of Lemmas 5.7--5.10 we will need ``generic small
perturbations of orientation preserving congruences of $H^2$''.
For this to make sense, we have to give a topology on the
set of orientation preserving congruences of $H^2$.
(Below $\varphi _0, \varphi $ will be orientation preserving
congruences of $H^2$.)
This will
be the {\it{topology of uniform convergence on compacta in
$H^2$}}. This has as (open)
neighbourhood base at some $\varphi _0$
the following: 
$$
\left\{ 
\{ \varphi \mid x \in K \Longrightarrow d(\varphi _0 x,
\varphi x) < \varepsilon \}
\mid K \subset H^2 {\text{ is compact, and }} \varepsilon > 0
\right\} .
\tag 5.3
$$
(Since such a $K$ is bounded in the metric of $H^2$, rather
than all such $K$, we may take only closed balls of some
fixed centre and integer radii.) With this topology the maps
$\varphi \mapsto \varphi x$, where $x \in H^2$, are
continuous.

By \cite{2}, 
Theorem 2,
the set of all congruences of $H^2$ (orientation preserving and
reversing), with the topology defined by the same formula
\thetag{5.3},
is a topological group (i.e., the operations $(\varphi _1,
\varphi _2) \mapsto \varphi _1 \varphi _2$ and $\varphi
\mapsto \varphi ^{-1}$ are continuous). This implies the same
for its subgroup consisting of all orientation-preserving
congru-

\newpage

\noindent
ences of $H^2$.

$$
{\text{Denote by }} {\text{Isom}}^+(H^2) {\text{ the set of
orientation-preserving congruences of }} H^2 .
\tag 5.4
$$
Now recall that in the conformal model, as a subspace of
${\Bbb{C}}$, the set ${\text{Isom}}^+(H^2)$ consists of 
the maps of the following form:
$$
\cases
z \mapsto f_{c, \alpha } (z) :=
c (z - \alpha )/(1 - \overline {\alpha } z) ,
{\text{ where }} | z | < 1 , 
\\
{\text{and }} c \in \{ z \in {\Bbb{C}}
\mid | z | = 1 \} = S^1 , {\text{ and }}
\alpha \in {\text{int}}\, B^2 \subset {\Bbb{C}}
\endcases
\tag 5.5
$$
(cf. \cite{8}).

This is a real analytic manifold, with
real analytic coordinates
$c$, $a := {\text{Re}}\, \alpha $ and $b :=
{\text{Im}}\, \alpha $. 
The usual formulas for
composition and inverse make this real analytic
manifold a topological group, where the composition and the
inverse are real
analytic (a real analytic Lie group with these
coordinates). Moreover, the action $(\varphi , z) \mapsto
\varphi z$ is real analytic (in
$c$, $a$, $b$,
${\text{Re}}\, z $, ${\text{Im}}\, z $),
and for given $c$ and $\alpha $
it is complex analytic in $z$.
Actually, rather than $z \in {\text{int}}\, B^2$, we may take
$z \in B^2$, and then the composition, inverse, and the action
still are continuous (actually real analytic).

In the following Lemma 5.6
we will show that the above geometric
definition of convergence on ${\text{Isom}}^+(H^2)$ is
equivalent to a definition using the analytic functions from
\thetag{5.5}. This will be necessary in Lemma 5.8.
This statement is
surely known, but we could not locate a proof for it.
Therefore we give its simple proof.


\proclaim{Lemma 5.6}
A net $f_{c(\gamma ), \alpha (\gamma )} (z)$
({\rm{defined in \thetag{5.5}}})
converges to $f_{c, \alpha } (z)$ 
in ${\text{\rm{Isom}}}^+$
\newline
$ (H^2)$
in the
topology of uniform convergence on compacta in $H^2$, if and
only if $\left( c(\gamma ),
\alpha (\gamma ) \right) \to
(c, \alpha )$ in $S^1 \times {\text{\rm{int}}}\,B^2$.
\endproclaim


\demo{Proof}
{\bf{1.}}
First we show that the map $S^1 \times {\text{int}}\, B^2 \ni
(c, \alpha ) \mapsto f_{c, \alpha } (z) \in
{\text{Isom}}^+(H^2)$ (with $f_{c, \alpha } (z)$ from
\thetag{5.5}), 
with ${\text{Isom}}^+(H^2)$ endowed with
the topology of uniform
convergence on compacta in $H^2$, is continuous. We write
$a := {\text{Re}}\, \alpha $ and $b :=
{\text{Im}}\, \alpha $.

Since $S^1 \times {\text{int}}\, B^2$ is a metric space,
therefore it suffices to consider only sequences
$\left( c(n), \alpha (n) \right) \in
S^1 \times {\text{int}}\, B^2$,
rather than general nets.
Observe that $S^1 \times {\text{int}}\, B^2 \ni
\left( c(n), \alpha (n) \right) , (c, \alpha )$
and $\left( c(n), \alpha (n) \right) \to (c, \alpha )$
imply
$\alpha (n) , \alpha \in rB^2$ for some $r \in (0,1)$. Hence,
letting $d_{\Bbb{C}}(\cdot , \cdot )$ being the distance in
${\Bbb{C}}$, we have
$f_{c(n), \alpha (n)} (z) \to f_{c, \alpha } (z)$ in
${\Bbb{C}}$, i.e.,
$d_{\Bbb{C}}
\left( f_{c(n), \alpha (n)} (z) , f_{c, \alpha } (z)
\right) \to 0$,
pointwise for $z \in B^2$, but even uniformly on $B^2$.
Namely,
$| (\partial /\partial c) f_{c, \alpha } (z) |,
| (\partial /\partial a) f_{c, \alpha } (z) |,
| (\partial /\partial b) f_{c, \alpha } (z) |
\le {\text{const}} / (1 - | \alpha |)^2 \le {\text{const}}
/ (1-r)^2$,
uniformly for all $z \in B^2$.
Therefore the same
unifom convergence holds for any
compact subset $C \subset {\text{int}}\, B^2$.
However, by Part I (i.e., \cite{5}), 
(0.5), for distinct points of $C$,
the distance $d(\cdot , \cdot )$
in $H^2$, and the
distance $d_{\Bbb{C}}(\cdot , \cdot )$ in the conformal
model as a
subspace of ${\Bbb{C}}$,
have a 
quotient bounded from below and above. Hence \thetag{5.3}
holds for $d(\cdot , \cdot )$, since
from above it holds for the distance
$d_{\Bbb{C}}(\cdot , \cdot )$. This shows that our map from
$S^1 \times {\text{int}}\, B^2$ to ${\text{Isom}}^+(H^2)$,
endowed with the topology of uniform
convergence on compacta in $H^2$, is continuous.

{\bf{2.}}
Second we show that also conversely, convergence of
a net $f_{c(\gamma ), \alpha (\gamma )}$ to $f_{c, \alpha }$,
in the topology of uniform
convergence on compacta in $H^2$, implies
$( c(\gamma ), \alpha (\gamma ) ) \to $

\newpage

\noindent
$(c, \alpha ) $ in $S^1
\times {\text{int}}\, B^2$.

Let us choose the compact set $\{ 0 \} \subset
{\text{int}}\, B^2$. Then
$$
- c(\gamma ) \alpha (\gamma ) =
f_{c(\gamma ), \alpha (\gamma )} (0) \to
f_{c, \alpha } (0) = - c \alpha .
\tag 5.6.1
$$
Now observe that also 
$f_{c(\gamma ), \alpha (\gamma )} ^{-1} \to
f_{c, \alpha }^{-1}$, by the
topological group property of ${\text{Isom}}^+(H^2)$. Let us
choose
the compact set $\{ 0, w \} \subset {\text{int}}\, B^2$,
where $w \ne 0$.
Then 
$$
\alpha (\gamma ) =
f_{c(\gamma ), \alpha (\gamma )} ^{-1} (0) \to 
f_{c, \alpha } ^{-1} (0) = \alpha ,
\tag 5.6.2
$$
and
$$
\left( w + c (\gamma ) \alpha (\gamma ) \right) /
( c (\gamma ) + {\overline{\alpha (\gamma )}} w ) )
= f_{c(\gamma ), \alpha (\gamma )} ^{-1} (w)
\to f_{c, \alpha } ^{-1} (w) =
(w + c \alpha ) / (c + {\overline{\alpha }} w) .
\tag 5.6.3
$$

Then, unless  $\alpha = 0$, we have by \thetag{5.6.1} and
\thetag{5.6.2}, choosing $\gamma $
sufficiently large, that
$$
c(\gamma ) =  c(\gamma ) \alpha (\gamma ) / \alpha (\gamma )
\to c \alpha / \alpha = c .
\tag 5.6.4
$$

There remains the case $\alpha = 0$. In this case,
taking in consideration $| c | = 1 > | w | > 0$, 
we can rewrite \thetag{5.6.3} as
$$
\left( w + c (\gamma ) \alpha (\gamma ) \right) /
( c (\gamma ) + {\overline{\alpha (\gamma )}} w ) ) \to 
( w + c \alpha ) / ( c + {\overline{\alpha }} w ) = w / c .
\tag 5.6.5 
$$
If $c (\gamma )$ had not converged to $c$, then by
compactness of $S^1$, for some subnet it would converge to
some
$c^* \in S^1 \setminus \{ c \} $, and then the limit in
\thetag{5.6.5} for this subnet would be $w/c^* \ne w/c$, a
contradiction.

This ends the proof of the lemma.
$\blacksquare $
\enddemo


\proclaim{Lemma 5.7}
Assume \thetag{*} with $d = 2$ and $X = H^2$.
Suppose {\rm{(1)}} of Theorem {\rm{5}}.
Let the connected boundary components of
both $K$ and $L$ be
hypercycles or straight lines. 
Moreover, suppose {\thetag{5.2}}, case {\rm{(1)}} from above.
Then {\rm{(E)}} of {\rm{(3)}} of Theorem {\rm{5}} holds.
\endproclaim


\demo{Proof}
{\bf{1.}}
By \thetag{5.2}, (1), each of the graphs of $\varphi K$ and $\psi L$ 
contains a path of length
$2$ or a $2$-cycle.
The corresponding boundary components of $\varphi K$ and $\psi L$ are
denoted by $\varphi K_1,\varphi K_2$, and $\psi L_1,\psi
L_2$, with $\varphi K_2$ following $\varphi K_1$ 
on bd\,$(\varphi K)$, and $\psi L_2$ following $\psi L_1$ on
bd\,$(\psi L)$,
according to the positive orientation. (If one of the graphs is a $2$-cycle, then
this does not work. Then we fix some
notation, cf.\ below.)

{\it{We use the conformal model.}}
Recall that any two or three distinct points on the boundary
of the model can be
taken by (the extension of) an orientation-preserving
congruence to any other two or three distinct 
boundary points of the model, of the same orientation.
Therefore we may suppose the following. The considered
common vertex of $\varphi K_1$ and $\varphi K_2$ is
$1 \in {\Bbb{C}}$. If their other vertices coincide, let them
be $-1$. If their other vertices are different, let them be
very close, in $B^2$,
to $-1$. Hence all other
boundary components of $\varphi K$ (if any)
are very close in $B^2$,
to $-1$, as well. Moreover,
the considered common vertex of $\psi L_1$ and $\psi L_2$ is
$i \in {\Bbb{C}}$. If their other

\newpage

\noindent
vertices coincide, let them
be $-i$. If their other vertices are
different, let them be
very close, in $B^2$,
to $-i$. Hence all other
boundary components of $\psi L$ (if any)
are very close in $B^2$,
to $-i$, as well. 

Further,
$\varphi K_2$ follows $\varphi K_1$
on bd\,$(\varphi K)$
at $1$ in the positive sense, and similarly for $\psi L_2,\psi
L_1$ at $i$. Then we have that the distance of $0$ to
$(\varphi K) \cap (\psi L)$ is small (possibly $0$). Moreover,
$(\varphi K) \cap (\psi L)$ is a convex arc-quadrangle,
bounded by arcs of
$\varphi K_1, \psi L_1, \varphi K_2, \psi L_2$, in this cyclic
order, in the positive sense. In fact, all other boundary
components, both of $\varphi K$ and of $\psi L$, are, in $B^2$,
very close to the boundary of the model, hence
cannot cut off parts of the above arc-quadrangle,
which arc-quadrangle is not close to the boundary of the
model.

Thus $(\varphi K) \cap (\psi L)$ is a compact arc-quadrangle.
Its vertices are its only non-smooth points, hence they are
preserved by any congruence admitted by this arc-quadrangle.
Analogously as in the proof of Lemma 5.5, by
Part I (i.e., \cite{5}), 
\thetag{2}, 
\thetag{3} 
and \thetag{4} 
we have that
the possible
non-trivial congruences admitted by this arc-quadrangle
are the following: 
\newline
(a)
two fourfold rotations, and one central symmetry
with a unique centre of
symmetry, namely the midpoint of any
of the two diagonals (which therefore coincide); 
\newline
(b)
axial symmetries w.r.t.\ diagonals (which therefore
bisect the angles at their endpoints);
\newline
(c)
axial symmetries w.r.t.\ the common orthogonal bisector
straight lines of two opposite arc-sides (which are also common
orthogonal straight lines to the base-lines of the hypercycles
containing these arc-sides), which simultaneosly
interchange the remaining two opposite arc-sides (hence also
the hypercycles containing them, and also their base-lines).

If we have a congruence which is a fourfold
rotation, then we
also have a congruence which
is a central symmetry. Hence we need
not exclude the case of a fourfold rotation,
exclusion of a central symmetry will suffice.

Observe that a non-trivial congruence admitted by
$(\varphi K) \cap (\psi L)$
is a non-trivial congruence admitted by
$(\varphi K_1) \cup (\varphi K_2)
\cup (\psi L_1) \cup (\psi L_2)$ as well.

{\bf{2.}}
Next we investigate, for later use, two special cases, and we
prove for them (E) of (3) of our theorem.

First suppose that both $K$ and $L$ are parallel domains of
some straight lines. Thus the (positive)
curvatures of $K_1$ and $K_2$
are equal, and also the (positive)
curvatures of $L_1$ and $L_2$ are equal. We are going to show
that also the curvatures of $K_i$ and $L_i$ coincide, i.e.,
that (E) of (3) of our theorem holds. {\it{Suppose, for
contradiction,
that, e.g., the (positive) curvature of $K_i$ is less than
the (positive) curvature of $L_i$.}}

Observe that $K$, or $L$, 
is symmetrical w.r.t.\ any straight line
orthogonal to the common base line of the
$K_i$'s, or of the $L_i$'s,
resp. 
Hence any two congruent copies 
of $K$ (and of $L$)
are directly congruent if and only if they are
indirectly congruent. (Thus $K$ and $L$ do not have
any specific orientation.) {\it{Therefore we need not care
about
orientation preserving/reversing property of $\varphi $ and
$\psi $, it suffices to derive a contradiction for any 
$\varphi $ and $\psi $.}}
(Recall that in \thetag{*} we had
$\varphi , \psi $ as direct congruences.)

Now we choose a new pair $\varphi , \psi $ (different from
that in {\bf{1}}). This we make so that the 
infinite points
of each of $\varphi K_1,\varphi K_2,\psi L_1,\psi L_2$ are
$\pm 1$, and $\varphi K_1,\psi L_1$ are in the

\newpage

\noindent
open lower
half-plane, and $\varphi K_2,\psi L_2$ are in the open upper
half-plane. Then we
rotate $\psi L$ about the infinite point $1$, 
counterclockwise, by a
rotation of a small measure, obtaining a new $\psi 'L$.
Then, in the conformal model, 
the Euclidean 
tangents of the closures in $B^2$ of
$\psi ' L_1$ and $\psi ' L_2$, at $1$, are the
same as those for $\psi L_1$ and $\psi L_2$.
Therefore, in the new position, 
$(\varphi K) \cap (\psi ' L)$ is an arc-triangle, bounded by
$\varphi K_2,\psi ' L_2,\varphi K_1$, in this positive
order.
This has a
unique infinite point, namely $1$,
and does not contain a paracircle with
this infinite point.
Hence by Part I (i.e., \cite{5}), 
\thetag{2} 
and \thetag{3}, 
a non-trivial congruence, admitted by it, is 
an axial symmetry w.r.t.\ a unique
axis passing through $1$.

Observe that this axial symmetry of $(\varphi K) \cap
(\psi ' L)$
is also an axial symmetry of $[(\varphi K) \cap (\psi ' L)]
\cap P = (\varphi K) \cap P$, for $P$ a paracircle with
centre $1$, which in $B^2$
is a sufficiently small circle
(satisfying $P \cap (\varphi K) \cap
(\psi ' L_2) = \emptyset $).
We have
symmetry of $(\varphi K) \cap P$,
w.r.t.\
the common base line of $\varphi K_1$ and $\varphi K_2$,
i.e., the part of the
real axis in the model. Even this is the unique axial symmetry
of $(\varphi K) \cap P$, with axis of symmetry having $1$ as
an infinite point.

Then $(\varphi K) \cap (\psi ' L)$ and $(\varphi K) \cap P$
coincide in some neighbourhood of $1 \in {\Bbb{C}}$.
Hence also the unique axis of symmetry of
$(\varphi K) \cap (\psi ' L)$,
having $1$ as an infinite point, 
is the part of the real axis in the model. Then
this axial symmetry should preserve the arc-side of
$(\varphi K) \cap (\psi ' L)$ on $\psi ' L_2$.
Therefore the transversal intersection of  
the real axis with $\psi ' L_2$ in the conformal model 
should be an
orthogonal intersection.
However then in $B^2$,
the real axis should intersect the
closure in $B^2$ of the hypercycle
$\psi ' L_2$, at $1$, at a right angle as well, which is
impossible.
Hence, in the new position,
$(\varphi K) \cap (\psi ' L)$ admits no non-trivial
congruence, which is 
a contradiction. Therefore our indirect hypotesis that, e.g.,
the
curvature of $\varphi K_i$ was less than that of $\psi L_i$,
was false.
This proves (E) of (3) of our theorem.

{\bf{3.}}
Next, as a second special case, 
we suppose that $K$ and $L$ are congruent, and both are bounded by two hypercycles with common base-line.

By a suitable notation, 
we have that the curvatures of $K_1$ and $ 
L_1$
are equal, and also
that the curvatures of $K_2$ and $L_2$
are equal.
We are going to show that 
$K_1$ and $K_2$, and then also $L_1$ and
$L_2$, each has the same
curvature, i.e., that (E) of (3) of our theorem holds. 
{\it{Suppose, for contradiction,
that, e.g., the (non-negative) curvature of $K_1$ is less than
the (positive) curvature of $K_2$.}}

Again, like in {\bf{2}}, $K$, or $L$ 
is symmetrical w.r.t.\ any straight line
orthogonal to the common base line of the
$K_i$'s, or of the $L_i$'s, resp. Therefore
we derive a contradiction for any 
$\varphi $ and $\psi $, independently of their
orientation preserving/reversing property.

Now we choose a new pair $\varphi , \psi $ (different from
that in {\bf{1}}).
Let us fix $\varphi K$ so that its points at infinity are $\pm 1$, and $\varphi K_1$
lies in the closed
lower half-plane, and $\varphi K_2$ lies in the open upper
half-plane (observe that possibly $\varphi K_1$ is a straight
line -- but $\varphi K_2$ is not).
Let us
obtain $\psi L$ by rotating $\varphi K$ about the infinite
point
$1$ in positive sense a bit,
with the image of $\varphi K_i$ being $\psi L_i$. Then
$(\varphi K) \cap (\psi L)$ is
bounded by two arcs,
one lying on $\varphi K_1$, the other one lying on
$\psi L_2$. Then $(\varphi K) \cap (\psi L)$ has one
non-smooth point, at $(\varphi K_1) \cap (\psi L_2)$. 
Hence a non-trivial congruence admitted by
$(\varphi K) \cap (\psi L)$
is an axial symmetry, w.r.t.\ the angle bisector of the angle of

\newpage

\noindent
$(\varphi K) \cap (\psi L)$ at this point. However, this
axial symmetry
interchanges the two portions of the boundary, lying on
$\varphi K_1$ and $\psi L_2$.
This contradicts our indirect hypothesis
that the curvature of $\varphi K_1$
is less than the curvature of $\varphi K_2$,
i.e., the curvature of $\psi L_2$.
This contradiction shows (E) of (3) of our theorem.

{\bf{4.}}
We turn to the general case.
We begin with case (a), i.e., when $(\varphi K) \cap (\psi L)$
has a central symmetry. In this case
we will show (E) of (3) of our theorem.

Then $\varphi K_1$ and $\varphi K_2$ have two common infinite
points (images of each other by this symmetry), and the same statement holds
for $\psi L_1$ and $\psi L_2$, hence the graphs of both 
$\varphi K$ and $\psi L$ 
are $2$-cycles. Clearly then the curvatures of $\varphi K_1$ and $\varphi K_2$, as well as 
those of $\psi L_1$ and $\psi L_2$, are equal, and are
positive. That is, we have the case settled in {\bf{2}}, and
hence (E) of (3) of our theorem holds.

{\bf{5.}}
Each of the cases (b) and (c) in {\bf{1}}
describes two axial symmetries. Therefore, unless we have (E)
of (3) of our theorem,
we have to exclude simultaneously four possible axial
symmetries.
First observe that our question about $(\varphi K) \cap
(\psi L)$, for arbitrary $\varphi $ and $\psi $, is equivalent
to the same question with $\psi = $ identity.
$$
{\text{So let }} \psi = {\text{identity.
(This hypothesis will be applied in the remainder
of this proof.)}}
\tag 5.7.1
$$

Then we
are going to show that unless we have (E) of (3) of our
theorem, we have the following.
A generic arbitrary small perturbation $\varphi '$
of $\varphi $ (in the topology on ${\text{Isom}}^+(H^2)$,
given by \thetag{5.3}, cf.\ also \thetag{5.4})
simultaneously destroys all four possible axial
symmetries from above.

Actually, in {\bf{6}}, in case $(\alpha )$ (including the
case of symmetry w.r.t.\ $l_1$), (E) of (3) of our
theorem will be proved directly.

In the
remaining three cases, i.e., symmetry w.r.t.\ $l_2$, in case
$(\beta )$, and symmetries w.r.t.\ $l_3$ and $l_4$,
we show the following. The $i$'th
possible congruence (for $1 \le i \le 3$)
admitted by $(\varphi ' K) \cap (\psi L) = (\varphi ' K)
\cap L$
does not exist for $\varphi '$ in a set of the form
$U_i(\varphi ) \setminus N_i(\varphi )$. Here  $U_i(\varphi )$
is an open
neighbourhood of $\varphi $, and $N_i(\varphi )$ is a
relatively
nowhere dense, relatively closed subset of $U_i(\varphi )$.

By elements of topology, this can be rewritten as follows. 
The $i$'th possible congruence ($1 \le i \le 3$) admitted by
$(\varphi ' K) \cap L$
does not exist for $\varphi '$ in a set of the form
$U_i(\varphi ) \setminus M_i$. Here  $U_i(\varphi )$
is as above, and $M_i$ is a
nowhere dense, closed subset of ${\text{Isom}}^+(H^2)$.
Then none of
the investigated three
axial symmetries hold on $\cap _{1 \le i \le 3}
[U_i(\varphi ) \setminus M_i] =
[\cap _{1 \le i \le 3}
U_i(\varphi )] \setminus [
\cup _{1 \le i \le 3} M_i] $.
Recall that a finite union of nowhere dense closed subsets
is also 
nowhere dense and closed. Hence none of the investigated three
axial symmetries hold on a set of the form
$[\cap _{1 \le i \le 3} U_i(\varphi )] \setminus N(\varphi )$,
where $N(\varphi )$
is a relatively nowhere dense, relatively closed subset of 
$\cap _{1 \le i \le 3} U_i(\varphi )$.
Therefore,
for a generic small perturbation
$\varphi '$ of $\varphi $, the set $
(\varphi 'K) \cap L$ will have none of the investigated three
possible axial symmetries. 

{\bf{6.}}
We turn to the case when $(\varphi K) \cap L$
has an 
axial symmetry w.r.t.\ a straight line $l$ spanned by
some diagonal. Then symmetry w.r.t.\ $l$
interchanges $(\varphi K_1) \cup (\varphi K_2)$ and
$L_1 \cup L_2$. Therefore these last two sets are also
axially symmetric images of each other. 

Hence the graphs of $\varphi K$ and $L$ either 
\newline
($\alpha $)
are simultaneously $2$-cycles, with edges
$\varphi K_1,\varphi K_2$, and $L_1,L_2$, or

\newpage

\noindent
($\beta $)
simultaneously contain paths of length $2$, namely
$(\varphi K_1)(\varphi K_2)$ and $L_1 L_2$,
resp. 

Observe that if the symmetry w.r.t.\ $l$ interchanges two
arc-sides
of our arc-quadrangle, then it interchanges
the entire hypercycles
containing them, and also their base-lines. Hence
in case ($\alpha $) we have that ${\text{bd}}\,(\varphi K)
= (\varphi K_1) \cup (\varphi K_2)$ and ${\text{bd}}\,L
= L_1 \cup L_2$ are symmetric images of
each other w.r.t.\ $l$. Hence also $\varphi K$ and $L$ 
are symmetric images of each other w.r.t.\ $l$.
However then we have the situation discussed in {\bf{3}}, and
hence (E) of (3) of our theorem holds.

Recall that in {\bf{1}}
of this proof,
$\varphi K$ and $\psi L = L$ were chosen as follows. The considered
common infinite point of $\varphi K_1$ and $\varphi K_2$ was
$1$, the considered
common infinite point of $L_1$ and $L_2$ was $i$, and
the other end-points of
$\varphi K_1$ and $\varphi K_2$ (possibly coinciding)
were close to $-1$, and the
other end-points of
$L_1$ and $L_2$ (possibly coinciding) were close to $-i$.

By {\bf{1}} of this proof $(\varphi K) \cap L$ is a
compact arc-quadrangle, such that
the distance of $0$ to it is small (possibly is $0$). 
It has two diagonals. One of these diagonals
connects the intersection
points $(\varphi K_1) \cap L_1$ and 
$(\varphi K_2) \cap L_2$, and spans the straight line $l_1$. 
The other diagonal connects the
intersection points $(\varphi K_2) \cap L_1$ and 
$(\varphi K_1) \cap L_2$, and spans the straight line $l_2$.
(It will be convenient to identify a point $x$ with the
one-point set $\{ x \} $.)

{\bf{A.}}
If $l_1$ is an axis of symmetry of
$(\varphi K) \cap L$, 
then also the entire hypercycles $\varphi K_1$ and $L_1$,
as well as $\varphi K_2$ and $L_2$ are symmetric images of
each other w.r.t.\ $l_1$.
Thus the common infinite point $1$ of $\varphi K_1$ and
$\varphi K_2$ has as symmetric image w.r.t.\ $l_1$
a common infinite point of
$L_1$ and $L_2$, close to $-i$ in ${\Bbb{C}}$, hence
different from their common infinite point 
$i$. Hence $L$ is bounded by $L_1$ and $L_2$. Similarly,
$\varphi K$ is bounded by $\varphi K_1$ and $\varphi K_2$.
Then $\varphi K_1$ and $L_1$
are symmetric images of each other w.r.t.\
$l_1$. Since both $\varphi K$ and $L$ are axially
symmetric, therefore by {\bf{2}}
$\varphi K_1$ and $L_1$ are also
directly congruent.
Thus we have case $(\alpha )$, which above was shown to
imply (E) of (3) of our theorem.

{\bf{B.}}
If $l_2$ is an axis of symmetry of
$(\varphi K) \cap L$, 
then also the entire hypercycles $\varphi K_1$ and $L_2$,
as well as $\varphi K_2$ and $L_1$, are symmetric images of
each other w.r.t.\ $l_2$. In case $(\alpha )$ again
we have (E) of (3) of our theorem.

Therefore suppose case $(\beta )$.
First we show that some generic small perturbation $\varphi '$
has the property, that the analogously defined $l_2'$
(i.e., the straight line
spanned by the vertices
$(\varphi ' K_2) \cap L_1$ and $(\varphi ' K_1) \cap L_2$)
is not
an axis of symmetry of $(\varphi ' K) \cap L$. More exactly,
we will show that $l_2'$ is not
an angle bisector of the angle
of our perturbed arc-quadrangle at its vertex at
$(\varphi ' K_2) \cap L_1$. 

The arc-sides of our perturbed
arc-quadrangle on $\varphi ' K_2$ and on $L_1$ (which
by the conformal model intersect
transversally at the vertex $(\varphi ' K_2) \cap L_1$,
hence the convex angle
determined by them is unique) determine uniquely the
angle-bisector $a'$ of our perturbed
arc-quadrangle at this vertex.
(For $\varphi K_2$ and/or $L_1$ straight line/s
we have to consider also on which side of this/these straight
lines our arc-quadrangle lies.)

Now we apply to $\varphi K$ a
small non-zero translation along the base line 
of $\varphi K_2$,
obtaining $\varphi ' K$. This translation
preserves $\varphi K_2$ setwise -- hence $a'$ does not change
by this translation. Moreover,
it rotates $\varphi K_1$ non-trivially
about the common infinite point $1$ of $\varphi K_2$ and
$\varphi K_1$. Then $(\varphi ' K_1) \cap L_2 \not\subset
a'$. Namely else $\emptyset \ne (\varphi ' K_1) \cap (L_2 \cap
a') \subset (\varphi ' K_1) \cap (\varphi K_1) = \emptyset $.
(Since different rotated copies of $\varphi K_1$,

\newpage

\noindent
about $1$, are disjoint.)
In other words, the perturbed diagonal and
angle-bisector are different, hence we do not have axial
symmetry of $(\varphi ' K) \cap L$
w.r.t.\ the perturbed straight line
$l_2'$.
Hence the investigated symmetry does not exist
for the perturbed arc-quadrangle.

Let us choose a small connected open neighbourhood
$U(\varphi )$ of $\varphi $.
We may assume that the infinite points of $\varphi ' K_1,
\varphi ' K_2, L_1, L_2$ (altogether six infinite points)
are different, and cyclically follow each other in the same
positive sense as those of $\varphi K_1,
\varphi K_2, L_1, L_2$.
As shown above, for some
$\varphi ' \in U(\varphi )$ we have that $l_2'$
does not bisect the investigated angle of our
arc-quadrangle. Therefore the
equality of the angles into which $l_2'$
divides our angle, for some
$\varphi ' \in U(\varphi )$,
is expressed by an analytic equation for $\varphi '$, which is
not an identity. Hence,
by the principle of analytic continuation, this equality is
satisfied only for $\varphi '$ in a relatively nowhere
dense relatively closed subset of $U(\varphi )$.
Hence, all
such $(\varphi ')$'s, with $(\varphi ' K) \cap L$
satisfying the symmetry property in {\bf{B}}, in case
$(\beta )$, form
a nowhere dense closed set in $U(\varphi )$.

{\bf{7.}} 
There remained the case, when
$(\varphi K) \cap L$ has 
an axial symmetry w.r.t.\ a
common orthogonal bisector straight line of two opposite
edges. If these edges lie on
$\varphi K_1$ and $\varphi K_2$, then let this line be $l_3$. 
Hence $l_3$ is a common orthogonal straight line to the
base-lines of $\varphi K_1$ and $\varphi K_2$.
Moreover, this axial symmetry interchanges the arc-sides
on $L_1$ and $L_2$, hence also the base-lines of $L_1$ and
$L_2$ (case (c) in {\bf{1}}).

Observe that a common orthogonal straight line
to the parallel base-lines of $\varphi K_1$ and $\varphi K_2$ 
exists only if these base-lines coincide. In this case
the graph of $\varphi K$ is a $2$-cycle, and
$\varphi K$ is
bounded just by $\varphi K_1$ and $\varphi K_2$. Let the other common infinite point of
$\varphi K_1$ and $\varphi K_2$ be $-1$
(and let $\psi L = L$ be as described in the third paragraph
of {\bf{1}} of the proof of this lemma).
The axis $l_3$ of our
symmetry is a straight line orthogonal to 
this common 
base line, i.e., to the real axis.

On the other hand, $l_3$ is the unique axis of symmetry
interchanging the parallel base-lines of
$L_1$ and $L_2$, hence contains $i$.
(Since our axial symmetry w.r.t.\ $l_3$, admitted by
$(\varphi K) \cap L$, yields
an axial symmetry, admitted by each of $\varphi K_1$,
$\varphi K_2$, and $L_1 \cup L_2$.)

Thus $l_3$
is the imaginary axis. Hence the common base-line of 
$\varphi K_1$ and $\varphi K_2$,
and the unique axis of symmetry of
$L_1 \cup L_2$, interchanging $L_1$ and $L_2$,
are orthogonal.
However, all such $\varphi $'s form
a closed set in our topology
(given by \thetag{5.3}), which is moreover nowhere dense.
Namely,
a small rotation of $\varphi K$, about the
point of intersection of the above two lines,
destroying this orthogonality,
is possible.

For the case of axial symmetry, w.r.t.\ the common orthogonal
bisector straight line $l_4$ to the opposite edges on $L_1$ and 
$L_2$, we proceed analogously. (Observe that the roles of
$\varphi $ and $\psi $ can be changed: we may fix $\varphi $
as identity, and then $\psi $ will vary.)

{\bf{8.}}
Hence, by {\bf{5}} of this proof, 
we obtain the following. In {\bf{6}}, case
$(\alpha )$ was settled (which
included the case of symmetry w.r.t.\ $l_1$). In the
remaining cases we had the following.
Case
$(\beta )$, with symmetry w.r.t.\ $l_2$, was investigated
in {\bf{6}}.
The cases of
symmetries w.r.t.\ $l_3$  and $l_4$ were investigated
in {\bf{7}}.
In all these remaining three cases
we obtained the following.
For a generic small perturbation $\varphi '$ of
$\varphi $, we had that $(\varphi ' K) \cap L$ had
none of the remaining three axial
symmetries.
Summing up:

\newpage

\noindent
we have proved (E) of (3) of our theorem.
$\blacksquare $
\enddemo


\proclaim{Lemma 5.8}
Assume \thetag{*} with $d = 2$ and $X = H^2$.
Suppose {\rm{(1)}} of Theorem {\rm{5}}.
Let the connected boundary components
both of $K$ and $L$ be
hypercycles or straight lines. 
Moreover, suppose {\thetag{5.2}}, case {\rm{(2)}}.
Then this leads to a contradiction.
\endproclaim


\demo{Proof}
{\bf{1.}}
Let, e.g., the graph of $\varphi K$ 
consist of vertex disjoint edges, whose number is at least $2$. Let us
choose two vertex-disjoint edges of this graph, $\varphi K_1$,
$\varphi K_2$, say. Further, let the graph of $\psi L$ 
contain a path of length $2$ or a $2$-cycle, 
consisting of $\psi L_1$ and $\psi L_2$, 
where $\psi L_2$ follows $\psi L_1$ in the positive orientation
on ${\text{bd}}\,(\psi L)$
(if we have
a $2$-cycle, then their numeration is done in some way).
We are going to show that this case cannot occur.

We fix $\varphi K$ and thus $\varphi K_1$ and $\varphi K_2$,
and will choose $\psi L$ in the
following way. The set $\varphi K$ lies in the closed
convex set bounded by $\varphi K_1$ and
$\varphi K_2$. Then we have relatively open arcs $I_1$
and $I_2$
of the boundary of the (conformal) model, bounded by 
the infinite points of $\varphi K_1$ and $\varphi K_2$,
resp., and
lying outside the closure in the closed unit circle $B^2$
of the closed
convex set bounded by $\varphi K_1$ and
$\varphi K_2$. We choose the (considered)
common infinite point of $\psi L_1$
and $\psi L_2$ close to the midpoint $1$ of $I_1$,
with $\psi L_2$
following $\psi L_1$ there in the positive sense,
and the other infinite
points of $\psi L_1$ and
$\psi L_2$ (possibly coinciding) close to the midpoint of
$I_2$. We suppose that $\varphi K_1$ separates $1$
and $\varphi K_2$.

Then $(\varphi K) \cap (\psi L)$ is contained in 
a {\it{compact arc-quadrangle $Q$, bounded by arcs
lying on $\varphi K_1, \psi L_2, \varphi K_2, \psi L_1$, 
in this order, say}}. 

Observe that all boundary components of $\psi L$, other than
$\psi L_1$ and $\psi L_2$, 
are in the model very close to the boundary of the model, hence
cannot cut off parts of this arc-quadrangle $Q$, which
arc-quadrangle is
not close to this boundary. 
So these boundary components have no arcs on 
bd\,$[ (\varphi K) \cap (\psi L) ] $.
So we need not deal with these 
boundary components.

However, there may exist several boundary components
$\varphi K_i$ of $\varphi K$,
with $i\ne 1,2$,
which cut off parts of this arc-quadrangle, hence have
disjoint non-trivial closed arcs on
bd\,$[ (\varphi K) \cap (\psi L) ] $. However, there
can be only finitely many of these, since each compact set
in $H^2$ meets only finitely many $(\varphi K_i)$'s.

Since we investigate case \thetag{5.2}, (2), we have that the 
$(\varphi K_i)$'s have no common endpoints.
Of course, $\psi L_1$ and $\psi L_2$ have at least one
common endpoint. However, by
con\-struc\-tion, 
neither $\psi L_1$ and any $\varphi K_i$ (including
$\varphi K_1$ and $\varphi K_2$),
nor $\psi L_2$ and any $\varphi K_i$ (including $\varphi K_1$
and $\varphi K_2$), have any
common endpoint.

{\bf{2.}}
We are going 
to show that
$$
\cases
{\text{any non-trivial congruence admitted by }}
(\varphi K) \cap (\psi L) {\text{ is a congruence}}
\\
{\text{admitted by }} Q {\text{ as well; moreover, it is a
congruence mapping each}}
\\
{\text{arc-side of }} Q {\text{ either to itself,
or to the opposite arc-side of }} Q.
\endcases
\tag 5.8.1
$$

We have that bd\,$[ (\varphi K) \cap (\psi L) ] $
consists of arcs, following
each other, in the positive sense, lying on $\varphi K_1,
\psi L_2,\varphi K_{i(1)},\psi L_2,\varphi K_{i(2)},...,
\varphi K_{i(j)},\psi L_2,\varphi K_2,\psi L_1,$
\newline
$\varphi K_{i(j+1)},\psi L_1,\varphi K_{i(j+2)},
...,\varphi K_{i(k)},\psi L_1$, in this order, say. From all of
these arcs
only those lying on $\psi L_1$ and $\psi L_2$ lie on different hypercycles, which have
at least one infinite end-point in common. 

\newpage

\noindent
Let us introduce a symmetric relation $\Cal R$ on the
arc-sides of $(\varphi K) \cap (\psi L)$. For two such
arc-sides $S_1,S_2$ 
we write $S_1 {\Cal R} S_2$, 
if the hypercycles spanned by these sides are
different, and have at least one common end-point. Clearly 
any non-trivial congruence $\chi $
of $(\varphi K) \cap (\psi L)$ preserves this
relation $\Cal R$ (i.e., $S_1 {\Cal R} S_2
\Longleftrightarrow
\chi (S_1) {\Cal R} \chi (S_2)$),
hence also the set of arc-sides ${\Cal S} := 
\{ S_1 \mid \exists S_2
$ such that $S_1 {\Cal R} S_2 \} $. Observe that the relation
$\Cal R$ induces a complete bipartite graph on the vertex set
${\Cal S}$, with
classes ${\Cal L}_i$, for 
$i=1,2$, where ${\Cal L}_i$ is the set of arc-sides
of $(\varphi K) \cap (\psi L)$, lying on $\psi L_i$, for $i=1,2$.

Therefore, each non-trivial congruence admitted by
$(\varphi K) \cap (\psi L)$
preserves the two-element set
$\{ {\Cal L}_1, {\Cal L}_2 \} $. Of course, also
the cyclic order of the arc-sides of $(\varphi K) \cap (\psi L)$
is preserved, 
up to inversion. Let the first end-point of the first arc-side and 
last endpoint of the last arc-side in ${\Cal L}_i$ 
(i.e., lying on $\psi L_i$), in
the positive order,
be $v_{i,1}$ and $v_{i,2}$.
Then 
the set $\{ \{ v_{1,1},v_{1,2} \} , \{ v_{2,1},v_{2,2} \} \} $ is preserved by
each non-trivial congruence admitted by $(\varphi K) \cap (\psi L)$,
as well. So \thetag{5.8.1} is proved.

By {\thetag{5.8.1}}
we need to discuss only the congruences of the
arc-quadrangle
$Q$, more exactly only those of them, which
satisfy the property in the second half of \thetag{5.8.1}.
Therefore, like in the
proof of Lemma 5.7, last paragraph of {\bf{1}}, 
the possible non-trivial congruences admitted by $Q$, to be
investigated
and to be excluded, are central symmetry,
and axial symmetries w.r.t.\ common orthogonal side-bisector
straight lines 
of opposite sides.

Below, {\it{if there is an axial symmetry interchanging
$\psi L_1$
and\,\ $\psi L_2$}} (which is unique, if it exists, and
is determined by $\psi L$), {\it{then
its axis of
symmetry will be denoted by $a(\psi L)$}}.
Like in the proof of Lemma 5.7, {\bf{3}}, {\it{we set $\psi =$
identity}}.

{\bf{3.}}
We begin with the case of central symmetry of $Q$.
The central symmetry interchanges the arc-sides lying on
$\varphi K_1$ and $\varphi K_2$,
hence also $\varphi K_1$ and $\varphi K_2$, hence also
the infinite points of $\varphi K_1$,
and the infinite points of $\varphi K_2$. Thus its
centre must be the intersection of the straight lines connecting the
interchanged end-points. (This determines the interchanged pairs of infinite
points uniquely: these pairs separate each other on the
boundary of the model.)
We denote this centre of symmetry, which is determined
by $\varphi K$, by $c(\varphi K)$.
Also, by central symmetry,
the graph of $\psi L = L$ must be a $2$-cycle, with
$L_1$ and
$L_2$ having the same
curvatures. Then $c(\varphi K)$ lies on
the common base line of $L_1$ and $L_2$, which is
$a(L)$.
(Else the
distances of $c(\varphi K)$ from $L_1$ and $L_2$ would
be different.) Hence
$c(\varphi K)$ {\it{lies on}} $a(L)$.
Clearly the set
of such $\varphi $'s is a nowhere dense closed set.

{\bf{4.}}
We continue with the case of axial symmetry of $Q$
w.r.t.\ the common orthogonal
bisector of the arc-sides lying on $L_1$ and
$L_2$.
Such a common orthogonal
straight
line is orthogonal to the base lines of $L_1$ and
$L_2$ as well, hence it
exists only if the base lines of $L_1$ and $L_2$
coincide (they cannot be
parallel but different), i.e., the graph of $L$ is a
$2$-cycle. 
Then this symmetry interchanges $\varphi K_1$ and $\varphi K_2$,
hence
its axis is the axis of the unique axial
symmetry interchanging 
$\varphi K_1$ and $\varphi K_2$.
{\it{We denote
this axis of symmetry, which is determined
by $\varphi K$, by $a'(\varphi K)$.}}
Then
reflection w.r.t.\ $a'(\varphi K)$
preserves the common base line of
$L_1$ and $L_2$. Therefore this common base line
either coincides with the axis of reflection, or intersects it
orthogonally.
The first case is impossible, since $L_1$ and
$L_2$, and hence
also their common base line is not contained in
the closed convex set bounded by $\varphi K_1$ and

\newpage

\noindent
$\varphi K_2$.
Hence $a'(\varphi K)$ intersects orthogonally
the common base line of $L_1$ and $L_2$, which is
$a(L)$. That is,
$a'(\varphi K)$ {\it{intersects orthogonally}} $a(L)$.
Clearly the set
of such $\varphi $'s is a nowhere dense closed set.

{\bf{5.}}
We continue with the case of axial symmetry of $Q$
w.r.t.\ the common orthogonal
bisector straight line
of the arc-sides lying on $\varphi K_1$ and $\varphi K_2$. This
axis is the unique
straight line orthogonal to $\varphi K_1$ and $\varphi K_2$
(and hence also
to their base lines).
We denote
this axis of symmetry, which is determined
by $\varphi K$, by $a''(\varphi K)$. 
Hence,
$a''(\varphi K)$ {\it{coincides with}} $a(L)$.
Clearly the set
of such $\varphi $'s is a nowhere dense closed set.

{\bf{6.}}
Considering all three possible cases in {\bf{3}}, {\bf{4}},
{\bf{5}}, we apply {\bf{3}} of the proof of Lemma 5.7.
Thus we obtain that for a generic small perturbation of
$\varphi $, we have that $(\varphi K) \cap L$ admits
none of the possible three
non-trivial congruences.
This is a contradiction, which proves the lemma.
$\blacksquare $
\enddemo


\proclaim{Lemma 5.9}
Assume \thetag{*} with $d = 2$ and $X = H^2$.
Suppose {\rm{(1)}} of Theorem {\rm{5}}.
Let
the connected boundary components both of $K$ and $L$ be 
hypercycles or straight lines. 
Moreover, suppose {\thetag{5.2}}, case {\rm{(3)}}.
Then this leads to a contradiction.
\endproclaim


\demo{Proof}
Let, e.g., the graph of $\varphi K$ 
contain a single edge $\varphi K_1$, i.e., this is the unique
boundary component of $\varphi K$,
and let the graph of $\psi L$ contain a
path of length $2$ or a $2$-cycle. 
We are going to show that this is impossible.

Let the graph of $\psi L$ contain two edges
$\psi L_1$ and $\psi L_2$ with a common vertex at $1$. 
Let the other endpoints of $\psi L_1$ and $\psi L_2$ be
close to $-1$.
We consider the conformal model.
We consider $\psi L_1$ and $\psi L_2$ as fixed, and $\varphi K$ 
as being in a small Euclidean neighbourhood of $1$.
Then $\varphi K$ does not intersect 
any other $\psi L_i$, so we may leave them out of
consideration. Moreover, let
$\varphi K_1$ intersect
both $\psi L_1$ and $\psi L_2$. Then $\varphi K$ lies on
that side of $\varphi K_1$, as $1$.

Then $(\varphi K) \cap (\psi L)$ is an arc-triangle, with one
infinite vertex at $1$. Hence any of its non-trivial
congruences is an axial symmetry,
w.r.t.\ an axis passing through $1$ (it cannot be a non-trivial
rotation about the
infinite point $1$, by Part I (i.e., \cite{5}), 
\thetag{3}). 
Moreover,
$\psi L_1$ and $\psi L_2$, as well as
their base lines,
are images of each other by this axial
symmetry. Then
this axis of symmetry, interchanging these base lines, is 
uniquely determined by $\psi L$ (else $(\varphi K) \cap
(\psi L)$ would admit a non-trivial
rotation about the
infinite point $1$, which is excluded by \thetag{3}). 
This axis
intersects 
the arc-side of $(\varphi K) \cap (\psi L)$ on $\varphi K_1$
orthogonally. However, a
small generic
rotation of $\varphi K$, and thus also of $\varphi K_1$,
about the intersection point of the above axis of symmetry
and $\varphi K_1$, destroys orthogonality of this intersection,
hence destroys also this unique congruence.
$\blacksquare $
\enddemo


In the proof of the next lemma, in case \thetag{5.2}, (4),
we will need boundaries of
closed convex sets $K,L$
with interior points,
with finitely many (but at least one) hypercycle or straight
line
boundary components
in $H^2$, and also their boundaries in the
(conformal or col\-li\-ne\-ar) model circle together
with its boundary $S^1$.
The second boundaries of closed convex sets
with interior points
are homeomorphic to $S^1$, and
contain the first ones. In case \thetag{5.2}, (4),
with the finiteness
hypothesis,
their set theoretical difference (second one minus first one)
is the union of at least

\newpage

\noindent
one but finitely many disjoint
closed proper circular subarcs of $S^1$ (in case
\thetag{5.2}, (4)).


\proclaim{Lemma 5.10}
Assume \thetag{*} with $d = 2$ and $X = H^2$.
Suppose {\rm{(1)}} of Theorem {\rm{5}}.
Let the connected boundary components
both of $K$ and $L$ be 
finitely many
hypercycles or straight lines. 
Moreover, suppose {\thetag{5.2}}, case {\rm{(4)}}.
Then this leads to a contradiction.
\endproclaim


\demo{Proof}
By hypothesis, both graphs are finite.
By \thetag{5.2}, (4), both graphs 
contain only vertex-disjoint edges,
and, e.g., the graph of $L$ has at least two edges.
We are going to show that this is impossible.

Let $\psi L_1$ and $\psi L_2$ denote two neighbourly boundary
components of $\psi L$,
with $\psi L_2$ following $\psi L_1$ in the positive
orientation. (That is, 
passing on the boundary of $\psi L$, 
meant in one of the model circles
together with its boundary circle,
from $\psi L_1$
to $\psi L_2$, in the positive sense, there are no other
connected components of 
bd\,$(\psi L)$ between them. Here we use
finiteness of the set of the boundary components of $L$.
Cf.\ also Remark 3.)
Then, denoting by $\psi l_1$ the second 
infinite point of $\psi L_1$, and by
$\psi l_2$ the first infinite point of $\psi L_2$ (both taken
in the positive
orientation), the counterclockwise open arc
${\widehat{(\psi l_1)(\psi l_2)}}$ of $S^1$
contains no infinite point of any boundary component of
$\psi L$.
Then
the base lines of $\psi L_1$ and $\psi L_2$ are symmetric
images of each other w.r.t.\ a unique
axis $a(\psi L)$, determined by $\psi L$. The infinite points
of $a(\psi L)$ are $p(\psi L)$ and $q(\psi L)$, where 
$p(\psi L)$ lies in the  counterclockwise open arc
${\widehat{(\psi l_1)(\psi l_2)}} \subset S^1$. We 
suppose that $a(\psi L)$ {\it{lies on the real axis, and then
$p(\psi L) = 1$ and $q(\psi L) = -1$}}.

Let $\varphi K_1$ be a boundary component of $\varphi K$. Let 
its infinite end-points be $\varphi k_1'$ and $\varphi k_1''$,
following each other in this
order in the positive sense, on the boundary of $\varphi K$,
meant in the
model together with its boundary circle. 
Let us begin with
the position, when $\varphi k_1'=\psi l_2$ and $\varphi k_1''
=\psi l_1$,
and $\varphi K$ lies on the same side of $\varphi K_1$, as
$1$.

Now let us translate $\varphi K$, and thus also
$\varphi K_1$, 
along the real
axis a bit, so that $\varphi k_1'$ and
$\varphi k_1''$ move on $S^1$ in
the direction from $1$ to $-1$, both on
the respective open arcs of $S^1$
from  $1$ to $-1$,
in which they are contained.
For the new
congruent copy of $K$ 
we will not apply a new notation, but
will preserve the old notation $\varphi K$. 
We want to determine the intersection
$(\varphi K) \cap (\psi L)$. 

Let the boundary components of $\varphi K$, or of $\psi L$,
be, in the
positive sense, $\varphi K_1,\varphi K_2,...,\varphi K_n$, or
$\psi L_1,...,\psi L_m$, resp. By \thetag{5.2}, (4),
using the
collinear model, we see that any of $\varphi K$ and $\psi L$
can be obtained from a convex
polygon, with all vertices at infinity
(not counted to the polygon), whose number of
vertices 
is even, in the following way. We have that
$\varphi K$ or $\psi L$ contains the polygon, and 
${\text{bd}}\,(\varphi K)$ or 
${\text{bd}}\,(\psi L)$ is obtained from the boundary of
the polygon the following way. We replace 
each second side of the polygon
by a hypercycle with this base line
(outwards of the polygon),
and remove each other 
second side of the polygon.
(Including the case when this convex polygon
is a $2$-gon,
i.e., a segment.)

Then we may suppose that all boundary components of the new 
$\varphi K$, except $\varphi K_1$ 
(if any), lie strictly on the other side of the 
straight
line $(\psi l_1)(\psi l_2)$, as the new $\varphi K_1$.
(Here we use finiteness of the
set of the boundary components of $K$.
Cf.\

\newpage

\noindent
also Remark 3.)
All these are
boundary components of $(\varphi K) \cap (\psi L)$ as well. 
There is still one boundary
component of $(\varphi K) \cap (\psi L)$. This begins at
$\psi l_2$, then passes on
$\psi L_2$, then passes on $\varphi K_1$, then on some
$\psi L_{i(1)}$, then once more on $\varphi K_1$,
then on some $\psi L_{i(2)}$, ..., then on some
$\psi L_{i(k)}$, then once more on
$\varphi K_1$, then on $\psi L_1$, and ends at $\psi l_1$, with
all these arc-sides being nondegenerate.
(There can be no common arc of $\varphi K_1$ and any $\psi
L_i$. For $\psi L_1$ and $\psi L_2$ this is evident. For
$L_i$, with $i \ge 3$,
this could occur only if 
we had $\varphi K_1 = \psi L_i$, and then $\psi L_i$
would intersect both $\psi L_1$ and $\psi L_2$, a
contradiction.)
Then
$$
\cases
{\text{any non-trivial congruence admitted by }} (\varphi K)
\cap (\psi L) {\text{ preserves}}
\\
{\text{this unique non-smooth boundary component }} C
{\text{ of }} (\varphi K) \cap (\psi L).
\endcases
\tag 5.10.1
$$

The extension
to $B^2$ of a congruence admitted by $C$, in particular,
preserves the pair of its infinite points.
That is, either each of them is mapped to itself, or they are
exchanged by the congruence.
If each of them is mapped to itself
by the congruence, then the semi-infinite hypercycle-arc side
of $C$ lying on $\psi L_2$ is also mapped to itself 
by the congruence. Then
by induction
one sees that each hypercycle-arc side of $C$ is mapped to
itself. Thus all three vertices of some triangle
(one infinite and two finite vertices of the first two
arc-sides of $C$) are preserved, hence the congruence
is trivial.

If the pair of infinite points of $C$ is exchanged by the
congruence, then
similarly one sees that all the hypercycle-arc
sides of $C$ are mapped to
all the hypercycle-arc sides of $C$, but their order is
reversed.
The congruence preserves also the property
on which side of $C$ is
$(\varphi K) \cap (\psi L)$ contained.
Hence the congruence is orientation-reversing, and
then by Part I (i.e., \cite{5}), 
\thetag{5} 
it is an axial symmetry. Then it exchanges
the hypercycle-arc sides of $C$ lying on
$\psi L_1$ and $\psi L_2$, hence it exchanges
$\psi L_1$ and $\psi L_2$ as well.
Therefore its axis is $a(\psi L)$, which lies on the
$\xi $-axis,
cf.\ the second paragraph of this proof. Moreover, it exchanges
the first and last arcs of $\varphi K_1$ on $C$ (these
may coincide).
Hence this congruence maps $\varphi K_1$ onto itself.
Therefore its axis $a(\psi L)$,
which intersects $\varphi K_1$ transversally,
is orthogonal to $\varphi K_1$. 

However, a
small rotation of $\varphi K$, and thus also of $\varphi K_1$,
about the intersection point of the axis of symmetry
$a(\psi L)$
and $\varphi K_1$, preserves 
the combinatorial structure of $(\varphi K) \cap (\psi
L)$ (only possibly the set of indices 
$\{ i(1),...,i(k) \} $ will
change, but this does not invalidate the above
considerations). Simultaneously, this small rotation
destroys orthogonality of the intersection of $a(\psi L)$
and $\varphi K_1$,
hence destroys also this unique non-trivial congruence.
$\blacksquare $
\enddemo


{\it{Proof of Theorem}} 5, {\bf{continuation.}}
{\bf{6.}}
Now the earlier parts of the proof of Theorem 5, and Lemmas
5.1--5.10 end the proof of Theorem 5.
$\blacksquare $
\enddemo


\demo{Proof of Theorem {\rm{6}}}
{\bf{1.}}
We have $(4) \Longrightarrow (1)$ by \cite{4}, 
Theorem 1.
We have evidently $(1) \Longrightarrow (2)$. There remains
to prove $(4) \Longrightarrow (3)$ and 
$(2) \Longrightarrow (4)$ and $(3) \Longrightarrow (4)$.
These we do in the next three lemmas.


\proclaim{Lemma 6.1}
In Theorem {\rm{6}} we have $(4) \Longrightarrow (3)$.
\endproclaim


\newpage

\noindent
\demo{Proof}
We copy the proof of \cite{4}, 
Lemma 1.3, which says just
$(4) \Longrightarrow (1)$ in our Theorem 6.

{\bf{1.}}
Let $K$ and $L$ be congruent balls (in $S^d$, ${\Bbb{R}}^d$
or  $H^d$).

For $\varphi K = \psi L$ we do not have Theorem 6, (2), (C),
hence the implication ${\text{(A)}} \land {\text{(B)}}
\land {\text{(C)}}
\Longrightarrow (3)$ is satisfied. Further let
$\varphi K \ne \psi L$.

We denote by $o$ the midpoint of the
segment $S$ with endpoints the centres of $\varphi K$ and
$\psi L$, and by $H$ the orthogonal bisector hyperplane of
$S$. Then by ${\text{int}}\,[(\varphi K) \cap
(\psi L)] \ne \emptyset $
the length of $S$ is less than the diameter of $K$, hence
$o \in {\text{int}}\,[(\varphi K) \cap (\psi L)]$.
Then $\varphi K$ and $\psi L$ are symmetric images of each
other w.r.t.\ $H$. Hence both $(\varphi K) \cap (\psi L)$ and
$S$ are symmetric w.r.t.\ $H$.
We claim that for $z \in (\varphi K) \cap (\psi L)$, the
maximum distance from ${\text{bd}}\,[(\varphi K) \cap
(\psi L)] $ is attained only for $z = o$.

We suppose that the image of $H$ in the conformal model lies
on the $x_1\dots x_{d-1}$-coordinate hyperplane,
the image of $S$
lies on the $x_d$-axis, whence the image of $o$ is $0$.
We say that
$\varphi K$ {\it{lies lower than}} $\psi L$.
We may
suppose that  $z \in (\varphi K) \cap (\psi L)$
lies on or above $H$, and
$z \not\in {\text{bd}}\,(\varphi K)$. Then the minimum
distance of $z$ and the part of ${\text{bd}}\,(\varphi K)$
above $H$ is realized by a geodesic segment, from $z$ to some
point of ${\text{bd}}\,(\varphi K)$ above $H$,
which orthogonally intersects
${\text{bd}}\,(\varphi K)$. Then this geodesic
segment lies on a radius of $\varphi K$. Fixing this
radius and varying $z$ on it,
the maximum distance is attained only for
$z \in H$.
Consider a ball $\varphi K_0$, concentric with $\varphi K$,
such that $o \in {\text{bd}}\,(\varphi K_0)$.

Then on every radius of $\varphi K$,
the signed distances of
the intersection points of this radius
with ${\text{bd}}\,(\varphi K)$ and
${\text{bd}}\,(\varphi _1 K)$ are constant.
By the conformal model, 
$(\varphi K_0) \setminus \{ o \} $ lies strictly below $H$.
Therefore the maximum distance of $z \in (\varphi K) \cap
(\psi L)$, lying
on or above $H$, to the part of ${\text{bd}}\,
(\varphi K)$ above $H$, is attained only for $z = o$.

This implies our claim that
for $z \in (\varphi K) \cap (\psi L)$, the
maximum distance from ${\text{bd}}\,[(\varphi K) \cap
(\psi L)] $ is attained only for $z = o$. It also follows
that
the farthest 
points of ${\text{bd}}\,[(\varphi K) \cap (\psi L)]$ from
$o$ are the intersections of the line, spanned by $S$, with
${\text{bd}}\,[(\varphi K) \cap (\psi L)]$.

Therefore $(\varphi K) \cap (\psi L)$
has a unique inball, with centre $o$.
So Theorem 6, (2), (B) is satisfied. Moreover,
Theorem 6, (2), (C) is also satisfied, with $\{ \varphi x,
\psi y \} $ being the intersection of
${\text{bd}}\,[(\varphi K) \cap (\psi L)]$
with
the straight line, spanned by $S$.
Last, Theorem 6, (3) is also satisfied, with the
orthogonal bisector hyperplane
of $S = [\varphi x, \psi y]$ being $H$.

{\bf{2.}}
Now let us consider any
compact intersection of two paraballs $\varphi K$ and $\psi
L$, with nonempty
interior (in $H^d$).
Then, like in \cite{4}, 
Lemma 1.3, their
infinite points, $\varphi k$ and $\psi l$, say, are
different. 
Let the other points of
${\text{bd}}\,(\varphi K)$ and ${\text{bd}}\,(\psi L)$
on the line
$(\varphi k)(\psi l)$ be $\varphi k'$ and $\psi l'$.
Like in \cite{4}, 
Lemma 1.3,
the order of these points on this line is  
$\varphi k, \psi l', \varphi k', \psi l$, and these points
are different. 
We denote by $o$ the midpoint of $S:=[\varphi k', \psi l']$
(then $o \in {\text{int}}\,[(\varphi K) \cap (\psi L)]$),
and by $H$ the orthogonal bisector hyperplane of
$S$. 
Then $\varphi K$ and $\psi L$ are symmetric images of each
other w.r.t.\ $H$.
Hence both $(\varphi K) \cap (\psi L)$ and
$S$ are symmetric w.r.t.\ $H$.
We claim that for $z \in (\varphi K) \cap (\psi L)$, the
maximum distance from ${\text{bd}}\,[(\varphi K) \cap
(\psi L)] $ is attained only for $z = o$. 

We suppose that the image of $H$ in the conformal model lies
on
the $x_1 \dots x_{d-1}$-coordinate hyperplane,
and the image of 
$\varphi k$ or $\psi l$ is $(0, \dots, 0, -1)$ or
$(0, \dots, 0,1)$,

\newpage

\noindent
resp., whence the image of $o$ is $0$.
We say that $\varphi K$
{\it{lies lower than}} $\psi L$. Like in {\bf{1}}, we may
suppose that $z \in (\varphi K) \cap (\psi L)$
lies on or above $H$, and
$z \not\in {\text{bd}}\,(\varphi K)$.
Then the minimum distance of $z$ and the part of
${\text{bd}}\,(\varphi K)$ above $H$ is realized by
a geodesic segment, from
$z$ to some point of ${\text{bd}}\,(\varphi K)$
above $H$,
which orthogonally intersects ${\text{bd}}\,(\varphi K)$.
Then this geodesic segment lies on a straight
line with infinite point
$\varphi k$. Fixing this line and varying $z$ on it,
the maximum distance is attained only for $z \in H$. 
Let us consider another paracycle $\varphi _1 K$, having
$\varphi k$ as its infinite point, and with $o \in
{\text{bd}}\,(\varphi _1 K)$ (this is unique).

Then on every
straight line with $\varphi k$ as an infinite point, the signed
distances of
the intersection points of this line
with ${\text{bd}}\,(\varphi K)$ and
${\text{bd}}\,(\varphi _1 K)$ are constant. (The 
intersection point with ${\text{bd}}\,(\varphi _1 K)$ lies
between $\varphi k$ and the intersection point with 
${\text{bd}}\,(\varphi K)$.) However, by the
conformal model, 
$(\varphi _1 K) \setminus \{ o \} $ lies strictly below $H$.
Therefore the maximum distance of $z \in (\varphi K) \cap
(\psi L)$, lying
on or above $H$, to the part of ${\text{bd}}\,
(\varphi K)$ above $H$, is attained only for $z = o$.

This implies our claim that
for $z \in (\varphi K) \cap (\psi L)$, the
maximum distance from ${\text{bd}}\,[(\varphi K) \cap
(\psi L)] $ is attained only for $z = o$. 
It also follows that the farthest points of
${\text{bd}}\,[(\varphi K)\cap (\psi L)]$ from $o$ are
$\varphi k'$ and $\psi l'$.

Therefore $(\varphi K)\cap (\psi L)$ has a unique inball,
with centre $o$. So
Theorem 6, (2), (B) is satisfied. Moreover,
Theorem 6, (2), (C) is also
satisfied, with $\{ \varphi x, \psi y \} $ (where $x = k'$
and $y = l'$) being the
intersection of ${\text{bd}}\,[(\varphi K) \cap (\psi L)]$
with the straight
line $(\varphi k)(\psi l)$. Last, Theorem 6, (3) is
also satisfied, with the orthogonal bisector hyperplane of
$S = [\varphi k', \psi l'] = [\varphi x, \psi y]$ being $H$.

{\bf{3.}} 
There remains the case when the connected components of both
$K$ and $L$ are congruent hyperspheres, with common positive
distance, $\lambda $, say, from their respective
base hyperplanes (in $H^d$). In \cite{4}, 
proof of Lemma 1.3 it was
shown that different hypersphere boundary components of
$\varphi K$, or of $\psi L$, have a distance at least
$2 \lambda $, and supposing ${\text{diam}}[(\varphi K) \cap
(\psi L)] < 2 \lambda $, for some 
boundary components $\varphi K_i$ of $\varphi K$ and 
$\psi L_j$ of $\psi L$, we have $(\varphi K) \cap
(\psi L) = (\varphi K_i^*) \cap (\psi L_j^*)$. Here 
$\varphi K_i^*$, or $\psi L_j^*$, is the closed convex set
bounded by $\varphi K_i$, or by $\psi L_j$, resp.
(Later we will use only compactness of
$(\varphi K_i^*) \cap (\psi L_j^*)$.)
By \cite{4}, 
Lemma 1.1, the parallel domain of the base
hyperplane of $K_i$, or of $L_j$, with distance $\lambda $,
lies in $K$, or in $L$, resp.

Suppose that the base hyperplanes of
$\varphi K_i$ and $\psi L_j$ have a common finite point. Then
$\emptyset \ne 
{\text{int}}\,[(\varphi K_i^*) \cap (\psi L_j^*)]$. Moreover,
$(\varphi K_i^*) \cap (\psi L_j^*)$
has some infinite
point,
hence ${\text{diam}}\,[(\varphi K) \cap (\psi L)] = \infty $,
a contradiction.

Suppose that the 
base hyperplanes of
$\varphi K_i$ and $\psi L_j$ have a (unique)
common infinite point $p$.
Then in the conformal model, as a subspace of ${\Bbb{R}}^d$,
the stereoangles of 
$\varphi K_i^*$ and $\psi L_j^*$ at $p$
are greater than $1/4$ times
the total stereoangle, while the conformal model ball,
as a subspace of ${\Bbb{R}}^d$, has at
$p$ a stereoangle $1/2$ times the total stereoangle.
Hence $\varphi K_i^* \cap \psi L_j^*
\ne \emptyset $, and it has $p$ as an infinite point, again
implying ${\text{diam}}\,[(\varphi K) \cap (\psi L)] =
\infty $, a contradiction.

Hence the base hyperplanes
of $\varphi K_i$ and $\psi L_j$ have 
no common finite or infinite point. Let $S$ be the segment
realizing the distance of these two base hyperplanes,  
$o$ be its midpoint, and $H$ be its
orthogonal bisector hyperplane.
We say that the base hyperplane
of $\varphi K_i$ {\it{lies below
that of}} $\psi L_j$.
Then $\varphi K_i$ {\it{lies above its
base}}

\newpage

\noindent
{\it{hyperplane}},
and 
$\psi L_j$ {\it{lies below its base hyperplane}}.
Namely, in the remaining
three cases, 
$(\varphi K_i^*) \cap (\psi L_j^*)$ would have infinite points,
and then ${\text{diam}}\,[(\varphi K) \cap (\psi L)] =
\infty $, a contradiction.
Then $\varphi K$ and $\psi L$ are symmetric images of each
other w.r.t.\ $H$. Hence both $(\varphi K) \cap (\psi L)$ and
$S$ are symmetric w.r.t.\ $H$.
Also, using \thetag{*},
${\text{int}}\,[(\varphi K) \cap (\psi L)] \ne \emptyset
$ is symmetric w.r.t.\ $o$, hence
$o \in {\text{int}} [(\varphi K) \cap (\psi L)]$.
We claim that for $z \in (\varphi K) \cap
(\psi L)$, the maximum distance from ${\text{bd}}\,
[(\varphi K) \cap (\psi L)]$ is attained only for $z = o$.

We suppose that the image of $H$ in the conformal model
lies on the $x_1 \dots x_{d-1}$-hyperplane,
and the image of $S$ lies on the $x_d$-axis,
whence the image of $o$ is $0$.
Like in {\bf{1}} and {\bf{2}}, we may
suppose that $z \in (\varphi K) \cap
(\psi L)$ lies on or above $H$,
and $z \not\in {\text{bd}}\,(\varphi K_i)$.
Then the minimum distance of $z$ and the part of
$\varphi K_i$ above $H$ is realized by
a geodesic segment, from
$z$ to some point of $\varphi K_i$ above $H$,
which orthogonally intersect $\varphi K_i$.
Thus this geodesic segment lies on a straight
line orthogonal to the
base hyperplane of $\varphi K_i$. Fixing this line and
varying $z$ on it,
the maximum distance is attained only for $z \in H$. 
Let us consider another hypersphere ${\varphi {\tilde{K_i}}}$,
with the same base hyperplane as $\varphi K_i$, such that
$o \in \varphi {\tilde{K_i}}$ (this is unique).

Then on every straight line, orthogonal to our base hyperplane,
the signed distances of the intersection points of this line
with $\varphi K_i$ and $\varphi {\tilde{K_i}}$ are constant. 
However, by the
conformal model, $(\varphi {\tilde{K_i}})
\setminus \{ o \} $
lies strictly below $H$. Therefore the maximum distance of 
$z \in (\varphi K) \cap (\psi L)$, lying on or above $H$,
to the part of $\varphi K_i$ above $H$,
is attained only for $z = o$.

This implies our claim that for $z \in (\varphi K)
\cap (\psi L)$, the maximum distance from ${\text{bd}}\,[
(\varphi K) \cap (\psi L)]$ is attained only for $z = o$.  
It also follows that the farthest points of
${\text{bd}}\,[\varphi K) \cap (\psi L)]$ from $o$ are the 
two points of intersection of ${\text{bd}}\,[(\varphi K) \cap
(\psi L)]$ and the line spanned by $S$.

Therefore $(\varphi K) \cap (\psi L)$ has a unique inball,
with centre $o$. So Theorem
6, (2), (B) is satisfied. Moreover, Theorem
6, (2), (C) is also satisfied, with $\{ \varphi x, \psi y \} $
being the intersection of ${\text{bd}}\,[\varphi K) \cap
(\psi L)]$ with the straight line spanned by $S$. Last, 
Theorem 6, (3) is also satisfied, with the orthogonal bisector
hyperplane of $[\varphi x, \psi y]$ being $H$.
$\blacksquare $
\enddemo


We note that in \cite{4}, 
as a part of the proof of
Lemma 1.6,
uniqueness of the inball
of $(\varphi K) \cap (\psi L)$, for sufficiently small
diameter,
was proved for $K,L$ being
balls in $S^d$ and ${\Bbb{R}}^d$, and for $K,L$ having
congruent hypersphere boundary components (with
hyperplanes excluded) in $H^d$. In the special case of our
Lemma 6.1,
the proof in this paper
proves a bit more,
is shorter, and all three cases {\bf{1}}-{\bf{3}} are treated 
analogously.

\vskip.1cm


\proclaim{Lemma 6.2}
In Theorem {\rm{6}} we have $(2) \Longrightarrow (4)$.
\endproclaim


\demo{Proof}
Observe that \cite{4}, 
proof of Theorem 1, $(1)
\Longrightarrow (2)$ (i.e., our Theorem 6, $(1)
\Longrightarrow (4)$) did not use central symmetry of all
intersections $(\varphi K) \cap (\psi L)$ of
``small'' diameter, but only those
which satisfied (B) and (C) of Theorem 6, (2). (Moreover,
in \cite{4}, 
Theorem 1 there occured $\varepsilon =
\varepsilon (K,L)$, while in our Theorem 6, (2) we have
$\varepsilon = \varepsilon (K,L,x,y)$. This difference will be
discussed later.)

In fact, in \cite{4}, 
Lemma 1.6 (2), and in
its proof (just below
(1.6.1))
the centre of symmetry
$c$ of $(\varphi K) \cap (\psi L)$ was identified as the point
$O$ such that for some

\newpage

\noindent
$\varepsilon _1 > 0$ we had that
$B(O, \varepsilon _1)$ was the unique inball of 
$(\varphi K) \cap (\psi L)$. This is just our Theorem 6, (2)
(B).
(Observe that in case of \cite{4}, 
Lemma 1.6 (1) we
had that $X = S^d$ and $K,L$ were halfspheres, and
then Theorem 6 (4) held.)

To show Theorem 6, (2) (C), we have to recall the construction
of $B(O, \varepsilon _1)$ from \cite{4}. 
In \cite{4}, 
(A) (which is \thetag{**} in this paper), for $x \in
{\text{bd}}\,K$ and $y \in {\text{bd}}\,L$, there was
asserted the existence of two balls, with the properties
stated there. However, if the radii of these balls were chosen
less than $\varepsilon _1(x)$,
or $\varepsilon _1(y)$, resp.\ (cf.\ \cite{4}, 
(1.5.7)),
then these balls, except their points $x$ or $y$, lied
actually
in ${\text{int}}\,K$, or ${\text{int}}\,L$, resp. We supposed
that the radii of these new balls were equal to some 
$\varepsilon _1 \in (0, \min \{ \varepsilon _1 (x),
\varepsilon _1 (y) \} /2)$. 

In the first paragraph of the
proof of \cite{4} 
Lemma 1.5, it was proved that
either $X = S^d$ and $K,L$ were halfspheres, and then 
Theorem 6 (4) held, or, e.g., $K$ had an
exposed point. Suppose this second case. Then in
\cite{4} 
Lemma 1.7 (2) it was proved that both $K$ and
$L$ were strictly convex.
(\cite{4} 
Lemma 1.7 (1) asserted the first case from the
beginning of this paragraph.)
In \cite{4} 
(1.5.4), we chose a point
$O \in X$ (``origin''), and chose $\varphi _0$ and $\psi _0$
so that $\varphi _0 x = O$ and $\psi _0 y = O$, and 
$\varphi _0 K$ and $\psi _0 L$ touched each other from
outside, at $O$ (recall $(**) \Longrightarrow C^1$). 

Then in \cite{4} 
(1.5.6) and (1.5.7) we moved
$\varphi _0 K$ and $\psi _0 L$, and with them the new balls
as rigidly attached to them, toward each other. The
centres of the new balls moved along the common normal line
of $\varphi _0 K$ and $\psi _0 L$
at $O$, through the distance $\varepsilon _1$, while we
allowed any simultaneous, independent rotations of
$\varphi _0 K$ and $\psi _0 L$ about their original common
normal. The new positions of $\varphi _0 K$ and $\psi _0 L$
were denoted by $\varphi K$ and $\psi L$.
Then the
centres of the new balls moved 
to their common new position $O$.
From above, originally these new
balls, except for their points $\varphi _0x$ and $\psi _0y$,
resp., were contained in ${\text{int}}\,(\varphi _0 K)$ and
${\text{int}}\,(\psi _0 L)$, resp. Therefore in their common
new position these new
balls lied in 
$[{\text{int}}\,(\varphi K)] \cap [{\text{int}}\,
(\psi L)] =
{\text{int}}\,[(\varphi K) \cap (\psi L)]$, except for two
points of their common boundary, namely 
$\varphi x$ and $\psi y$. Moreover, these two
points
lied on the original common
normal (hence were antipodal), and were the
centres of the new balls before the translation and
rotations.

The later steps of the proof only used central symmetry in the
above described positions. Moreover, 
observe that in \cite{4}, 
proof of Theorem 1, $(1)
\Longrightarrow (2)$ in (1) there occured the hypothesis
$\varepsilon = \varepsilon (K,L)$, stronger than in our
Theorem 6, (2). However, actually in the
proof only one $x \in {\text{bd}}\,K$ and one
$y \in {\text{bd}}\,L$ were used at a time,
to prove {\it{equality of all
sectional curvatures of\, ${\text{\rm{bd}}}\,K$ at
$x$ with all
sectional curvatures of ${\text{\rm{bd}}}\,L$ at $y$ (and
their positivity for $X = {\Bbb{R}}^d$ and $X = H^d$)}}.
Hence we may admit that $\varepsilon $ depends also on $x$
and $y$, as in our Theorem 6, (2).

In fact, recall that $\varphi K$ and
$\psi L$ were not uniquely determined, but they could be
independently rotated about the original common normal line.
(This line contained $O$. Suppose that the image of $O$ in
one of the 
models is $0$. Then the images of these rotations in the model
are Euclidean rotations about the image of this line in the
model.)
{\it{Therefore any two-dimensional normal section of
${\text{\rm{bd}}}\,
(\varphi K)$ could become the centrally symmetric image of
any two-dimensional normal section of ${\text{\rm{bd}}}\,
(\psi L)$,
with the central symmetry interchanging $\varphi x$ and
$\psi y$.
Hence the
centre of symmetry was the midpoint of
$[\varphi x, \psi y]$.}} (For $S^d$ another centre of
symmetry was the point of $S^d$ antipodal
to this midpoint. However, the
central symmetries w.r.t.\ these two

\newpage

\noindent
centres as maps coincide,
so we may leave the other centre of
symmetry
out of consideration.)
From above the curvature of any two-dimensional normal section
of ${\text{bd}}\,(\varphi K)$ at $\varphi x$, and
the curvature of any two-dimensional normal section of
${\text{bd}}\,(\psi L)$ at $\psi y$ coincide, as claimed
above. What was written above, readily implies also the 
{\it{local rotational symmetry of both ${\text{\rm{bd}}}\,K$
and
${\text{\rm{bd}}}\,L$, at any of their points $x$ and $y$
(w.r.t.\ the normal of ${\text{\rm{bd}}}\,K$ at $x$, and
that of ${\text{\rm{bd}}}\,L$ at $y$ -- which exist by
\thetag{**})}}.

In \cite{4}, 
Lemma 1.5 (1)
constantness and equality (and for ${\Bbb{R}}^d$ and $H^d$
also positivity)
of all
sectional curvatures of ${\text{bd}}\,K$ and
${\text{bd}}\,L$ at all their
points $x$ and $y$ were proved.
In \cite{4}, 
Lemma 1.5 (2) the
local rotational symmetry of both ${\text{bd}}\,K$ and
${\text{bd}}\,L$, at any of their points $x$ and $y$
was proved.
From these the global statement of \cite{4}, 
Theorem 1,
(2), i.e., 
the global statement of (4) of our theorem, was
deduced in \cite{4}, 
Lemmas 1.8 and 1.9. This ends the
proof of $(2) \Longrightarrow (4)$ in Theorem 6.
$\blacksquare $
\enddemo


\proclaim{Lemma 6.3}  
In Theorem {\rm{6}} we have $(3) \Longrightarrow (4)$.
\endproclaim


\demo{Proof}
We repeat the proof of Lemma 6.2, except that we replace the
first sentence in italics, in its last but three
paragraph, by the following.
{\it{Therefore any two-dimensional normal section of
${\text{\rm{bd}}}\,
(\varphi K)$ could become the symmetric image of
any two-dimensional normal section of ${\text{\rm{bd}}}\,
(\psi L)$,
via the symmetry w.r.t.\ the orthogonal bisector hyperplane
of $[\varphi x, \psi y]$.}}
$\blacksquare $
\enddemo


{\it{Proof of Theorem}} {\rm{6}}, {\bf{continuation.}}
{\bf{2.}}
Now the proof of Theorem 6 follows from the earlier parts of
the proof of Theorem 6, and from Lemmas 6.1-6.3.
$\blacksquare $
\enddemo


The proof of Theorem 7 will use several ideas from the
proof of
Theorem 4 in \cite{4}. 


\demo{Proof of Theorem {\rm{7}}}
{\bf{1.}}
Suppose (7) of Theorem 7. By \thetag{****}, for $S^d$
the centres of 
$\varphi K$ and $\psi L$ are not antipodal. Hence 
the midpoint of the (shorter) segment joining these centres 
is a centre of symmetry of
$$
M :=
{\text{conv}}\, [ (\varphi K) \cup (\psi L) ] .
\tag 7.1
$$
Hence Theorem 7,
$(7) \Longrightarrow (1)$ holds. 
Moreover, if we suppose (7) of Theorem 7, then either
$\varphi K = \psi L$ and then (2) of Theorem 7 holds vacuously,
or else there is a unique diametral segment $D$ of
${\text{conv}}\, [ (\varphi K) \cup (\psi L) ]$, and
the orthogonal bisector hyperplane of $D$
is the orthogonal bisector hyperplane of the
segment with endpoints the centres of $K$ and $L$, and $M$
is symmetrical w.r.t.\ it. Hence Theorem 7,
$(7) \Longrightarrow (2)$ holds.

The implications Theorem 7, $(1) \Longrightarrow (3)$ and
$(2) \Longrightarrow (4)$ are evident. 

There remains to show Theorem 7, $(3) \Longrightarrow (5)$,
$(4) \Longrightarrow (6)$, $(5) \Longrightarrow (7)$ and
$(6) \Longrightarrow (7)$.
 

\proclaim{Lemma 7.1}
Assume \thetag{****}.
Assume the hypotheses of Theorem {\rm{7}}. 
Then in Theorem {\rm{7}} we have $(3) \Longrightarrow (5)$ and
$(4) \Longrightarrow (6)$.
\endproclaim


\demo{Proof}
Suppose (3) or (4) of Theorem 7, resp.
As in part {\bf{2}} of the proof of Theorem 4 in
\cite{4} 
(this part does not use smoothness of $K$ and $L$), we have
the following.
For any $x \in {\text{bd}}\,K$ and $y \in {\text{bd}}\,L$,
there exists a ball $B(o,R) \subset X$ (for
$S^d$ we may

\newpage

\noindent
suppose even that $B(o,R)$ lies in the open
southern
hemisphere, which implies that $B(o,R)$ does not contain any
antipodal pair of points of $S^d$), with centre
$o$ and radius $R$ (less than $\pi /2$ but
arbitrarily close to $\pi /2$ for $S^d$, and
arbitrarily large for ${\Bbb{R}}^d$ and $H^d$), and there exist
$\varphi $ and $\psi $ with the following properties. 
\newline
(1)
$\varphi x, \psi y$ are antipodal points of
${\text{bd}}\,B(o,R)$ (thus $\varphi x \ne \psi y$;
but, in the case of $S^d$, they
are not antipodal in $S^d$),
\newline
(2) 
$\varphi K, \psi L \subset B(o,R)$; even $\varphi K$ and
$\psi L$ are contained in balls of radius some $r < R$,
contained in $B(o,R)$, and tangent to $B(o,R)$ at $\varphi x$
and $\psi y$, resp.,
\newline
(3)
by (2), $\left( {\text{bd}}\,B(o,R) \right) \cap
[ {\text{bd}}\,(\varphi K) ]
= \{ \varphi x \} $, and 
$\left( {\text{bd}}\,B(o,R) \right)
\cap [ {\text{bd}}\,(\psi L) ]
= \{ \psi y \} $. Then, in the collinear model,
any point of ${\text{bd}}\,B(o,R)$, except $\varphi x$ and
$\psi y$,
can be separated by a hyperplane, intersecting
${\text{int}}\,B(o,R)$, from both of $\varphi K$ and $\psi L$,
hence also from $M$.
Therefore $\left( {\text{bd}}\,B(o,R) \right)
\cap ( {\text{bd}}\,M ) = \{
\varphi x , \psi y \} $.
\newline
(4) $B(o,R)$ is the unique ball of minimal radius containing $M$.
\newline
Then (4) implies that any congruence admitted by 
$M$ leaves $B(o,R)$ invariant. This holds in particular for
central symmetry, and for symmetry w.r.t.\ the orthogonal
bisector hyperplane of the unique diametral segment $D :=
[\varphi x, \psi y]$ of $M$.
Hence for central symmetry the centre of symmetry
is $o$, and both of these symmetries
interchange $\varphi x$ with $\psi y$. (For $S^d$ another
centre of symmetry is $-o$, but the central symmetries
w.r.t.\ $o$ and $-o$, as set maps, are identical, so we may
leave $-o$ out of consideration.)

By its definition, $D(\varepsilon )$ is also centrally
symmetric w.r.t.\ $o$, and also is symmetric w.r.t.\ the
orthogonal bisector hyperplane of $D$.
Therefore in case (3) of Theorem 7, we have that
$M \cap D(\varepsilon )$ is
also centrally
symmetric (w.r.t.\ $o$). Similarly, in case (4) of
Theorem 7,
we have that $M \cap D(\varepsilon )$ is
also symmetric w.r.t.\ the
orthogonal bisector hyperplane of $D$.
Hence (5) or (6)
of Theorem 7 holds, resp., for any $\varepsilon (K,L,x,y,
\varphi , \psi , S)
\in (0, \pi /2)$ for $S^d$, and for any $\varepsilon (K,L,x,y,
\varphi , \psi )
\in (0, \infty )$ for ${\Bbb{R}}^d$ and $H^d$.
$\blacksquare $
\enddemo


So there remains to prove $(5) \Longrightarrow (7)$
and $(6) \Longrightarrow (7)$ of Theorem 7.


\proclaim{Lemma 7.2}
Assume \thetag{****}. 
Assume the hypotheses of Theorem {\rm{7}}.
Then both of
{\rm{(5)}} and {\rm{(6)}}
of Theorem {\rm{7}} imply that both of $K$ and $L$
are $C^1$.
\endproclaim


\demo{Proof}
It suffices to give the proof for $K$, the other
case is analogous.

Suppose that $x \in {\text{bd}}\,K$ is not a
smooth point of ${\text{bd}}\,K$. Recall that convex bodies are
almost everywhere smooth (cf.\ \cite{6}). 
(This takes over from ${\Bbb{R}}^d$
to $S^d$ and $H^d$ via the collinear models.) Then let
$y \in {\text{bd}}\,L$ be a smooth point of ${\text{bd}}\,L$.
Further, let $o$ be the midpoint of $D$ (from Theorem 7,
(3), (B)).

Now we consider the collinear model, so that $o$ becomes the
centre $0$ of the model. Then the image of $B(o,R) \subset X$
from the proof of Lemma 7.1 
is a
Euclidean ball in the collinear model, with centre $0$.
The images of points, sets in this collinear model, as a subset of ${\Bbb{R}}^d$, will be
denoted by $'$, like $(\varphi x)'$, $\left( B(o,R) \right)'$,
etc. The radius
of $\left( B(o,R) \right)'$ is denoted by $R'$, and $D(\varepsilon )'$ is
in this model a right cylinder, of height $2R'$, and has bases
$(d - 1)$-balls, of radius $\varepsilon '$. Further,
diam$'\,( \cdot )$ denotes
the diameter in the Euclidean sense, in the collinear model, as a subset of ${\Bbb{R}}^d$. Moreover, we suppose that 
(1)-(4) from the proof of Lemma 7.1 hold.

Let $(\varphi x)'$ be the north
pole and $(\psi y)'$ the south pole of
$\left( B(o,R) \right)'$. Then the
horizontal supporting planes of $\left( B(o,R) \right)'$
at its poles are
supporting hyperplanes

\newpage

\noindent
of 
$M'$ (cf.\ \thetag{7.1}) and $M' \cap D(\varepsilon )'$
as well. At $(\psi y)'$ this is the unique supporting hyperplane
of $(\psi L)'$,
but at $(\varphi x)'$ there are also other supporting planes of
$(\varphi K)'$, even there are ones which are not horizontal,
but are arbitrarily close to horizontal. 
Let $H'$ be one such supporting plane. 

We have diam$'\, [ (\varphi K)' ] ,$
diam$'\, [ (\psi L)' ]
<$ diam$'\, [ \left( B(o,R) \right)' ] $. Therefore 
the highest point of $(\psi L)'$ lies strictly below the
height of the north pole. By a suitable, almost horizontal
choice of the
supporting plane $H'$ we can attain that the ball segment cut
off by $H'$ from $\left( B(o,R) \right)'$,
lying ``above'' $H'$,
lies entirely above the highest point of
$(\psi L)'$. Then both $(\varphi K)'$ and $(\psi L)'$ lie
``below'' $H'$, hence also 
$M'$ and $M' \cap D(\varepsilon )'$ lie ``below'' $H'$,
and then $H'$ is a supporting hyperplane of
$M' \cap D(\varepsilon )'$ at $(\varphi x)'$.

In case of
central symmetry of $M' \cap D(\varepsilon )'$ w.r.t.\
$o' = 0$ (for Theorem 7, (5))
then $-H'$ is a supporting hyperplane of
$M'\cap D(\varepsilon )'$ at $(\psi y)' = -
(\varphi x)'$. In case of symmetry of
$M' \cap D(\varepsilon )'$ w.r.t.\ the orthogonal bisector
hyperplane of $D$, we have that the symmetrical image
${\tilde{H}}'$ of $H'$ w.r.t.\ 
this orthogonally bisecting hyperplane (for Theorem 7, (6))
is a
supporting hyperplane of
$M'\cap D(\varepsilon )'$ at $(\psi y)' = -
(\varphi x)'$.
Since being
a supporting plane is a local property, therefore
then $-H'$, or ${\tilde{H}}'$
is a supporting hyperplane of $M'$ at $(\psi y)'$
as well, resp.
Then $-H'$, or ${\tilde{H}}'$ is
also a supporting hyperplane of $(\psi L)' \subset M'$
at $(\psi y)'$, resp.,
which is not horizontal, a contradiction to smoothness of
${\text{bd}}\,L$ at $y$.
$\blacksquare $
\enddemo

In Lemma 7.3 smoothness of $K$ and $L$ follows from Lemma 7.2.
However, also for the proof of Theorem 8 we will use
Lemma 7.3, where smoothness of $K$ and $L$ will follow from the
hypotheses of Theorem 8, for $S^2$, and from Lemma 8.2 for
${\Bbb{R}}^2$ and $H^2$. On the other hand,
strict convexity of $K$ and $L$ follows
both from the hypotheses of Theorem 7 and Theorem 8.


\proclaim{Lemma 7.3}
Assume \thetag{****}.
Assume the hypotheses of Theorem {\rm{7}}.
Let $K$ and $L$ be $C^1$, and let
$x \in {\text{\rm{bd}}}\,K$ and $y \in {\text{\rm{bd}}}\,L$.
Then, in the situation
described in
{\rm{(1)}} to {\rm{(4)}} of the proof of Lemma {\rm{7.1}},
and with $M$ from \thetag{7.1}, we have the following.
In some open
neighbourhood of
$\varphi x$, we have that ${\text{\rm{bd}}}\, (\varphi K)$
and  
${\text{\rm{bd}}}\,M$ coincide. Also, in some open
neighbourhood of
$\psi y$, we have that ${\text{\rm{bd}}}\, (\psi L)$
and  
${\text{\rm{bd}}} \,M$ coincide.
\endproclaim


\demo{Proof}
It suffices to give the proof for
${\text{\rm{bd}}}\,(\varphi K)$, the 
case of ${\text{\rm{bd}}}\,(\psi L)$ is analogous.

It suffices to show that in some neighbourhood of
$(\varphi x)'$ we have that ${\text{bd}}\,[(\varphi K)']$ and  
${\text{bd}}\, M'$ coincide.

We use the collinear model, and the 
coordinate system from the proof of Lemma 7.2 on it. We use
$'$ as in the proof of Lemma 7.2 (e.g., $o' = 0$).

Let $(\varphi x)' = R' e_d$; then $(\psi y)' = - R' e_d$.
Denoting by $h(M',u)$ the support function of $M'$, etc.,
for $u \in S^{d - 1}$,
we have
$$
\cases
h(M',e_d) = h \left( (\varphi K)', e_d \right) = R'
{\text{ and }}
\\
h((\psi L)', e_d) \le
{\text{\rm{diam}}}' \,[ (\psi L)' ] - R' =:
R' - c' < R' .
\endcases
\tag 7.3.1
$$
Since $M', (\varphi K)', (\psi L)'
\subset \left( B(o,R) \right)' = B(0,R')$,
therefore for their support functions we have
for $u_1,u_2 \in S^{d - 1}$ that
$|h(M',u_1) - h(M',u_2)| \le R' \cdot \| u_1 - u_2 \| $,

\newpage

\noindent
etc.

Let $u \in S^d$, with $R' \cdot \| u - e_d \| < c'/2$.
That is, $u$ can vary in the {\it{open spherical cap
$C' \subset S^{d - 1}$,
of centre $e_d$ and spherical radius
$2 \arcsin \,[c'/(4R')]$}}.
Then
$h \left( (\varphi K)', u \right) >
h(\left (\varphi K)', e_d \right) - c'/2 = R' - c'/2$, and
$h \left( (\psi L)', u \right) <
h \left( (\psi L)', e_d \right) + c'/2 < R' - c'/2$.
This implies that 
$$
\cases
u \in C' \Longrightarrow
h(M',u) = h \big( {\text{conv}} \,[ (\varphi K)' \cup
(\psi L)' ] , u \big) 
\\
 = \max \{ h \left( (\varphi K)', u
\right) , h \left( (\psi L)', u \right) \} =
h \left( (\varphi K)', u \right) .
\endcases
\tag 7.3.2
$$

By the supporting sphere
hypothesis, $K$ and $L$, hence also $(\varphi K)'$ and
$(\psi L)'$
are strictly convex, and by the hypothesis of this lemma they
are $C^1$.
Then their spherical maps (i.e., $(\varphi
x^* )' \in
{\text{bd}} \,[ (\varphi K)' ] $ is mapped to the
outer unit normal of $(\varphi K)'$
at $(\varphi x^* )'$, and similarly for $(\psi L)'$), are
homeomorphisms
${\text{bd}} \,[ (\varphi K)' ] \to S^{d - 1}$,
and ${\text{bd}} \,[ (\psi L)' ] \to S^{d - 1}$,
resp.

It will be convenient to identify a singleton set $\{ a \} $
with $a$.
The support function of a convex body $A \subset {\Bbb{R}}^d$
can be extended from $S^{d - 1}$ to ${\Bbb{R}}^d$,
in a positively homogeneous way (at $0$ it is defined as $0$).
Now recall the following (cf.\ \cite{1}, 
\cite{9}).  
\newline
(1) Suppose that for the supporting
hyperplane $H(A,u)$ of $A$,
with outer unit normal $u \in S^{d - 1}$, the
intersection $A \cap H(A,u)$ is a singleton. Then the gradient
of the extended support function of $A$ at $u$
exists, and equals this singleton point.
\newline
(2)
Conversely, if this gradient exists at $u \in S^{d - 1}$, then
$A \cap H(A,u)$ is a singleton.

Applying (1) to the strictly convex body
$(\varphi K)'$, we get that the gradient of
the extended support function of $(\varphi K)'$ exists
at any $u \in S^{d - 1}$, and its value is the singleton
$(\varphi K)' \cap H \left( (\varphi K)', u \right) $.
However, by \thetag{7.3.2}, on the open
subset $C' \subset S^{d - 1}$ the 
extended support function of $(\varphi K)'$ equals the 
extended support function of $M'$. Hence the gradient of
this second function also
exists on $C'$, and its value is the same
singleton $(\varphi K)'
\cap H \left( (\varphi K)', u \right)$.
Then applying (2) to the convex body $M'$, we get that
for $u \in C'$ the set $M' \cap H(M',u)$ is a singleton.
Once more applying (1), but now for the convex body
$M'$, we get that
this singleton equals the gradient of the extended support
function of $M'$ at $u \in C'$.
Thus $(\varphi K)' \cap H \left( (\varphi K)', u \right) $ and
$M' \cap H(M',u)$ are the (existing)
gradients of the same function, at
the same point $u \in C'$, hence they are equal.

Summing up: for $u \in C'$ the singleton support sets of 
$(\varphi K)'$ and $M'$ at $u$ coincide. The union of all 
these singleton sets, for all $u \in C'$, are the
inverse spherical images of the open subset
$C' \subset S^{d - 1}$, for $(\varphi K)'$ and $M'$. But the
spherical map and the inverse spherical map for
$(\varphi K)'$ are homeomorphisms, hence the inverse
spherical image $I'$ of $C' \subset S^{d - 1}$ is open in
${\text{bd}} \,[ (\varphi K)' ] $. By the invariance
of domain theorem for ${\text{bd}}\,(M')$ we have that $I'$
is open in ${\text{bd}}\,(M')$ as well.

Now we prove local coincidence of 
${\text{bd}} \,[ (\varphi K)' ] $ and 
${\text{bd}}\,(M')$ at $(\varphi x)'$.
In fact, let $B' \left( (\varphi x)', \delta ' \right) $
be a ball of centre $(\varphi x)'$ and radius $\delta '$
(meant in ${\Bbb{R}}^2$),
contained in ${\text{int}}\,B(0,1)$. Let
$\delta ' > 0$ be so small, that $I' \supset
[{\text{int}}\,B' \left( (\varphi x)', \delta ' \right) ]
\cap {\text{bd}}\, [(\varphi K)']$ and $I' \supset
[{\text{int}}\,B' \left( (\varphi x)', \delta ' \right) ]
\cap {\text{bd}}\,(M')$. Then
$$
\cases
z' \in [{\text{\rm{int}}}\,B' \left( (\varphi x)', \delta '
\right) ] \cap {\text{bd}}\, [(\varphi K)'] \Longrightarrow z'
\in I' \subset {\text{bd}}\, (M').
\\
{\text{Similarly, }} z' \in [{\text{int}}\,B' \left(
(\varphi x)', \delta ' \right) ] \cap {\text{bd}}\,(M')
\Longrightarrow z' \in {\text{bd}}\, [(\varphi K)']
\endcases
\tag 7.3.3
$$

\newpage

\noindent
$\blacksquare $
\enddemo


The next lemma will be used also in the proof of Theorem 8.
Observe that, under the hypotheses of Lemma 7.2, both $K$ and
$L$ are $C^1$. In the proof of Theorem 8, the $C^1$ property
of $K$ and $L$ is a hypothesis of Theorem 8 for $S^2$,
while for ${\Bbb{R}}^2$ and $H^2$, under the hypotheses of
Lemma 8.2, it will follow from Lemma 8.2.


\proclaim{Lemma 7.4}
Assume \thetag{****}.
Assume the hypotheses of Theorem {\rm{7}}.
Let $K$ and $L$ be $C^1$, and let
$x \in {\text{\rm{bd}}}\,K$ and $y \in {\text{\rm{bd}}}\,L$.
Then, in the situation
described in
{\rm{(1)}} to {\rm{(4)}} of the proof of Lemma {\rm{7.1}},
in some neighbourhood of $\varphi x$, and in some
neighbourhood of $\psi y$,
we have that ${\text{\rm{bd}}}\, M$ and  
${\text{\rm{bd}}}\,[M \cap D(\varepsilon )]$ coincide. 
\endproclaim


\demo{Proof}
We use the collinear model, and the 
coordinate system from the proof of Lemma 7.2 on it. We use
$'$ as in the proof of Lemma 7.2.

By the supporting sphere hypothesis, both $K$ and $L$ are
strictly convex. Moreover, by the hypothesis of this lemma,
they are also $C^1$.
Then
$(\varphi K)'$ and $(\psi L)'$ are also strictly convex and
$C^1$.

For our lemma it suffices to show that
$$
\cases
{\text{in some neighbourhood of }} (\varphi x)',
{\text{ and of }} (\psi y)', {\text{ we}}
\\
{\text{have that }}
{\text{bd}}\,(M') {\text{ and }} {\text{bd}}
\big[ [ M \cap D(\varepsilon ) ] ' \big] {\text{ coincide.}}
\endcases
\tag 7.4.1
$$

Let $[D_h(\varepsilon )]'$
be the closed parallel strip,
whose boundary hyperplanes are spanned by the two (horizontal)
bases of the cylinder $[D(\varepsilon )]'$. Let
$[D_v(\varepsilon )]'$ be the closed
both-way infinite cylinder,
with directrices spanned by the (vertical) directrices of
the cylinder $[D(\varepsilon )]'$.

We assert that
$$
\cases
{\text{the neighbourhood of }} (\varphi x)', {\text{ or of }}
(\psi y)', {\text{ whose existence}}
\\
{\text{is claimed in
\thetag{7.4.1}, can be chosen as }} 
[{\text{int}}\,D_v(\varepsilon )]'.
\endcases
\tag 7.4.2
$$

First observe that $\{ \varphi x, \psi y \} ' \subset
{\text{bd}}\,M'$ and $\{ \varphi x, \psi y \} ' \subset
{\text{bd}}\,[M \cap
D(\varepsilon )]'$. Therefore it remains to show that
$[({\text{bd}}\,M) \setminus \{ \varphi x, \psi y \} ]'$
and
$\big[ \big[ {\text{bd}}\,[M \cap D(\varepsilon )] \big]
\setminus \{ \varphi x,
\psi y \} \big]'$ coincide in some neighbourhood of 
$(\varphi x)'$, or of $(\psi y)'$, resp., where this 
neighbourhood is 
$[{\text{int}}\,D_v(\varepsilon )]'$.
That is, we have to prove
$$
\cases
[({\text{bd}}\,M) \setminus \{ \varphi x, \psi y \} ]' \cap
[{\text{int}}\,D_v(\varepsilon )]' =
\\
\big[ \big[ {\text{bd}}\,[M \cap D(\varepsilon )] \big]
\setminus \{ \varphi x,
\psi y \} \big]' \cap [{\text{int}}\,D_v(\varepsilon )]'.
\endcases
\tag 7.4.3
$$
We are going to prove that both sides of \thetag{7.4.3} contain
the other one. 

Recall the formulas ${\text{bd}}\,(A \cup B) \subset ({\text{bd}}\,A) \cup
({\text{bd}}\,B)$ and ${\text{bd}}\,(A \cap B) \subset ({\text{bd}}\,A) \cup
({\text{bd}}\,B)$. Then by the proof of Lemma 7.1, (2) and
(3) we have 
$$
\cases
[ ({\text{bd}}\,M) \setminus \{ (\varphi x), (\psi y) \} ]',
\big[
\big[ {\text{bd}}\,[M \cap D(\varepsilon )] \big] \setminus
\{ \varphi x, \psi y \} \big]'
\\
\subset [{\text{int}}\,B(o,R)]' \subset
[{\text{int}}\, D_h(\varepsilon )]',
\endcases
\tag 7.4.4
$$

\newpage

\noindent
which implies
$$
\cases
[({\text{bd}}\,M) \setminus \{ \varphi x, \psi y \} ]' =
[({\text{bd}}\,M) \setminus \{ \varphi x, \psi y \} ]' \cap
[{\text{int}}\,D_h(\varepsilon )]' {\text{ and}}
\\
\big[ \big[ {\text{bd}}\,[M \cap D(\varepsilon )] \big]
\setminus
\{ \varphi x, \psi y \} \big]' = 
\big[ \big[
{\text{bd}}\,[M \cap D(\varepsilon )] \big] \setminus
\{ \varphi x, \psi y \} \big]' \cap [ {\text{int}}\,
D_h(\varepsilon ) ]'.
\endcases
\tag 7.4.5
$$
Intersecting the right hand side
sets in \thetag{7.4.5} with $[{\text{int}}\,
D_v(\varepsilon )]'$, we get
$$
\cases
[({\text{bd}}\,M) \setminus \{ \varphi x, \psi y \} ]'
\cap
[ {\text{int}}\,D(\varepsilon ) ]' {\text{ and }}
\\
\big[ \big[ {\text{bd}}\,[M \cap D(\varepsilon )] \big]
\setminus \{ \varphi x, \psi y \} \big]' \cap
[ {\text{int}}\,D(\varepsilon ) ]', {\text{ resp.,}}
\endcases
\tag 7.4.6
$$
and we have to show that both of these sets contains the
other one.

We have
$$
\cases
[({\text{bd}}\,M) \setminus \{ \varphi x, \psi y \} ]'
\cap [ {\text{int}}\,D(\varepsilon ) ]'
= \big[ \big[ {\text{bd}} \big[ [M \cap D(\varepsilon )] \cup
[M \setminus D(\varepsilon )] \big] \big] \setminus
\\
\{ \varphi x, \psi y \} \big]'
\cap [ {\text{int}}\,D(\varepsilon ) ]'
\subset \big\{
\big[ \big[ {\text{bd}}\, [M \cap D(\varepsilon )] \big]
\setminus
\{ \varphi x, \psi y \} \big]'
\cap [ {\text{int}}\,D(\varepsilon ) ]' \big\} \cup
\\
\big\{ \big[
{\text{cl}}\,[X \setminus D(\varepsilon )] \big]'
\cap [ {\text{int}}\,D(\varepsilon ) ]' \big\} = 
\big[ \big[ {\text{bd}}\,[M \cap D(\varepsilon )] \big]
\setminus \{ \varphi x, \psi y \} \big]' \cap
[ {\text{int}}\,D(\varepsilon ) ]',
\endcases
\tag 7.4.7
$$
i.e., the first set in \thetag{7.4.6} is contained in the second
set in \thetag{7.4.6}.

On the other hand,
$$
\cases
\big[ \big[ {\text{bd}}\,[M \cap D(\varepsilon )] \big]
\setminus \{ \varphi x, \psi y \} \big]' \cap
[ {\text{int}}\,D(\varepsilon ) ]' \subset
\big[ \big[ ({\text{bd}}\,M) \cup
[{\text{bd}}\,D(\varepsilon )] \big]
\\
\setminus
\{ \varphi x, \psi y \} \big]' \cap
[ {\text{int}}\,D(\varepsilon ) ]' \subset
\big\{ \big[ ({\text{bd}}\,M) \setminus
\{ \varphi x, \psi y \} ]' \cap
[ {\text{int}}\,D(\varepsilon ) ]' \big\} \cup
\\
\big\{ [{\text{bd}}\,D(\varepsilon )]' \cap
[ {\text{int}}\,D(\varepsilon ) ]' \big\} =
[ ({\text{bd}}\,M) \setminus
\{ \varphi x, \psi y \} ]' \cap
[ {\text{int}}\,D(\varepsilon ) ]',
\endcases
\tag 7.4.8
$$
i.e., the second set in \thetag{7.4.6} is contained in the first
set in \thetag{7.4.6}.
This ends the proof that the two sets in \thetag{7.4.6}
coincide. Then the two sets in \thetag{7.4.3} also
coincide, and \thetag{7.4.1} is proved, and hence the
statement of the lemma is proved.
$\blacksquare $
\enddemo


\proclaim{Lemma 7.5}
Assume \thetag{****}.
Assume the hypotheses of Theorem {\rm{7}}. Let $K$ and $L$ be
$C^1$.
Then both of {\rm{(5)}} and {\rm{(6)}}
of Theorem {\rm{(6)}} imply that
both of ${\text{\rm{bd}}}\,K$ and ${\text{\rm{bd}}}\,L$ are at
each of their
boundary points twice differentiable {\rm{(cf.\ 6.
Preliminaries)}}.
Moreover, in the
situation described in {\rm{(1)}} to {\rm{(4)}}
in the proof of Lemma {\rm{7.1}},
for any $2$-plane $P$ containing $D = [\varphi x, \psi y]$,
the curvatures of $(\varphi K) \cap P$ at $\varphi x$ and
that of $(\psi L) \cap P$ at $\psi y$ coincide.
\endproclaim


\demo{Proof}
As stated in 6.
Preliminaries, convex surfaces are
almost everywhere twice differentiable, in the sense as given
there.
If, say,
${\text{bd}}\,K$ is not twice differentiable at its point
$x$, then we will
choose $y \in {\text{bd}}\,L$ so that ${\text{bd}}\,L$
is twice differentiable at $y$.

{\it{Now let us consider the situation described in {\rm{(1)}}
to {\rm{(4)}} in the proof of Lemma}} 7.1.
Then, by Lemma 7.3, in some neighbourhood of $\varphi x$ we
have that ${\text{bd}}\,(\varphi K)$ and ${\text{bd}}\,M$
coincide. Then 
${\text{bd}}\,M$ is not twice differentiable at
its point $\varphi x$.
Similarly, in some neighbourhood of $\psi y$ we
have that ${\text{bd}}\,(\psi L)$ and ${\text{bd}}\,M$
coincide. Then 
${\text{bd}}\,M$ is twice differentiable at
its point $\psi y$.

However, the central symmetry of $M \cap D(\varepsilon )$
about the centre $o$ (in case of (5) of Theorem 7),
or its symmetry w.r.t.\ the orthogonal
bisector hyperplane of $D$ (in

\newpage

\newpage

\noindent
case of (6) of
Theorem 7), interchanges $\varphi x$ and $\psi y$, as well as
some of their open neighbourhoods in ${\text{bd}}\,M$.
Then ${\text{bd}}\,M$ cannot
be not twice differentiable at $\varphi x$, but 
twice differentiable at $\psi y$. This contradiction proves
the first statement of the lemma.

The second statement of the lemma follows analogously: two
congruent,
twice differentiable curves
have at their points $\varphi x$
and $\psi y$,
corresponding to each other by the congruence, equal
curvatures.
$\blacksquare $
\enddemo


The following lemma is an analogue of \cite{4}, 
Lemma 1.8.


\proclaim{Lemma 7.6}
Assume \thetag{****}.
Assume the hypotheses of Theorem {\rm{7}}.
Then both {\rm{(5)}} and {\rm{(6)}}
of Theorem {\rm{7}} imply that
both ${\text{\rm{bd}}}\,K$ and ${\text{\rm{bd}}}\,L$
are, at each of their boundary points, locally,
relatively open subsets of spherical surfaces,
of some fixed sectional curvature $\kappa > 0$ for
$S^d$ and ${\Bbb{R}}^d$, and $\kappa > 1$ for $H^d$.
For $X = S^d$ the radius of these spherical surfaces is
less than $\pi /2$.
\endproclaim


\demo{Proof}
Again consider the situation described in (1) to (4) of the
proof of Lemma 7.1. Like in the proof of Lemma 6.2,
observe that $\varphi K$ and $\psi L$ are
not uniquely determined, but we may apply any rotations about
the line spanned by $D$ to any of them, independently of each
other.
Thus varying $\varphi $
and $\psi $ in this way, by Lemma 7.3, the central symmetry 
of $M \cap D(\varepsilon )$ about the centre $o$
(in case of Theorem 7, (5)), or its
symmetry w.r.t.\ the orthogonal bisector hyperplane of $D$
(in case of Theorem 7, (6)), may interchange locally a 
a $2$-dimensional normal
section of $M$, hence of $\varphi K$ at $\varphi x$, with any 
$2$-dimensional normal
section of $M$, hence of $\psi L$ at $\psi y$.
This immediately implies that both
${\text{bd}}\,(\varphi K)$ and ${\text{bd}}\,(\psi L)$
are locally rotationally symmetric about the line spanned by
$D$. By Lemma 7.5, the curvatures of the $2$-dimensional
normal sections of $\varphi K$ at $\varphi x$, and of
$\psi L$ at $\psi y$, exist and coincide.  

Therefore, all these (congruent)
$2$-dimensional normal sections have curvatures
at $\varphi x$ and $\psi y$, resp., and all these
curvatures are equal to some number $\kappa \ge 0$. Then the
supporting sphere hypothesis implies that $\kappa > 0$ for
$S^d$ and ${\Bbb{R}}^d$, and $\kappa > 1$ for $H^d$.

We still have to show that these $2$-dimensional normal 
sections are circular arcs. We claim that the 
$2$-dimensional normal sections $N$ of
$[ {\text{bd}}\,(\varphi K)] \cap
D(\varepsilon )$ containing
$\varphi x$ (or
of $[ {\text{bd}}\,(\psi L) ]
\cap D(\varepsilon )$ containing
$\psi y$) are $2$-dimensional
normal sections at any of their other points $\varphi x^*$
(or $\psi y^*)$, lying in the same open halfspace bounded by
the orthogonal bisector hyperplane of $D$, as $x$ (or $y$).
This means just that the normal of ${\text{bd}}\,(\varphi K)$
at $\varphi x^*$ lies in the $2$-plane spanned by $N$, which
is evident for $d = 2$. 
However, for $d \ge 3$, local rotational
symmetry of ${\text{bd}}\,(\varphi K)$ at $\varphi x$,
which necessarily has
as rotation axis the normal
of ${\text{bd}}\,(\varphi K)$ at $\varphi x$,
implies local symmetry of ${\text{bd}}\,(\varphi K)$
w.r.t.\ any hyperplane containing this axis
and $\varphi x^*$.
Thus the uniquely determined normal of
${\text{bd}}\,(\varphi K)$ at $\varphi x^*$ lies in any of
these hyperplanes, so it lies in their intersection, which is
just the $2$-plane spanned by this rotation axis and
$\varphi x^*$. Thus this intersection is 
spanned by $N$,
proving our claim.

Therefore $N$ has at each of its points $\varphi x^*$,
lying in the same open halfspace bounded by
the orthogonal bisector hyperplane of $D$, as $x$,
the
constant curvature $\kappa $. Thus $N$ locally is a cycle.
The already proved inequalities for $\kappa $ in this Lemma 
imply that $N$ locally is a circular arc. Therefore 
${\text{bd}}\,(\varphi K)$ is locally a spherical

\newpage

\noindent
surface, of
sectional curvature $\kappa $. For $X = S^d$, by the
supporting sphere hypothesis, this spherical surface has a
radius less $\pi /2$.
$\blacksquare $
\enddemo


{\it{Proof of Theorem}} 7, {\bf{continuation.}}
{\bf{2.}}
We apply Lemma 1.9 of \cite{4}. 
Its hypothesis, for the case
of the existence of the
supporting spheres at any boundary point, is just the
statement of our
Lemma 7.6. (Additionally there is as hypothesis (1) of
Theorem 1
or (1) of Theorem 4 of \cite{4}. 
However, these are actually
not
used in the proof of Lemma 1.9 of \cite{4}. 
Observe still that
the $C_{\alpha }$'s in the proof of Lemma 1.9 of
\cite{4} 
are correctly {\it{unions}} of equivalence classes there.)
Its conclusion implies that the connected components
both of $K$ and $L$ are congruent spheres (for $X = S^d$ of
radius at most $\pi /2$), paraspheres or
hyperspheres, but in case of congruent
spheres $K$ and $L$ are congruent
balls. However, the inequalities for $\kappa $ in our
Lemma 7.6 exclude paraspheres and hyperspheres. So
$K$ and $L$ are congruent balls, for $X = S^d$ of radius
less than $\pi /2$, by Lemma 7.6. This ends the
proof of Theorem 7.
$\blacksquare $
\enddemo


\demo{Proof of Theorem {\rm{8}}}
{\bf{1.}}
We have the evident implications $(7) \Longrightarrow (1)$,
and $(1) \Longrightarrow (2) \Longrightarrow (4)$, and
$(1) \Longrightarrow (3) \Longrightarrow (4)$, and 
$(5) \Longrightarrow (6)$.
Therefore we have to prove only the
implications $(3) \Longrightarrow (5)$, and
$(4) \Longrightarrow (6)$, and $(6) \Longrightarrow (7)$.

We will use the notation $M$ from \thetag{7.1}.


\proclaim{Lemma 8.1}
Assume \thetag{****} with $d = 2$.
Assume the hypotheses of Theorem {\rm{8}}.
Then in Theorem {\rm{8}} we have $(3) \Longrightarrow (5)$, and
$(4) \Longrightarrow (6)$.
\endproclaim


\demo{Proof}
We begin with copying
the proof of Lemma 7.1. We may suppose that there
is a circle $B(o,R)$, satisfying (1)-(4) from the proof of
Lemma
7.1, hence also (B) from Theorem 7 (3), and also satisfying
(A) from Theorem 7 (3). Hence the
hypotheses of Theorem 8, (3) and (4) are satisfied.

Observe that any congruence admitted by $M$
is a congruence admitted by its unique diametral segment
$D$ as well.
Therefore the 
non-trivial
congruences admitted by $M$
form a subset of $\{ $central symmetry w.r.t.\ $o$, 
axial symmetry w.r.t.\ the straight line spanned by $D$,
axial symmetry w.r.t.\ the orthogonal bisector line of $D \}$.
Observe that
$D(\varepsilon )$ admits all these three
non-trivial congruences. 

Hence, in case of Theorem 8 (3), axial symmetry of  
$M$ has as axis one
of the above two axes. Moreover, $D(\varepsilon )$ has the
same axis of symmetry (in case of both axes). Therefore 
$M \cap D(\varepsilon )$ also has the same axis of symmetry,
i.e., (5) is satisfied.

In case of Theorem 8 (4), in an analogous way, any non-trivial
congruence admitted by
$M$ is a non-trivial congruence admitted by
$D(\varepsilon )$ as well, hence also
$M \cap D(\varepsilon )$ admits this non-trivial congruence,
i.e., 
(6) is satisfied.
$\blacksquare $
\enddemo


{\it{Proof of Theorem}} 8, {\bf{continuation.}}
{\bf{2.}}
Observe that if for $K$ and $L$
we have supporting circles, of radius some
$r$ (with $r < \pi /2$ for $S^2$), then we have supporting
circles of any larger radius (also
less than $\pi /2$ for $S^2$) as well.
Now we define $R$ as $(r + \pi /2)/2 < \pi /2$
for $S^2$, and as $2r$ for
${\Bbb{R}}^2$ and $H^2$. If for $S^2$ we have that
$r < \pi /2$
is arbitrarily close to $\pi $, or for ${\Bbb{R}}^2$ and
$H^2$ we have that $r$ is arbitrarily large, then the same
statements hold for $R$ rather than $r$.
Occasionally we will use some other $R > r$.
Then {\it{we suppose}} (1)-(4) {\it{in
the proof of Lemma}} 7.1, 
{\it{meant with such an $R$ as
written above}}. 
In the proof of the implication
$(6) \Longrightarrow (7)$ of this theorem we will use only
sets $M$
with diameter $2R$ (which for $S^2$ can be arbitrarily
close to $\pi $, and

\newpage

\noindent
for ${\Bbb{R}}^2$ and $H^2$ can be
arbitrarily
large), and having unique diametral segments,
as will be clear from the proof. 
Hence the hypothesis Theorem 7, (3)
(A) will be satisfied. 
Also, by (3) of the proof of Lemma 7.1, the hypothesis
Theorem 7, (3) (B) will be satisfied.


\proclaim{Lemma 8.2}
Assume \thetag{****} with $d = 2$.
Assume the hypotheses of Theorem {\rm{8}}. Then
for $X = {\Bbb{R}}^2$ and $X = H^2$, {\rm{(6)}} of Theorem
{\rm{8}} 
implies that $K$ and $L$ are $C^1$.
\endproclaim


\demo{Proof}
{\bf{1.}}
Let us suppose that, e.g., $K$ is not smooth, and {\it{for
$x \in
{\text{\rm{bd}}}\,K$ it has an inner angle $\alpha \in
(0, \pi )$}}. 
Since ${\text{bd}}\,L$ is almost everywhere smooth, {\it{we
choose $y \in {\text{\rm{bd}}}\,L$ as a smooth point of
of ${\text{\rm{bd}}}\,L$}}.
{\it{Now let us consider, for this $x$ and this $y$,
the situation described in {\rm{(1)}}
to {\rm{(4)}} in the proof of Lemma}} 7.1. 
We suppose that the centre $o$ of $B(o,R)$ is mapped to
the point $0$
in the collinear model. Further we suppose, as in the proof of
Lemma 7.3, that $(\varphi x)' = R'e_2$ and $(\psi y)' =
-R'e_2$.
We will investigate conv\,$[ (\varphi K) \cup
(\psi L) ] $ in a sufficiently small
neighbourhood of $\varphi x$.
Since $D(\varepsilon ) \cap B(o,R)$ is a fixed neighbourhood of
$D = [\varphi x, \psi y]$
relative to $B(o,R)$, therefore for any fixed $\varepsilon $ 
a sufficiently small neighbourhood of $\varphi x$ relative to
$B(o,R)$ is
contained in $D(\varepsilon )$. We are going to ensure that
conv\,$[ (\varphi K) \cup (\psi L) ] $
is not even locally axially symmetric at $\varphi x$
w.r.t.\ the straight line
spanned by $D$. This will then imply that  
$\left[ {\text{conv}}\,[ (\varphi K) \cup (\psi L) ] \right]
\cap D(\varepsilon )$ is not 
axially symmetric w.r.t.\ this line.

{\bf{2.}}
Since each congruence admitted by
conv\,$[ (\varphi K) \cup (\psi L) ] $
preserves the endpoints of its unique diametral segment
$D$ (whose midpoint is $o$),
therefore the only possible non-trivial congruences admitted
by
conv\,$[ (\varphi K) \cup 
(\psi L) ] $ are central symmetry
w.r.t.\ $o$ and axial symmetry
w.r.t.\ the orthogonal bisector of $D$ (these
interchange $\varphi x $ and $\psi y$), and axial symmetry
w.r.t.\ the line spanned by $D$. We will consider $\varphi x$
and $\psi y$ as the north and south pole of $B(o,R)$ (thus
the images in the collinear model of
$\varphi x$ and $\psi y$ lie on the positive
and negative vertical
coordinate axis, resp.).

{\bf{3.}}
We are going to show that we can avoid axial symmetry of
conv\,$[ (\varphi K) \cup (\psi L) ] $
w.r.t.\ the line spanned by $D$, by suitably choosing
$\varphi $ and $\psi $. Even, this line will not be locally, at
$\varphi x$, an axis of symmetry of
conv\,$[ (\varphi K) \cup (\psi L) ] $.
For this aim we will have to
choose a sufficiently large $R$, greater than the fixed
radius $r$,
in dependence of $K$ and $L$.
In the situation described in {\rm{(1)}}
to {\rm{(4)}} in the proof of Lemma 7.1, we have that
$\varphi K$ and $\psi L$ are contained in {\it{balls
$B_{\varphi x}$ and
$B_{\psi y}$, say, of radius $r$,
both contained in $B(o,R)$, and containing $\varphi x$ and
$\psi y$, resp.}} Then the angle subtended by $\psi L$
at $\varphi x$ is at most the angle subtended by $B_{\psi y}$
at $\varphi x$, which is $\beta := 2 \arcsin [ r/(2R - r)
] $ or $\beta := 2 \arcsin [ {\text{sinh}}\,r/
{\text{sinh}}\,(2R - r) ] $ (by the law of sines),
for ${\Bbb{R}}^2$ or $H^2$,
resp. {\it{We choose $R > r$
so large, that $\beta < \alpha $
should be satisfied.}} Even actually both legs of the
angle subtended by $\psi L$ at $\varphi x$
enclose an angle at most $\beta /2 < \alpha /2 < \pi /2$
with the unit vector $-e_2$ (meant in the collinear model).

Our aim is to show that $\varphi K$ can be rotated through
a small angle about $\varphi x$, while preserving $\psi L$,
so that
conv\,$[ (\varphi K) \cup (\psi L) ] $ 
becomes at $\varphi x$
not even locally axially symmetric w.r.t.\ the straight line
spanned by $D$. 

We will consider the angular domains
subtended by $\varphi K$ and $\psi L$ at $\varphi x$.
These  
will be considered as subsets of the lower semicircle of
$S^1$, which we identify as usual with the angular interval
$[- \pi , 0]$. Let these angular intervals for 
$\varphi K$ and $\psi L$ be $[u_1, u_2]$ and $[v_1, v_2]$,
resp. Here $- \pi \le u_1 < u_2 \le 0$, but $- \pi =
u_1$ and $u_2 = 0$ cannot hold simultaneously, by $\alpha \in
(0, \pi )$. Further,

\newpage

\noindent
$$
- \pi < - \pi /2 - \beta /2 \le v_1 < - \pi /2 <  v_2 \le
- \pi /2 + \beta /2  < 0.
\tag 8.2.1
$$
Here the first and sixth 
inequalities hold by $\beta < \alpha < \pi $.
The second and fifth
inequalities hold since both legs of the
angle, subtended by $\psi L$ at $\varphi x$, enclose an angle
at most $\beta /2$ with $-e_2$ (in the collinear model).
The third and fourth
inequalities hold by smoothness of ${\text{bd}}\,(\psi L)$
at $\psi y$, which implies $D \cap {\text{int}}\,(\psi L) \ne
\emptyset $.
Moreover,
$\pi > \angle (u_1,u_2) = \alpha > \beta \ge \angle (v_1,v_2)
> 0$.
Then the collinear model shows that the
angular domain subtended at $\varphi x$ by  
conv\,$[ (\varphi K) \cup (\psi L) ] $
is
$$
A := [\min \{ u_1, v_1 \} , \max \{ u_2,v_2 \} ].
\tag 8.2.2
$$
(The subtended angular domain clearly contains $A$.
Conversely, there are straight
lines in $X$ through $\varphi x$, of 
direction vectors corresponding to
$\min \{ u_1, v_1 \} , \max \{ u_2,$
\newline
$v_2 \} $,
such that the following holds.
The images of these lines in the collinear model have 
the images both of $\varphi K$ and
$\psi L$ in the collinear model in the lower
side of these image non-vertical lines.)
Since $- \pi \le u_1 < u_2 \le 0$, but $- \pi = u_1$ and
$u_2 = 0$ cannot hold simultaneously, and
$- \pi < v_1 < - \pi /2 < v_2 < 0$,
therefore the angle of the angular
domain $A$ lies in $(0, \pi )$. That is, {\it{$\varphi x$ is
not a smooth point of}} conv\,$[ (\varphi K)
\cup (\psi L) ] $.

Then local symmetry, at $\varphi x$, of
conv\,$[ (\varphi K) \cup (\psi L) ] $,
w.r.t.\ the straight line spanned by $D$,
would imply
symmetry of the angular domain
$A$ w.r.t.\ the same line
({\it{i.e., as a subset of $[- \pi , 0]$,
w.r.t.\ $- \pi /2$}}).
We will investigate the possibilities about the relative
position of the angular intervals $[u_1,u_2]$ and $[v_1,v_2]$,
and in each case we will find a small rotation of $\varphi K$
about $\varphi x$, destroying the symmetry of $A$ 
w.r.t.\ the straight line spanned by $D$.

We may
assume that $u_2 < 0$, since else $\min \{ u_1,v_1 \} > - \pi
$, and $\max \{ u_2, v_2 \} = 0$, and then $A$ is not
symmetric w.r.t.\ the straight line spanned by $D$.
Analogously we  may assume that $u_1 > - \pi $.
So in the
following we assume $- \pi < u_1 < u_2 < 0$.

By $\angle (u_1,u_2) = \alpha > \beta \ge \angle (v_1,v_2)$ we
have either
$u_2 > v_2$, or $u_1 < v_1$. These two cases are symmetrical,
so it will suffice to deal with the case $u_2 > v_2$.
Then we distinguish the following cases: 
\newline
(A) $u_1 < v_1$,
\newline
(B) $ v_1 \le u_1$.

{\bf{4.}}
In Case (A) we have $A = [u_1,u_2] \subset (\pi, 0)$.
Then
we will show that one can rotate $\varphi K$ about
$\varphi x$ through any sufficiently small angle (thus
preserving the equality $A = [u_1,u_2]$), and axial
symmetry w.r.t.\ the straight line spanned by $D$ occurs only
for one position, namely
when the straight line spanned by
$D$ is the angle-bisector of the angular domain $A$. So
generically $A$ will not be axially symmetric w.r.t.\ the
straight line spanned by $D$.
By this small rotation we preserve the inclusion
$\varphi K \subset \left(
{\text{int}}\,B(o,R) \right) \cup \{ \varphi x \} $.

We consider the chords of $B_{\varphi x}$
from $\varphi x$, in the directions of $u_1$ and $u_2$, and
consider the intersection $C$
of $B_{\varphi x}$ with the angular domain
with vertex at $\varphi x$ and legs of directions $u_1,u_2$.
Then $\varphi K \subset C$ (by the collinear model).
Now let us choose some $r^* \in (r,R)$, and let 
$B'_{\varphi x}$ be a circle of radius $r^*$, containing
$\varphi x$ in its

\newpage

\noindent
boundary, and let $B'_{\varphi x} \subset
\left( {\text{int}}\,B(o,R) \right) \cup \{ \varphi x \} $.
Then all boundary points of $C$, lying on
${\text{bd}}\,B_{\varphi x}$, except $\varphi x$,
originally lied in ${\text{int}}\,B'_{\varphi x}$, and
had a distance from
${\text{bd}}\,B'_{\varphi x}$
at least some positive number. So 
for sufficiently small rotations about $\varphi x$ they
still will remain in ${\text{int}}\,B'_{\varphi x}$.
Therefore the rotated copy of $C$, being
the convex hull of the rotated 
above mentioned boundary
points of $C$ and of $\varphi x$, lies in
$\left( {\text{int}}\,B'_{\varphi x} \right)
\cup \{ \varphi x \} $. Hence we have the 
situation as described in (1) to (4) in the proof of Lemma 
7.1, however with $r$ replaced by $r^*$.
Generically the straight line spanned by $D$ is not the
angle-bisector of the angular domain $A$. Hence generically
we do not have at $\varphi x$ even
a local symmetry of
${\text{conv}}\,[(\varphi K) \cup (\psi L)]$
w.r.t.\ the straight line spanned by $D$, so
also we do not have its
global symmetry w.r.t.\ the same line. This
ends the proof of Case (A).

{\bf{5.}}
In case (B) we have $A = [v_1, u_2]$. Here
$v_1 \ge - \pi /2 - \beta /2$, and
$u_2 = u_1 + \alpha > v_1 + \beta \ge - \pi /2 + \beta /2$,
hence $(v_1 + u_2)/2 > - \pi /2$. Therefore
the angular domain
$A$, as a subset of $[- \pi , 0]$, is not symmetric w.r.t.\
$- \pi /2$, hence this angular domain
is not symmetric 
w.r.t.\ the straight
line spanned by $D$ (cf.\ {\bf{3}}).
Hence also
${\text{conv}}\,[(\varphi K) \cup (\psi L)]$ is not 
symmetric w.r.t.\ this line.
This ends the proof of Case (B).

{\bf{6.}}
As proved in {\bf{4}} and {\bf{5}}, there are $\varphi $ and
$\psi $, such
that conv\,$[ (\varphi K) \cup (\psi L) ] $
is not axially symmetric w.r.t.\
the straight line spanned by $D$.
Therefore, by {\bf{2}}, the only possible
non-trivial congruences
admitted by this set
are central symmetry
w.r.t.\ the midpoint $o$ of $D$, and axial symmetry
w.r.t.\
the orthogonal bisector of $D$. Both of these symmetries 
interchange $\varphi x $ and $\psi y$.

In {\bf{3}} we have seen that $\varphi x$ is not a
smooth point of 
conv\,$[ (\varphi K) \cup (\psi L) ] $.
On the other hand, by the choice of $y$ in {\bf{1}}
we have that
$\psi y \in {\text{bd}}\,(\psi L)$ is a smooth point of
$\psi L$. Thus the inner angle of $\psi L$ at $\psi y$ is
$\pi $, which is at most the inner angle of 
conv\,$[ (\varphi K) \cup (\psi L) ] $
at $\psi y$, which is at most $\pi $. Hence $\psi y$ is a
smooth point of 
conv\,$[ (\varphi K) \cup (\psi L) ] $.

However, as above shown,
the only possible non-trivial congruences admitted by 
conv\,$[ (\varphi K) \cup (\psi L) ] $
should interchange $\varphi x$ and $\psi y$. This is a
contradiction, proving Lemma 8.2.
$\blacksquare $
\enddemo


{\it{Proof of Theorem}} 8, {\bf{continuation.}}
{\bf{3.}}
We use the sign $'$ as introduced in the proof of Lemma
7.2, for the images in the collinear model, as a subset of
${\Bbb{R}}^2$. We use, e.g.,
$\left( B(o,R) \right)'$ and
$D(\varepsilon )'$ and ${\text{diam}}'( \cdot )$
as in the proof of Lemma 7.2. As in the proof of Lemma 8.2,
{\bf{1}}, we suppose $o' = 0$ and $(\varphi x)' = R'e_2$ and
$(\psi y)' = -R'e_2$.

For any $x \in {\text{bd}}\,K$ 
we have by (2) in the proof of Lemma 7.1 that $\varphi K$
is included not only in the circle $B(o,R)$ of radius $R$,
but even in a circle $B_{\varphi x}\,\, (\subset B(o,R))$,
of radius
$r\,\,(<R)$, Moreover, $B_{\varphi x}$
also contains $\varphi x$ in its boundary, and is
tangent at $\varphi x$ to $B(o,R)$. The analogous statement
holds also for $\psi L$.

Recall that both $K$ and $L$ are $C^1$, for $S^2$ by
hypothesis of Theorem 8, and for ${\Bbb{R}}^2$ and $H^2$ by
Lemma 8.2. They are also strictly convex, by the supporting
sphere hypothesis of Theorem 8. Hence the spherical maps
${\text{bd}}\,(\varphi K)' \to S^1$ and
${\text{bd}}\,(\psi L)' \to S^1$ are homeomorphisms. 

Combining Lemmas 7.3 and 7.4 (for $d = 2$), we get that

\newpage

\noindent
$$
\cases
{\text{in the situation described in {\rm{(1)}} to {\rm{(4)}}
of the proof of Lemma 7.1,}}
\\
{\text{in a neighbourhood of }} \varphi x, {\text{ or
of }} \psi y, {\text{ we have that }} {\text{\rm{bd}}}\,
[M \cap
\\
D(\varepsilon )] {\text{ coincides with }} {\text{\rm{bd}}}
\, (\varphi K), {\text{ or with }} {\text{\rm{bd}}} \,
(\psi L) , {\text{ resp.}}
\endcases
\tag 8.1
$$

By Theorem 8, (6),
$M \cap D(\varepsilon )$ admits
some non-trivial congruence. Such
a congruence also
preserves the unique diametral segment $D$ of 
$M \cap D(\varepsilon )$, hence also its midpoint $o$.
Thus
$$
\cases
{\text{a non-trivial congruence admitted by }} M \cap
D(\varepsilon ) {\text{ is either a}}
\\
{\text{central symmetry w.r.t.\ }} o, {\text{ or an axial
symmetry w.r.t.\ the straight}}
\\
{\text{line spanned by }} D {\text{ or w.r.t.\ the orthogonal
bisector line of }} D .
\endcases
\tag 8.2
$$

We say that {\it{$x \in {\text{\rm{bd}}}\, K$ is a
point of local symmetry of \,${\text{\rm{bd}}}\, K$}} if the
following holds. {\it{The point $x$
has some open arc neighbourhood relative to
${\text{\rm{bd}}}\, K$,
which is symmetric w.r.t.\ the normal
of \,${\text{\rm{bd}}}\, K$ at $x$.}}
Analogously we define points of local
symmetry of ${\text{bd}}\, L$.

We make a case distinction. Either
\newline
(1)
each $x \in {\text{bd}}\, K$, and each $y \in {\text{bd}}\, L$
is a point of local symmetry of ${\text{bd}}\, K$, and of
${\text{bd}}\, L$, resp., or
\newline
(2)
some $x_0 \in {\text{bd}}\, K$, or 
some $y_0 \in {\text{bd}}\, L$, is not a point of local
symmetry of ${\text{bd}}\, K$, or of
${\text{bd}}\, L$, resp.


In the second
statement of Lemma 8.3
compactness of $K$ and $L$ will not be required.


\proclaim{Lemma 8.3}
Assume \thetag{****} with $d = 2$.
Assume the hypotheses of Theorem {\rm{8}}, and {\rm{(6)}}
of Theorem {\rm{8}}. Let $K$ and $L$ be $C^1$.
Then
case {\rm{(1)}}
in {\bf{3}} of the proof of Theorem {\rm{8}}
implies {\rm{(7)}} of Theorem {\rm{8}}. More
exactly, suppose that $K,L$ are $C^1$, and
satisfy \thetag{*} -- except possibly
${\text{\rm{int}}}\,[(\varphi K) \cap
(\psi L)] \ne \emptyset $. Then case {\rm{(1)}} in {\bf{3}}
of the proof of Theorem {\rm{8}} is equivalent to
that each
connected component of the boundaries of $K$ and $L$ is a
cycle or a straight line. 
\endproclaim


\demo{Proof}
{\bf{1.}}
If each connected component of the boundaries of $K$ and $L$
is a cycle or a straight line, then case {\rm{(1)}} in
{\bf{3}} of the proof of Theorem {\rm{8}} holds evidently.

{\bf{2.}}
Conversely, suppose case {\rm{(1)}} in
{\bf{3}} of the proof of Theorem {\rm{8}}. Then
consider any connected component $K_i$ of
${\text{bd}}\, K$, with $K$ satisfying \thetag{****}. We will
show that it is a cycle, or a straight
line.
(The proof of this statement for $L$ is analogous.
Actually we need this statement for $K$
only in the compact case, when
${\text{bd}}\, K$ is connected.)

We claim that
$K_i$ is symmetric
w.r.t.\ the normal of $K_i$, at any $x \in 
K_i$. In fact, suppose the contrary. Then let
${\widehat{x_1 x_2}} \ne K_i$ be a maximal open
counterclockwise arc
of $K_i$, symmetric w.r.t.\ its normal at $x$.

Possibly
$x_1 = x_2$. However, this can occur only for $K$ compact.
Then $K_i = {\text{bd}}\, K$ equals the closure of
${\widehat{x_1 x_2}}$, hence it is symmetric w.r.t.\ the
normal of $K_i$ at $x$, as asserted. Hence further we may
assume $x_1 \ne x_2$.

The normals of $K_i$ at $x_1$ and $x_2$ are the one-sided
normals of the closure of the

\newpage

\noindent
open arc
${\widehat{x_1 x_2}} \subset K_i$
at $x_1$ and $x_2$.
Therefore these normals
are symmetric images of each
other w.r.t.\ the normal of $K_i$
at $x$. Observe that both $x_1, x_2$
are points
of local symmetry of $K_i$. Hence they are midpoints of
sufficiently short
counterclockwise open arcs ${\widehat{x_{1,1} x_{1,2}}}$ and 
${\widehat{x_{2,1} x_{2,2}}}$ of $K_i$, of equal lengths, and
symmetric w.r.t.\ the normals of
$K_i$ at $x_1$ and $x_2$, resp.
Then
$({\widehat{x_{1,1} x_{1,2}}}) \cup ({\widehat{x_1 x_2}})
\cup 
({\widehat{x_{2,1} x_{2,2}}})$ is an open arc of $K_i$, of
midpoint $x$, symmetric w.r.t.\ the normal of $K_i$
at $x$, strictly containing the maximal such open arc.
This 
is a contradiction,
proving symmetry of $K_i$ w.r.t.\ its normal at any
of its points $x \in K_i$.

We have that $K_i$
is twice differentiable at some of its points $x_0$. Let $x \in
K_i$ be arbitrary. Then let $x^*$ be the midpoint of (one of)
the arc(s) ${\widehat{x_0x}} \subset K_i$. Then the image of
$x_0$ w.r.t.\ the normal of $K_i$ at $x^*$ is $x$. Hence $K_i$
is twice differentiable at $x$, and the curvatures of 
$K_i$ at $x_0$ and $x$ coincide. Hence $K_i$ is of constant
curvature. That is, it is either a cycle, or a straight line.

Then by the supporting circle hypothesis
${\text{bd}}\,K = K_i$ is a circle, so $K$ is a circle.
$\blacksquare $
\enddemo


\proclaim{Lemma 8.4}
Assume \thetag{****} with $d = 2$.
Assume the hypotheses of Theorem {\rm{8}}, and {\rm{(6)}}
of Theorem {\rm{8}}. Let $K$ and $L$ be $C^1$.
Then
case {\rm{(2)}} in {\bf{3}} of the proof of Theorem {\rm{8}}
leads to a contradiction.
\endproclaim


\demo{Proof}
Suppose, e.g., that some point $x_0 \in {\text{bd}}\, K$ is
not a point of 
local symmetry of ${\text{bd}}\, K$. (The proof for $L$ is
analogous.)
Then consider for this
$x_0$ and any $y \in {\text{bd}}\, L$ the
situation described in (1) to (4) in the proof of Lemma 7.1.
By \thetag{8.1}
${\text{bd}}\, (\varphi K)$ and ${\text{bd}}\,[M \cap
D(\varepsilon )]$ at $\varphi x_0$ locally coincide. Therefore
$\varphi x_0$ is not a point of local symmetry of
${\text{bd}}\,[M \cap D(\varepsilon )]$ either.

Then ${\text{conv}}\,[M \cap D(\varepsilon )]$ is not
symmetric w.r.t.\ the
line spanned by $D$. Hence by \thetag{8.2}
it is symmetric either w.r.t.\ $o$, or
w.r.t.\ the
orthogonal bisector of $D$. In both cases the symmetric image
of $\varphi x_0$ is $\psi y$. Moreover, in both cases,
the symmetric image
of a short
open arc of ${\text{bd}}\,[M \cap D(\varepsilon )]$,
with
midpoint $\varphi x_0$, is a short open arc of
${\text{bd}}\,[M \cap D(\varepsilon )]$,
with midpoint $\psi y$. Then, using \thetag{8.1},
if ${\text{bd}}\, (\varphi K)$ at $\varphi x_0$ is not twice
differentiable, then by this symmetry ${\text{bd}}\, (\psi L)$
at any $\psi y \in  
{\text{bd}}\, (\psi L)$ is not twice differentiable either.
This is a
contradiction. Therefore ${\text{bd}}\, (\varphi K)$ at
$\varphi x_0$ is twice
differentiable, and, by this symmetry,
${\text{bd}}\, (\psi L)$ is everywhere
twice differentiable.
Moreover, by this symmetry,
the curvature of ${\text{bd}}\, (\varphi K)$ at $\varphi x_0$
coincides with the curvature of ${\text{bd}}\, (\psi L)$ at
any
$\psi y \in {\text{bd}}\, (\psi L)$. Therefore $\psi L$ has
constant curvature. Then, by the supporting circle
hypothesis, $L$ is a circle.

Then the fixed $\varphi x_0 \in {\text{bd}}\, [M \cap
D(\varepsilon )]$ and any $\psi y \in {\text{bd}}\, [M \cap
D(\varepsilon )]$ can become images of each other either
w.r.t.\, $o$, or w.r.t.\,
the orthogonal bisector of $D$. Moreover, sufficiently short
open arcs of ${\text{bd}}\, [M \cap D(\varepsilon )]$, of
midpoints $\varphi x_0$ and $\psi y$, become then 
images of each other w.r.t.\ one of these symmetries.
Then, on one hand, 
$\varphi x_0 \in {\text{bd}}\, [M \cap
D(\varepsilon )]$ is not a point of local symmetry of  
$M \cap D(\varepsilon )$. On the other hand, $\psi L$ and
$M \cap D(\varepsilon )$ coincide in a neighbourhood of 
$\psi y$. Hence, $\psi y$ is a point of local symmetry of  
$M \cap D(\varepsilon )$. This contradicts the fact that 
$\varphi x_0$ and $\psi y$, as well as the above short
open arcs of ${\text{bd}}\, [M \cap D(\varepsilon )]$,
can be symmetric images of each
other by such a symmetry of $M \cap D(\varepsilon )$.
Thus we have obtained a contradiction.
$\blacksquare $
\enddemo


{\it{Proof of Theorem}} {\rm{8}}, {\bf{continuation.}}
{\bf{4.}}
Now the proof of Theorem 8 follows 

\newpage

\noindent
from the earlier parts of
the proof of Theorem 8, and from Lemmas 8.1-8.4.
We only have to take
into consideration that the supporting circle hypothesis
excludes the paracycle and the hypercycle boundary
components of $K$ and $L$, so the connected boundary
components of $K$ and $L$ are circles. Then, by the last
paragraph of the proof of Lemma 1.7,
$K$ and $L$ are circles.
$\blacksquare $
\enddemo


\definition{Acknowledgements} 
The authors express their gratitude to I. B\'ar\'any, for
carrying the problem, and bringing the two authors
together; to B. Csik\'os for pointing out that Lemma 8.3 is
valid for $K_i \subset X$ being only a $C^1$ $1$-manifold,
which is
also a closed subset. Namely if it can be translated in itself
(this follows from symmetry w.r.t.\ the normal at any of its
points) then it is a cycle or a straight line. 
W
\enddefinition




\Refs

\widestnumber\key{W}


\ref 
\key 
\book 
\by  
\publ 
\publaddr 
\yr 
\endref 

\ref 
\key  
\by 
\paper  
\jour 
\pages  
\endref   

\ref
\key 
\by 
\paper 
\jour 
\vol 
\yr 
\pages  
\endref 



\ref 
\key 1 
\book Theorie der konvexen K\"orper, {\rm{Berichtigter
Reprint (Theory of convex bodies, corrected reprint, in
German)}}
\by T. Bonnesen, W. Fenchel
\publ Springer
\publaddr Berlin-New York
\yr 1974
\endref 



\ref 
\key 2 
\by J. Flachsmeyer
\paper On the convergence of motions
\jour
In: Gen. Top. Rel. Modern Anal. Alg. V (Proc. 5th Prague Top.
Symp., 1981), Sigma
Ser. Pure Math. {\bf{3}}, Heldermann, Berlin, 1983 
\pages 183-188 
\endref   

 
\ref
\key 3 
\by R. High
\paper Characterization of a disc, Solution to problem 1360 (posed by
P. R. Scott)
\jour Math. Magazine
\vol 64 
\yr 1991
\pages 353-354 
\endref 

\ref 
\key 4 
\by J. Jer\'onimo-Castro, E. Makai, Jr.
\paper 
Ball characterizations in spaces of constant curvature
\jour Studia Sci. Math. Hungar.
\vol 55
\yr 2018
\pages 421-478 
\endref 

\ref 
\key 5 
\by J. Jer\'onimo-Castro, E. Makai, Jr.
\paper 
Ball characterizations in planes and spaces of constant
curvature, I
\jour \,submitted 
\vol 
\yr 
\pages  
\endref 





\ref 
\key 6 
\book Convex bodies: the Brunn-Minkowski theory; 
Convex bodies: the Brunn-Minkowski theory, Second expanded
edition, {\rm{Encyclopedia of
Math. and its Appls., Vol.}} {\bf{44}}; {\bf{151}}
\by R. Schneider 
\publ Cambridge Univ. Press
\publaddr Cambridge
\yr 1993; 2014
\endref 





\ref
\key 7 
\paper Beltrami-Klein model
\jour Wikipedia
\endref


\ref
\key 8 
\paper M\"obius transformation, Ch. Subgroups of the M\"obius
group
\jour Wikipedia
\endref


\ref
\key 9 
\paper Support function
\jour Wikipedia
\endref

\endRefs
\enddocument